\documentclass[reqno,a4paper]{siamltex}

\usepackage{graphicx}
\usepackage{color}

\usepackage{amsmath,amsfonts}
\usepackage{mathrsfs}
\usepackage{a4wide}
\usepackage{lscape}
\ifx\theorem\undefined 
\newtheorem{theorem}{Theorem}

\fi

\newcommand{\figref}[1]{Figure~\ref{#1}}

\newcounter{remark}[section]
\def\theremark{\thesection.\arabic{remark}.}
\newenvironment{remark}%
    {\par\medbreak\refstepcounter{remark}%
         {\noindent\bf Remark~\theremark\ }}%
    {\par\medbreak}

\def\RR{\mathbb {R}}

\def\MP{\mathrm{MP}}
\def\wMP{w_\MP^{}}

\def\pref#1{(\text{\ref{#1}})}

\def\Mod#1{\left\|#1\right\|}

\def\R{\mathbb R}
\def\M{\mathscr{M}}
\def\G{\mathscr{G}}

\def\l{\lambda}

\def\Omq{{\Omega_\frac{1}{4}}}
\def\Xq{{X_\frac{1}{4}}}
\def\bA{\mbox{\boldmath $A$}}
\def\bAt{\mbox{\boldmath $\tilde A$}}
\def\bAx{{\bA_x^{}}}
\def\bAxT{{\bA_x^{T}}}
\def\bAy{{\bA_y^{}}}
\def\bAyT{{\bA_y^{T}}}
\def\bAxx{{\bA_{xx}^{}}}
\def\bAyy{{\bA_{yy}^{}}}
\def\bAxy{{\bA_{xy}^{}}}
\def\bAxyT{{\bA_{xy}^{T}}}
\def\bAxyt{{\bAt_{xy}^{}}}
\def\bAbih{{\bA_{\Delta^2}^{}}}
\def\bAbihinv{{\bA_{\Delta^2}^{-1}}}
\def\bAprod{{\bA_{\langle,\rangle}^{}}}
\def\bOne{\mbox{\boldmath $1$}}
\def\bNull{\mbox{\boldmath $0$}}
\def\bC{\mbox{\boldmath $C$}}
\def\bCf{{\bC_{\rm f}^{}}}
\def\bCb{{\bC_{\rm b}^{}}}
\def\bS{\mbox{\boldmath $S$}}
\def\bLambda{\mbox{\boldmath $\Lambda$}}
\def\bLambdaxx{{\bLambda_{xx}^{}}}
\def\bLambdayy{{\bLambda_{yy}^{}}}
\def\bLambdaxy{{\bLambda_{xy}^{}}}
\def\bLambdabih{{\bLambda_{\Delta^2}^{}}}

\def\bB{\mbox{\boldmath $B$}}
\def\bG{\mbox{\boldmath $G$}}
\def\bGt{\mbox{\boldmath $\tilde G$}}
\def\vKD{Von K\'arm\'an-Donnell}%
\def\XXint#1#2#3{{\setbox0=\hbox{$#1{#2#3}{\int}$}
     \vcenter{\hbox{$#2#3$}}\kern-.5\wd0}}

\def\frparbcenter#1{\framebox(43,60)[t]{\parbox{43mm}{\begin{center}#1\end{center}}}}

\author{Ji\v r\'I Hor\'ak\thanks{Universit\"at K\"oln, Germany}
\and Gabriel J. Lord\thanks{Heriot-Watt University, Edinburgh, United Kingdom}
\and Mark A. Peletier\thanks{Technische Universiteit Eindhoven, The Netherlands}}
\title{Numerical variational methods applied to cylinder buckling}
\date{DRAFT \today}

\begin{document}
\maketitle

\begin{abstract}
We review and compare different computational variational methods applied to a 
system of fourth order equations that arises as a model of cylinder buckling. 
We describe both the discretization and implementation, in particular how to 
deal with a 1 dimensional null space. We show that we can construct many 
different solutions from a complex energy surface.
We examine numerically convergence in the spatial discretization and
in the domain size.  
Finally we give a physical interpretation of some of the solutions found.
\end{abstract}

\section{Introduction}
We describe complementary approaches to finding solutions of
systems of fourth order elliptic PDEs. The techniques are applied to a
problem that arises in the classic treatment of an isotropic
cylindrical shell under axial compression but are also applicable to a
wide range of problems such as waves on a suspension bridge
\cite{Ho1,HoMcK}, the Fu\v{c}\'{\i}k spectrum of the Laplacian
\cite{HoRe}, or the formation of microstructure \cite{SH,FM}.

\begin{figure}[htbp]
  \centering
  \setlength{\unitlength}{1mm}
  \begin{picture}(140,78)
    \put(10,45){\includegraphics[width=4.9cm]{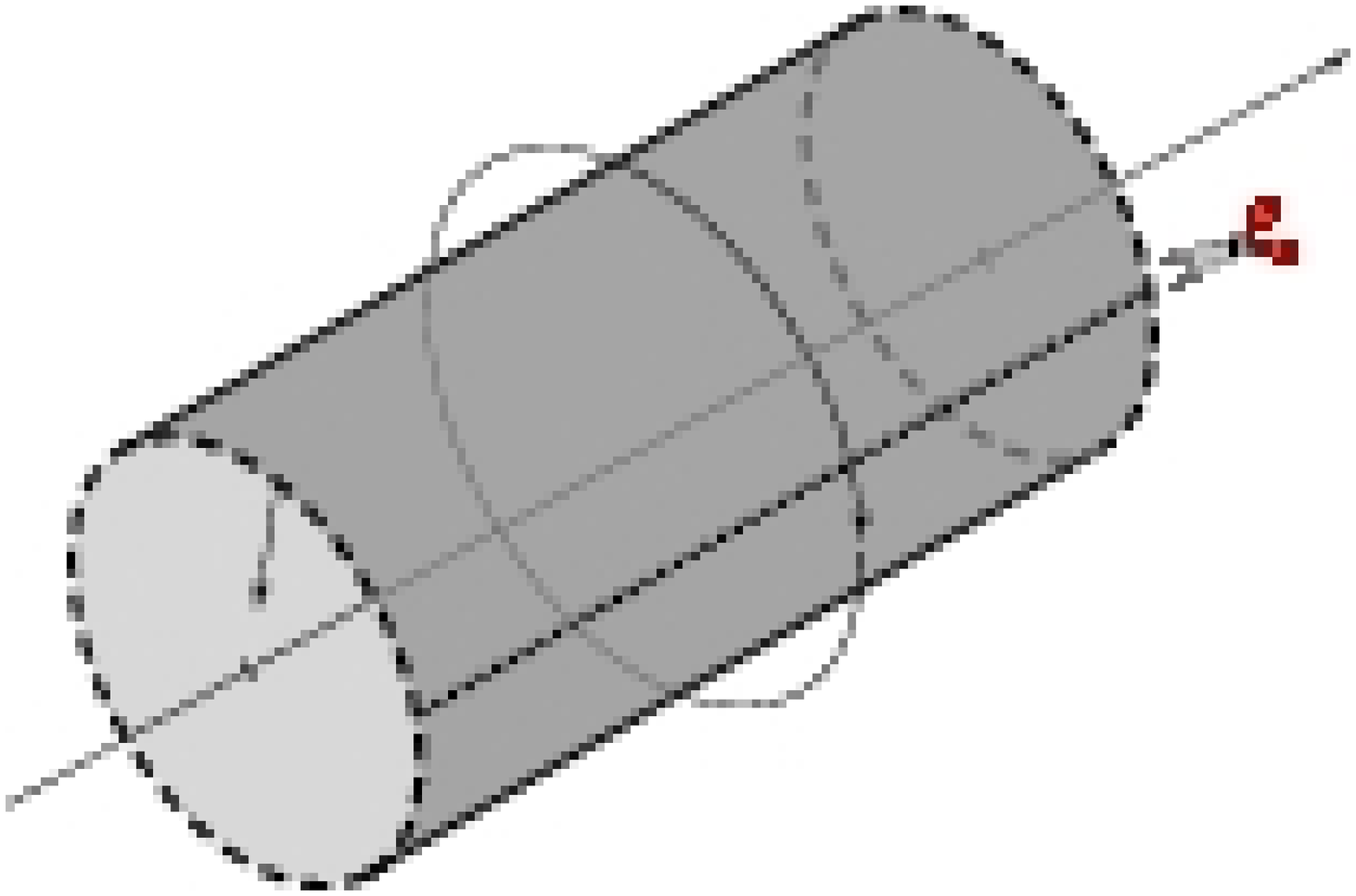}}
    \put(67,60){\includegraphics[width=1.35cm]{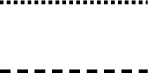}}
    \put(0,0){\includegraphics[width=6.75cm]{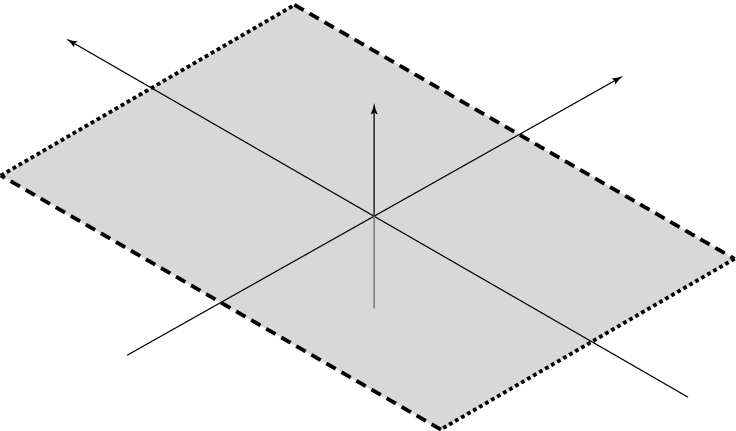}}
    \put(72,0){\includegraphics[width=6.75cm]{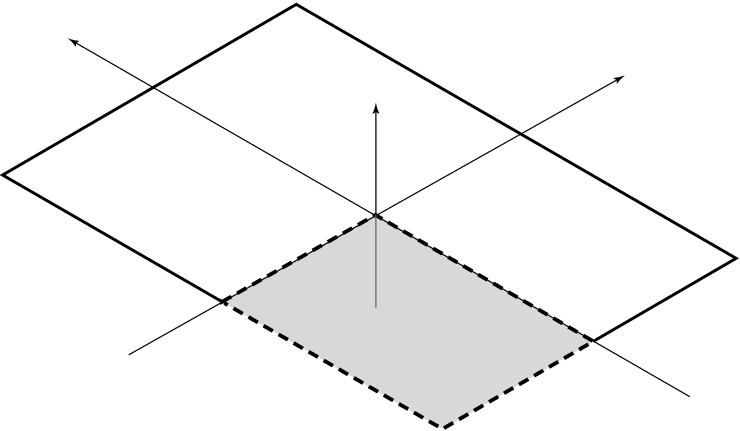}}
    \put(83,66){periodic boundary condition for $w,\phi$}
    \put(83,59.5){$w_\nu=(\Delta w)_\nu=\phi_\nu=(\Delta\phi)_\nu=0$}

    \put(10,72){(i)}
    \put(16,55){\small $w$}
    \put(58,72.5){\small $x$}
    \put(42,58.5){\small $y$}
    \put(16.5,50.5){\small $-a$}
    \put(46,65.5){\small $a$}

    \put(-2,36){(ii)}
    \put(35,29){\small $w$}
    \put(55,29){\small $x$}
    \put(8.5,36){\small $y$}
    \put(16.5,8){\small $-a$}
    \put(46.5,24){\small $a$}
    \put(51,4){\small $-b$}
    \put(12.5,28){\small $b$}
    \put(39,2){\small $\Omega$}

    \put(70,36){(iii)}
    \put(107,29){\small $w$}
    \put(127,29){\small $x$}
    \put(80.5,36){\small $y$}
    \put(88.5,8){\small $-a$}
    \put(123,4){\small $-b$}
    \put(109.7,3){\small $\Omq$}
  \end{picture}
  \caption{(i) The geometry of the cylinder, (ii) the computational domain and the boundary conditions, (iii) one quarter of the domain and the corresponding boundary conditions.}
  \label{fig:geom}
\end{figure}

The cylindrical shell offers a computationally 
challenging and physically relevant problem with a complex energy
surface.  We take as our model for the shell the \vKD{}
equations which can be rescaled \cite{HoLoPe1} to the form
\begin{align}
\label{eq:main_w2}
&  \Delta^2  w +  \l  w_{ x x} 
  - \phi_{ x x} - 2\,[ w,\phi] = 0, \\
&\Delta^2 \phi +  w_{ x x} + [ w, w] = 0,
\label{eq:main_phi2}
\end{align}
where the bracket is defined as
\begin{equation}
[u,v] = \frac12 u_{xx}v_{yy} + \frac 12 u_{yy}v_{xx} - u_{xy}v_{xy}.
\label{eq:bracket}
\end{equation}
The function $w$ is a scaled inward radial displacement measured from
the unbuckled (fundamental) state, $\phi$ is the Airy stress function,
and $\l\in(0,2)$ is a load parameter. The unknowns $w$ and $\phi$ are
defined on a two-dimensional spatial domain
$\Omega=(-a,a)\times(-b,b)$, where $x\in(-a,a)$ is the axial and
$y\in(-b,b)$ is the tangential coordinate. 
Since the $y$-domain $(-b,b)$ represents the circumference of the cylinder,
the following boundary conditions are prescribed:
\begin{subequations}
\label{def:BC}
\begin{align}
&w\text{ is periodic in $y$,\quad  and} \quad  w_x = (\Delta  w)_x = 0
\text{ at }x=\pm a,
\label{def:BCw} \\
&\phi\text{ is periodic in $y$,\quad  and} \quad  \phi_x = (\Delta  \phi)_x = 0
\text{ at }x=\pm a,
\label{def:BCphi}
\end{align}
\end{subequations}
as shown in Fig.~\ref{fig:geom}~(i), (ii).

\subsection{Functional setting}
We search for weak solutions $w$, $\phi$ of
(\ref{eq:main_w2}--\ref{def:BC}) in the space
\[
X=\left\{\psi\in H^2(\Omega) : \psi_x(\pm a,\cdot)=0,\ \psi \text{ is
    periodic in } y, \text{ and }\,\int_\Omega\psi = 0\right\}
\]
with norm
\[
  \Mod{w}^2_X = \int_\Omega\bigl(\,\Delta w^2 +\Delta\phi_1^2\,\bigr),
\]
where $\phi_1\in H^2(\Omega)$ is the unique solution of
\begin{equation}
  \label{eq:phi1}
  \Delta^2\phi_1 = -w_{xx}, \qquad \phi_1 \text{ satisfies~\pref{def:BCphi}, }
\qquad \text{and}\qquad \int_\Omega \phi_1 = 0.  
\end{equation}
This norm is equivalent to the $H^2$-norm on $X$, and with the
appropriate inner product $\langle\cdot,\cdot\rangle_X$ the space $X$
is a Hilbert space. Alternatively, if the load parameter $\l\in(0,2)$
is fixed, another norm
\[
  \Mod{w}^2_{X,\l} = \int_\Omega\bigl(\,\Delta w^2 +\Delta\phi_1^2 
  - \l w_x^2\,\bigr)
\]
can be used. Because of the estimate
\[
\int_\Omega w_x^2 = -\int_\Omega ww_{xx} = \int_\Omega w\Delta^2 \phi_1
= \int_\Omega \Delta w\Delta \phi_1 \leq \frac12 \int_\Omega \Delta w  ^2 + \frac12 \int_\Omega \Delta \phi_1^2 
= \frac12 \Mod{w}_X^2,
\]
it is equivalent to $\Mod{\cdot}_X$ and hence also to the $H^2$-norm
on $X$. The corresponding inner product will be denoted
$\langle\cdot,\cdot\rangle_{X,\l}$.

Equations~(\ref{eq:main_w2}--\ref{eq:main_phi2}) are related to the
stored energy $E$, the average axial shortening $S$, and the total
potential given by
\begin{equation}
\label{def:energy}
E(w) :=  \frac12\int_\Omega \left( \Delta  w^2 +  \Delta \phi^2\right), \qquad
\quad S(w):= \frac12 \int_\Omega  w_{ x}^2, \qquad\quad F_\l=E-\l S.
\end{equation}
Note that the function $\phi$ in~\pref{def:energy} is determined from
$w$ by solving~\pref{eq:main_phi2} with boundary
conditions~\pref{def:BCphi}. All the functionals $E$, $S$, and $F_\l$ belong to
$C^1(X)$, i.e., are continuously Fr\'echet differentiable.

The fact that~\pref{eq:main_w2} is a reformulation of the stationarity
condition $F_\l'=E'-\lambda S'=0$ will be important in
Sec.~\ref{sec:discretization}, and we therefore briefly sketch the
argument. It is easy to see that
\[
S'(w)\cdot h=-\int_\Omega w_{xx}h.
\]
For $E'(w)\cdot h$, let $w,\phi\in X$ solve \pref{eq:main_phi2} and
$h,\psi\in X$ solve $ \Delta^2 \psi = -h_{xx}-[h,h]-2[w,h] $. Then,
assuming sufficient regularity on $w$,
\begin{align}
E(w+h)-E(w) &= \int_\Omega\Delta w\Delta h +
\frac{1}{2}\int_\Omega(\Delta h)^2 + \int_\Omega\Delta\phi\Delta\psi +
\frac{1}{2}\int_\Omega(\Delta\psi)^2 \notag \\
&= \int_\Omega\left(h\Delta^2w - h\phi_{xx} - 2[w,h]\phi\right)
+ \frac{1}{2}\int_\Omega(\Delta h)^2 + \frac{1}{2}\int_\Omega(\Delta\psi)^2
- \int_\Omega[h,h]\phi\ , \notag
\end{align}
where we used integration by parts several times. The last three
integrals are $O(\Mod{h}_X^2)$ for $\Mod{h}_X\to 0$ and it can be
shown by integration by parts that
\begin{equation}
\int_\Omega [w,h]\,\phi = \int_\Omega h\,[w,\phi]\ .
\label{eq:w_h_phi}
\end{equation}
Therefore 
\[
F_\l'(w)\cdot h = E'(w)\cdot h - \lambda S'(w)\cdot h = 
\int_\Omega h\Bigl( \Delta^2  w +  \l  w_{ x x} 
  - \phi_{ x x} - 2\,[ w,\phi]\Bigr).
\]

\subsection{Review of some variational numerical methods}\label{sec:var_num}
We now describe the variational methods used to find numerical
approximations of critical points of the total potential $F_\l$. In
our numerical experiments these methods are accompanied by Newton's
method and continuation. The advantage of this approach is that it
combines the knowledge of global features of the energy landscape with
local ones of a neighborhood of a critical point. The details related
to spatial discretization will be discussed in
Sec.~\ref{sec:discretization}, the Newton-based methods in
Sec.~\ref{sec:Newton}.

\subsubsection{Steepest descent method (SDM)}\label{sec:SDM}
Let the load parameter $\l\in(0,2)$ be fixed; we work in a discretized
version of $(X,\langle\cdot,\cdot\rangle_{X,\l})$. We try to minimize
the total potential $F_\l$ by following its gradient flow. We solve the initial value
problem
\[
\frac{d}{dt}w(t)=-\nabla_\l F_\l(w(t))\,,\qquad w(0)=w_0\,,
\]
with a suitable starting point $w_0$ on some interval $(0,T]$. This
problem is then discretized in $t$.

In~\cite{HoLoPe1} it was shown that $w=0$ is a local minimizer of
$F_\l$. Indeed, if $\Mod{w_0}_{X,\l}$ is small, the numerical solution
$w(t)$ converges to zero as $t$ tends to infinity. If, on the other hand,
$\Mod{w_0}_{X,\l}$ is large, the numerical solution $w(t)$ stays
bounded away from zero. In most of our experiments, the numerical
algorithm did not converge for $t\to\infty$ in the large norm case. The only
exception for a relatively small value of $\l$ will be mentioned later
in Sec.~\ref{sec:domain_size}. Nevertheless, for a sufficiently large
computational domain $\Omega$ and a sufficiently large $t>0$ we obtain
$F_\l(w(t))<0$. Such a state $w(t)$ is needed for the mountain pass
algorithm as explained below. Existence of this state was also proved
in~\cite{HoLoPe1}.

\subsubsection{Mountain-pass algorithm (MPA)}\label{sec:MPA}
The algorithm was first proposed in \cite{ChMcK1} for a second order
elliptic problem in 1D and extended in \cite{HoMcK} to a
fourth-order problem in 2D. We give a brief description of the
algorithm here.

\begin{figure}[htbp]
  \begin{center}
    \setlength{\unitlength}{1mm}
    \begin{picture}(100,67)
      \put(0,0){{\includegraphics[width=10cm]{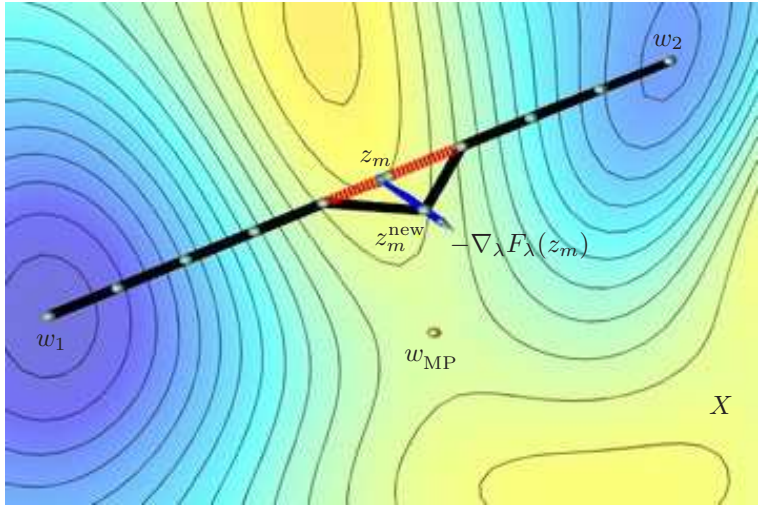}}}
      \put(4,21){$w_1$}
      \put(86,61){$w_2$}
      \put(47,45.5){$z_m$}
      \put(49,34.5){$z_m^\mathrm{new}$}
      \put(59,33.5){$-\nabla_\l F_\l(z_m)$}
      \put(53,19){$w_\MP^{}$}
      \put(93.5,12){$X$}
    \end{picture}
  \end{center}
  \caption{Deforming the path in the main loop of the mountain pass
    algorithm: point $z_m$ is moved a small distance in the direction
    $-\nabla_\l F_\l(z_m)$ and becomes $z_m^\mathrm{new}$. This step
    is repeated until the mountain pass point $\wMP$ is reached.}
  \label{fig:mpa}
\end{figure}

Let the load $\l\in(0,2)$ be fixed; we work again in a
discretized version of $(X,\langle\cdot,\cdot\rangle_{X,\l})$. We
denote $w_1=0$ the local minimum of $F_\l$ and take a point $w_2$ such
that $F_\l(w_2)<0$ (in practice this point is found using the SDM). We
take a discretized path $\{z_m\}_{m=0}^{p}$ connecting $z_0=w_1$ with
$z_p=w_2$. After finding the point $z_m$ at which $F_\l$ is maximal
along the path, this point is moved a small distance in the direction
of the steepest descent $-\nabla_\l F_\l(z_m)$. Thus the path has been
deformed and the maximum of $F_\l$ lowered. This deforming of the path
is repeated until the maximum along the path cannot be lowered any
more: a critical point $\wMP$ has been reached. Figure~\ref{fig:mpa}
illustrates the main idea of the method.

The mountain-pass algorithm is local in its nature. The numerical
solution $\wMP$ it finds has the mountain-pass property in a certain
neighborhood only. The choice of the path endpoint $w_2$ may influence
to which critical point the algorithm converges. Different choices of
$w_2$ are in turn achieved by choosing different initial points $w_0$
in the SDM.

\subsubsection{Constrained steepest descent method (CSDM)}
We fix the amount of shortening $S$ of the cylinder. This is often
considered as what actually occurs in experiments. We work now in a
discretized version of $(X,\langle\cdot,\cdot\rangle_X)$. Let $C>0$ be
a fixed number and define a set of $w$ with constant shortening
\begin{equation}
  \label{eq:constraint}
  \M=\{w\in X\,:\, S(w)=C\}\,.
\end{equation}
Critical points of $E$ under this constraint are critical points of
$F_\l$, where $\l$ is a Lagrange multiplier. The simplest such points
are local minima of the stored energy $E$ on~$\M$. We need to follow the gradient flow
of $E$ on~$\M$, hence we solve the initial value problem
\[
\frac{d}{dt}w(t)=-P_{w(t)}\nabla E(w(t))\,,\qquad w(0)=w_0\in\M\,,
\]
for $t>0$. $P_w$ denotes the orthogonal projection in $X$ on the tangent space of
$\M$ at $w\in\M$:
\[
P_w u= u - \frac{\langle\nabla S(w),u\rangle_X}{\Mod{\nabla S(w)}_X^2} \nabla S(w)\,.
\]
The details of the algorithm can be found in~\cite{Ho1}. The initial
value problem is solved by repeating the following two steps: given a
point $w\in\M$ find $\bar w=w - \Delta t P_w\nabla E(w)$, where
$\Delta t>0$ is small, and define $w_{\rm new}=c\bar w$, where the
scaling coefficient $c$ is chosen so that $w_{\rm new}\in\M$. The
algorithm is stopped when $\Mod{P_w\nabla E(w)}_X$ is smaller than a
prescribed tolerance. The corresponding load is given by
$$\l=\frac{\langle\nabla S(w),\nabla E(w)\rangle_X}{\Mod{\nabla
    S(w)}_X^2}.$$

\subsubsection{Constrained mountain-pass algorithm (CMPA)}
Let $C>0$ and~$\M$ be the set of $w$ with constant shortening given
in~\pref{eq:constraint}. We would like to find mountain-pass points of
$E$ on~$\M$. The method has been described in~\cite{Ho1} in detail. We
need two local minima $w_1, w_2$ of $E$ on~$\M$ which can be found
using the CSDM. The algorithm is then similar to the MPA. We take a
discretized path $\{z_m\}_{m=0}^{p}\subset\M$ connecting $z_0=w_1$
with $z_p=w_2$. After finding the point $z_m$ at which $E$ is maximal
along the path, this point is moved a small distance in the tangent
space to~$\M$ at $z_m$ in the direction of the steepest descent
$-P_{z_m}\nabla E(z_m)$ and than scaled (as in the CSDM) to come back to
$\M$. Thus the path has been deformed on~$\M$ and the maximum of $E$
lowered. This deforming of the path is repeated until the maximum
along the path cannot be lowered any more: a mountain-pass point of
$E$ on~$\M$ has been reached. The load $\l$ is computed as in the CSDM.

The choice of the end points $w_1$ and $w_2$ will in general influence to which
critical point the algorithm converges.

\subsection{Computational Domains}
We consider the problem on the domain $\Omega$
(\figref{fig:geom}~(ii)) both without further restraints and under a
symmetry assumption, which reduces the computational complexity. In
the latter case we assume
\begin{equation}
  \label{eq:sym}
  \begin{array}{c}
    w(x,y)=w(-x,y)=w(x,-y) \\
    \phi(x,y)=\phi(-x,y)=\phi(x,-y)
  \end{array}
  \quad \mbox{for}\ (x,y)\in\Omega\ .
\end{equation}
By looking for solutions $w,\phi\in X$ that satisfy~\pref{eq:sym} the
domain is effectively reduced to one quarter:
$\Omq=(-a,0)\times(-b,0)$ as shown in~\figref{fig:geom}~(iii). One needs to
solve (\ref{eq:main_w2}--\ref{eq:main_phi2}) only on $\Omq$ with the
boundary conditions
\begin{equation}
  \label{def:BCsym}
  w_\nu=(\Delta w)_\nu=\phi_\nu=(\Delta\phi)_\nu=0 \quad \mbox{on}\ \partial\Omq\ ,
\end{equation}
where $\nu$ denotes the outward normal direction to the boundary.
Hence we search for weak solutions of
(\ref{eq:main_w2}--\ref{eq:main_phi2}), \pref{def:BCsym} in the space
\[
\Xq=\left\{\psi\in H^2(\Omq) : \psi_\nu=0 \text{ on } \partial\Omq,
  \text{ and }\,\int_\Omq\psi = 0\right\}\ .
\]
We can then use \pref{eq:sym} to extend these functions to the whole
$\Omega$.

We have performed numerical experiments both with and without the
symmetry assumption. For the sake of simplicity we will give a
detailed description of the numerical methods for the second case only
where the boundary conditions are the same on all sides of $\Omq$. The
first case with periodic conditions on two sides of $\Omega$ is very
similar and will be briefly mentioned in~Remark~\ref{rem:full_domain}

\subsection{Solving the biharmonic equation}
In order to obtain $\phi$ for a given $w$, one has to
solve~\pref{eq:main_phi2}; to compute the norm of $w$, one has to
solve~\pref{eq:phi1}. Both problems are of the form
\begin{equation}
  \label{eq:biharm}
  \Delta^2\psi = f \text{ in } \Omq, \qquad \psi_\nu=(\Delta\psi)_\nu=0
  \text{ on } \partial\Omq, \qquad \int_\Omq\psi = 0,
\end{equation}
where $f\in L^1(\Omq)$ is given. If $\int_\Omq f=0$, then
\pref{eq:biharm} has a unique weak solution $\psi$ in $\Xq$. It is a
straightforward calculation to verify that the right-hand sides of
equations in \pref{eq:main_phi2} and \pref{eq:phi1} have zero average.

In the discretization described below the problem \pref{eq:biharm} is
treated as a system:
\begin{equation}
  \label{eq:biharm_syst}
  \begin{array}{c}
    -\Delta u = f \\ -\Delta v = u
  \end{array}
  \text{ in } \Omq, \qquad u_\nu=v_\nu=0
  \text{ on } \partial\Omq,
  \qquad \int_\Omq u = \int_\Omq v = 0.
\end{equation}
The system has a unique weak solution $(u,v)\in(H^1(\Omq))^2$. Since
the domain $\Omq$ has no reentrant corners, Theorem~1.4.5
of~\cite{MazjaNazPlam} guarantees that $v\in H^2(\Omq)$ and therefore
that the two formulations are equivalent.

\section{Finite difference discretization}\label{sec:discretization}
We discretize the domain $\Omq$ by a uniform mesh $(x_m,y_n)\in\Omq$
with $M$ points in the $x$-direction and $N$ points in the
$y$-direction:
\[
\begin{array}{ccc}
  x_m = -a + (m-\frac{1}{2})\Delta x, & \qquad & m\in\{1\ldots,M\},\\
  y_n = -b + (n-\frac{1}{2})\Delta y, & \qquad & n\in\{1\ldots,N\},
\end{array}
\]
where $\Delta x = a/M$, $\Delta y = b/N$. We represent the values of
some function $w$ on $\Omq$ at these points by a vector
$w=(w_i)_{i=1}^{MN}$, where $w_i=w(x_m,y_n)$ and $i=(n-1)M+m$. In our
notation we will not distinguish between $w$ as a function and $w$ as
a corresponding vector. The vector $w$ can also be interpreted as a
block vector with $N$ blocks, each containing $M$ values of a single
row of the mesh. We introduce the following conventions for notation:
\begin{itemize}
\item For two matrices $A^M=(a_{ij})_{i,j=1}^M$,
  $B^N=(b_{k\ell})_{k,\ell=1}^N$ we define $A^M \otimes B^N := (
  b_{k\ell}A^M)_{k,\ell=1}^N$, which is an $N\times N$ block matrix,
  each block is an $M\times M$ matrix.
\item For two vectors $u=(u_i)_{i=1}^{MN}$, $v=(v_i)_{i=1}^{MN}$ we
  define $u\odot v = (u_i v_i)_{i=1}^{MN}$, i.e., a product of
  the components.
\end{itemize}

To discretize second derivatives we use the standard central finite
differences (with Neumann boundary conditions \cite{RchtmyrMrtn}).
Let $Id^M$ denote the $M\times M$ identity matrix and define another
$M\times M$ matrix
\[
A_2^M=\left[
  \begin{array}{rrrrr}
    1 & -1 \\
    -1 & 2 & -1 \\
    & \ddots & \ddots & \ddots \\
    & & -1 & 2 & -1 \\
    & & & -1 & 1
  \end{array}
\right].
\]
The second derivatives $-\partial_{xx}$, $-\partial_{yy}$ and the
biharmonic operator $\Delta^2$ with the appropriate boundary
conditions are approximated by
\[
  \bAxx = \frac{1}{\Delta x^2}A_2^M\otimes Id^N, \qquad
  \bAyy = \frac{1}{\Delta y^2}Id^M\otimes A_2^N, \qquad
  \bAbih = (\bA_{xx} + \bA_{yy})^2,
\]
respectively.

\subsection{Discretization of $E$, $S$, and the bracket
  $[\cdot,\cdot]$}\label{sec:ES_disc}
Supposing that we can solve the discretized version of~\pref{eq:main_phi2}
\begin{equation}
  \label{eq:main_phi_d}
  \bAbih\phi-\bAxx w + [w,w]_2^{} = \bNull\,,
\end{equation}
we can also evaluate the energy $E$ and the shortening $S$:
\begin{equation}
E(w)=2\left(w^T \bAbih\, w + \phi^T \bAbih\,
  \phi\right)\Delta x\,\Delta y, \qquad S(w)=2\left(w^T \bAxx
  w\right)\Delta x\,\Delta y.
\label{eq:ES_disc}
\end{equation}

In order to solve~\pref{eq:main_phi_d} we need to be able to solve the
biharmonic equation and we need to choose a discretization of the
bracket $[\cdot,\cdot]$. This bracket appears in the equations in two
different roles: in equation~\pref{eq:main_phi2} the bracket is part
of the mapping $w\mapsto \phi$, and therefore of the definition of the
energy $E$; in equation~\pref{eq:main_w2}, which represents the
stationarity condition $E'-\lambda S'=0$, the bracket appears as a
result of differentiating $E$ with respect to $w$ and applying
partial integration. As a result, we need to use two different forms
of discretization for the two cases.

In both cases the bracket requires a discretization of the mixed derivative
$\partial_{xy}$. One choice is to use one-sided finite
differences. Define $M\times M$ matrices
\begin{equation}\label{eq:dx_RL}
A_{1L}^M=\left[
  \begin{array}{rrrr}
    0 \\
    -1 & 1 \\
    & \ddots & \ddots \\
    & & -1 & 1
  \end{array}
\right],\qquad
A_{1R}^M=\left[
  \begin{array}{rrrr}
    -1 & 1 \\
    & \ddots & \ddots \\
    & & -1 & 1 \\
    & & & 0
  \end{array}
\right].
\end{equation}
We choose either left or right-sided differences represented by these
matrices, respectively, let $A_1^M$ denote our choice
(cf.~Sec.~\ref{sec:bias}). The derivatives $\partial_x$, $\partial_y$,
and $-\partial_{xy}$ are approximated by
\begin{equation}\textstyle\label{eq:d_xy_RL}\textstyle
\bAx=\frac{1}{\Delta x}\,A_1^M\otimes Id^N, \qquad
\bAy=\frac{1}{\Delta y}\,Id^M\otimes A_1^N, \qquad
\bAxy=-\bAx\bAy\,.
\end{equation}

For the definition of $\phi$ in terms of $w$
(equation~\pref{eq:main_phi2}) we choose
\begin{equation}
  \label{eq:bracket_phi}
  [w,w]_2^{} = (\bAxx w)\odot(\bAyy w) - (\bAxy w)\odot(\bAxy w)\, ,
\end{equation}
and the corresponding choice for equation~\pref{eq:main_w2} is 
\begin{equation}
  \label{eq:bracket_w}
  \textstyle
  [w,\phi]_1^{}=\frac{1}{2}\bAyy\left\{(\bAxx w)\odot\phi\right\} +
  \frac{1}{2}\bAxx\left\{(\bAyy w)\odot\phi\right\} -
  \bAxyT\left\{(\bAxy w)\odot\phi\right\}\,.
\end{equation}
These two definitions are related in the sense given
in~\pref{eq:w_h_phi}: $[w,h]_2^T\phi= h^T[w,\phi]_1^{}$ for all $h$.

With these definitions the partial derivatives of discretized $E$
and $S$ with respect to the components of $w$ are given by
\begin{equation}
  \label{eq:fr_der}
  E'(w)=  \bAbih w +\bAxx\phi -2[w,\phi]_1^{}\,,\qquad
  S'(w)=  \bAxx w\,.
\end{equation}

\subsection{Solving the discretized biharmonic equation}
Matrix $\bAbih$ is symmetric and has a one-dimensional nullspace:
$\bAbih\bOne = \bNull$, where $\bOne$ and $\bNull$ are vectors with
$MN$-components which are all one and all zero, respectively. The same is
true for $\bAxx$ and $\bAyy$. For a given vector $f$ we would like to
solve
\begin{equation}
  \label{eq:biharm_disc}
  \bAbih\psi = f\,,\qquad \bOne^T\psi = 0\,.  
\end{equation}
A unique solution exists if and only if $f$ has zero average, i.e.,
$\bOne^T f=0$. So we must verify that the discretized versions of the
right-hand sides in \pref{eq:main_phi2}, \pref{eq:phi1} satisfy this
condition. Let $w\in\RR^{MN}$, then
\begin{align}
  \bOne^T \bAxx w &= 0\,, \notag \\
  \bOne^T [w,w]_2^{} &= (\bAxx w)^T (\bAyy w) - (\bAxy w)^T (\bAxy w) 
  = w^T \bAxx \bAyy w - w^T \bAxyT \bAxy w = 0\,, \label{eq:int_by_parts}
\end{align}
where the last equality holds because $\bAxT\bAx=\bAxx$ and
$\bAyT\bAy=\bAyy$, and because the $x$- and $y$-matrices commute. We
have, in fact, shown that the integration by parts formula from the
continuous case holds for our choice of spatial discretization. This is 
not true for an arbitrary discretization but is key for a successful scheme.

The inverse of matrix $\bAbih$ on the subspace of vectors with zero
average, denoted with a slight abuse of notation by $\bAbihinv$, can
be found, for example, using the fast cosine transform described below
in Sec.~\ref{sec:Fourier}.

\subsection{Computing the gradient}
The variational methods of this paper are based on a steepest descent
flow and modifications of this algorithm. The direction of the
steepest descent of $E$ at a point $w\in X$ is opposite to the
gradient of $E$ at $w$. The gradient is the Riesz representative of
the Fr\'echet derivative and hence we need to find a vector $u\in X$,
such that $E'(w)\cdot v=\langle u,v\rangle$ for all $v\in X$. The
inner product is either $\langle\cdot,\cdot\rangle_X$ or
$\langle\cdot,\cdot\rangle_{X,\l}$ and hence the gradient depends on
the choice of the inner product. We use the notation $u=\nabla E(w)$
for the gradient in $(X,\langle\cdot,\cdot\rangle_X)$ and $u=\nabla_\l
E(w)$ for the gradient in $(X,\langle\cdot,\cdot\rangle_{X,\l})$. To
find the discretized version of the gradient, we first need to
discretize the inner product.

Let $u,v\in\R^{MN}$, $\bOne^T u=\bOne^T v =0$. The inner product is
evaluated in the following way:
\begin{align}
  \langle u,v\rangle_{X,\l} &=4 \left(u^T \bAbih v + {\phi_1^u}^T \bAbih
    \phi_1^v -\l u^T \bAxx v\right) \Delta x\,\Delta y \notag \\
  &=4\Bigl( u^T\left(\bAbih +
      \bAxx\bAbihinv\bAxx -\l\bAxx\right)v\Bigr)\Delta x\,\Delta y \,, \notag
\end{align}
where $\phi_1^u, \phi_1^v$ are solutions of the discretized version
of~\pref{eq:phi1} with $w$ replaced by $u$ and $v$, respectively, and
we assume that we work on $\Omq$. For $w\in\R^{MN}$, $\bOne^T w=0$ the
Riesz representative of $E'(w)$ given in \pref{eq:fr_der} is computed
as
\begin{equation}\label{eq:grad}
  \nabla_\l E(w)=\left(\bAbih + \bAxx\bAbihinv\bAxx -
    \l\bAxx \right)^{-1} E'(w)\,.
\end{equation}
As in the case of $\bAbihinv$ we abused notation here since the
inverse only makes sense on a subspace of vectors with zero average.
It can be easily verified that $\bOne^T E'(w)=0$. The numerical
evaluation of $\nabla_\l S$ and of the $\langle
\cdot,\cdot\rangle_X$-gradients is similar.

\subsection{Fourier coordinates}\label{sec:Fourier}
In Fourier coordinates most of the finite difference operators become
diagonal matrices. This increases the efficiency of the numerical
algorithm and makes it possible to easily find the inverse of matrices
like $\bAbih$. See for example \cite{Canuto}.

On a uniform grid, it is a standard procedure to apply some form of
the fast Fourier transform to diagonalize finite difference matrices
like $A_2^M$ (see, e.g., \cite{Strang}). Due to the Neumann boundary
conditions \pref{def:BCsym} we need to employ the fast cosine
transform. We define $M\times M$ matrices
\[\textstyle
C_{\rm f}^M=\frac{1}{\sqrt{2M}}\left(2\cos\frac{(i-1)(2j-1)\pi}{2M}\right)_{i,j=1}^M , \qquad
C_{\rm b}^M=\frac{1}{\sqrt{2M}}\left[\left.
\begin{array}{c} 1 \\ \vdots \\1 \end{array}\right|
\left(2\cos\frac{(2i-1)(j-1)\pi}{2M}\right)_{i=1,j=2}^{M,M} \right],
\]
which have the following properties:
\[\textstyle
C_{\rm f}^M C_{\rm b}^M = Id^M, \qquad
C_{\rm f}^M A_2^M C_{\rm b}^M = \Lambda^M,
\]
where $\Lambda^M=\diag(2-2\cos\frac{(m-1)\pi}{M})_{m=1}^M$. Hence they
are inverses of each other and diagonalize $A_2^M$.

We further define matrices
\[
\bCf = C_{\rm f}^M\otimes C_{\rm f}^N, \qquad
\bCb = C_{\rm b}^M\otimes C_{\rm b}^N,
\]
which diagonalize $\bAxx$, $\bAyy$, and $\bAbih$:
\begin{equation}\label{eq:diag1}
\bCf\bAxx\bCb = \bLambdaxx, \qquad
\bCf\bAyy\bCb = \bLambdayy, \qquad
\bCf\bAbih\bCb = \bLambdabih,
\end{equation}
where the diagonal matrices are given by
\begin{equation}\label{eq:diag2}\textstyle
\bLambdaxx=\frac{1}{\Delta x^2}\Lambda^M\otimes Id^N, \qquad
\bLambdayy=\frac{1}{\Delta y^2}Id^M\otimes \Lambda^N, \qquad
\bLambdabih=(\bLambdaxx + \bLambdayy)^2.
\end{equation}
For a vector $w\in\R^{MN}$ we introduce its Fourier coordinates $\hat
w$ by
\[
\hat w = \bCf w \qquad\qquad w = \bCb \hat w \ .
\]
We note that $\bOne^T w=0$ if and only if the first component of $\hat
w$ is zero.

Most of the computations involved in the variational methods described
in~Sec.~\ref{sec:var_num} can be done in the Fourier coordinates. The
only time one needs to go back to the original coordinates is when
evaluating the brackets \pref{eq:bracket_phi} and \pref{eq:bracket_w},
because they are nonlinear and involve the discretized mixed
derivative operator $\bAxy$.

\subsection{Alternative discretization of $-\partial_{xy}$}\label{sec:d_xy}
The fast Fourier transform provides us with another discretization of
the mixed derivative which is not biased to the left or right. In an
analogy to~\pref{eq:diag2} and \pref{eq:diag1} we define
\[\textstyle
\bLambdaxy=\frac{1}{\Delta x\,\Delta y}\sqrt{\Lambda^M}\otimes
\sqrt{\Lambda^N}, \qquad \bAxyt=\bS \bLambdaxy \bCf\ ,
\]
where $\bS$ is the fast sine transform matrix
\[\textstyle
\bS = S^M\otimes S^N\,,\quad
S^M=\frac{1}{\sqrt{2M}}\left[\left.
\begin{array}{c} 0 \\ \vdots \\0 \end{array}\right|
\left(2\sin\frac{(2i-1)(j-1)\pi}{2M}\right)_{i=1,j=2}^{M,M} \right]\,.
\]
Property \pref{eq:int_by_parts} also holds with $\bAxy$ replaced by
$\bAxyt$.

\begin{remark}\label{rem:full_domain}
  When discretizing the problem on the full domain $\Omega$ with
  boundary conditions~\pref{def:BC}, we need to use different matrices
  in the $x$ and $y$-directions. In the $x$-direction we use the
  matrices described above, in the $y$-direction to discretize the
  second derivatives, for example, we use
  \[
  A_2=\left[
    \begin{array}{rrrrr}
      2 & -1 & & & -1 \\
      -1 & 2 & -1 \\
      & \ddots & \ddots & \ddots \\
      & & -1 & 2 & -1 \\
      -1 & & & -1 & 2
    \end{array}
  \right].
  \]
  In this direction the fast Fourier transform is used instead of the
  fast cosine/sine transform.
\end{remark}

\section{Newton's method}\label{sec:Newton}
We use Newton's method in two different ways. The first is to improve the numerical approximations
obtained by the variational numerical methods. Since these
are sometimes slow to converge, it is often faster to stop such an
algorithm early and use its result as an initial guess for  Newton's
method. The second use for Newton's method is as part of a numerical
continuation algorithm (see Sec. \ref{sec:continuation}).

\subsection{Newton's method for given load parameter $\l$}
This method can be used to improve solutions obtained by the MPA. Let
$\l\in(0,2)$ be given. We are solving
\begin{equation}
  \label{eq:main_newton_p}
  \G(w,\phi)=\left[
  \begin{array}{l}
    \G_1 \\ \G_2\rule{0pt}{2em}
  \end{array}
  \right]=
  \left[
    \begin{array}{l}
      \bAbih w- \l\bAxx w+ \bAxx\phi -2[w,\phi]_1 \\
      -\bAbih\phi + \bAxx w - [w,w]_2\rule{0pt}{2em}
    \end{array}
  \right]=
  \left[
    \begin{array}{l}
      \bNull \\ \bNull\rule{0pt}{2em}
    \end{array}
  \right]
\end{equation}
for $w$ and $\phi$ with zero average using Newton's method. The
matrix we need to invert is
\begin{equation}
  \label{eq:newton_matrix}
  \G'(w,\phi)=\left[
    \begin{array}{ll}
      \frac{\partial\G_1}{\partial w} &
      \frac{\partial\G_1}{\partial \phi} \\
      \frac{\partial\G_2}{\partial w} &
      \frac{\partial\G_2}{\partial \phi}\rule{0pt}{2em}
    \end{array}
  \right]=
  \left[
    \begin{array}{ll}
      \bAbih-\l\bAxx -2\bB_1 \quad & \bAxx -2\bB_2 \\
      \bAxx -2\bB_2^T & -\bAbih \rule{0pt}{2em}
    \end{array}
  \right]\,,
\end{equation}
where
\begin{align}
  \bB_1 = \frac{\partial}{\partial w}[w,\phi]_1 &=
  \frac{1}{2}\bAxx(\diag\phi)\bAyy + \frac{1}{2}\bAyy(\diag\phi)\bAxx -
  \bAxyT(\diag\phi)\bAxy\,, \notag \\
  \bB_2 = \frac{\partial}{\partial\phi}[w,\phi]_1 &=
  \frac{1}{2}\left(\frac{\partial}{\partial w}[w,w]_2\right)^T \notag \\
  &= \frac{1}{2}\bAxx(\diag\bAyy w) + \frac{1}{2}\bAyy(\diag\bAxx w) -
  \bAxyT(\diag\bAxy w)\,. \notag
\end{align}
From the properties of $\bAxx, \bAyy, \bAxy, \bAbih$ it follows that
matrix $\G'(w,\phi)$ is symmetric and singular, its nullspace is
spanned by $ \left[
  \begin{array}{c}
    \bOne\\ \bNull
  \end{array}
\right]
$, $
\left[
  \begin{array}{c}
    \bNull \\ \bOne
  \end{array}
\right]
$.

To describe how to find the inverse of $\G'(w,\phi)$ on a subspace
orthogonal to its nullspace, we introduce a new notation for the four
blocks of $\G'(w,\phi)$ from~\pref{eq:newton_matrix}:
\[
\G'(w,\phi)=\left[
  \begin{array}{cc}
    \bG_{11}^{} & \bG_{12}^{} \\ \bG_{12}^T & \bG_{22}^{}
  \end{array}
\right]
\]
For given vectors $u$, $\eta$ with zero average we need to find
vectors $v$, $\zeta$ with zero average such that
\begin{equation}\label{eq:newt_lin_syst}
\left[
  \begin{array}{cc}
    \bG_{11}^{} & \bG_{12}^{} \\ \bG_{12}^T & \bG_{22}^{}
  \end{array}
\right]
\left[
  \begin{array}{c}
    v \\ \zeta
  \end{array}
\right]=
\left[
  \begin{array}{c}
    u \\ \eta
  \end{array}
\right]\,.
\end{equation}
Let the tilde denote the block of the first $MN-1$ rows and columns of a
matrix and $MN-1$ components of a vector. The matrix $\left[
  \begin{array}{cc}
    \bGt_{11}^{} & \bGt_{12}^{} \\ \bGt_{12}^T & \bGt_{22}^{}
  \end{array}
\right]$ is symmetric, nonsingular, and sparse. It can be inverted by
a linear sparse solver. System~\pref{eq:newt_lin_syst} is then
solved in the following steps:
\[
\left[
  \begin{array}{c}
    r \\ \rho
  \end{array}
\right] := 
\left[
  \begin{array}{cc}
    \bGt_{11}^{} & \bGt_{12}^{} \\ \bGt_{12}^T & \bGt_{22}^{}
  \end{array}
\right]^{-1}
\left[
  \begin{array}{c}
    \tilde u \\ \tilde \eta
  \end{array}
\right]\,,\qquad
\begin{array}{ll}
  s :=-\frac{1}{MN}\tilde\bOne^T r\,, \qquad &
  v=\left[
    \begin{array}{c}
      r + s\tilde\bOne \\ s
    \end{array}
  \right]\,, \\
  \sigma := -\frac{1}{MN}\tilde\bOne^T \rho\,, &
  \eta=\left[
    \begin{array}{c}
      \rho + \sigma\tilde\bOne \\ \sigma
    \end{array}
  \right]\,.\rule{0pt}{2em}
\end{array}
\]

\subsection{Newton's method for given $S$}
This method can be used to improve solutions obtained by the CSDM and
the CMPA. Let $C>0$ be given. We are looking for numerical solutions
of~(\ref{eq:main_w2}--\ref{eq:main_phi2}) in the set~$\M$ defined
by~\pref{eq:constraint}. Hence we are solving
\begin{equation}
  \label{eq:main_newton_C}
  \left[
    \begin{array}{l}
      \bAbih w- \l\bAxx w+ \bAxx\phi -2[w,\phi]_1 \\
      -\bAbih\phi + \bAxx w - [w,w]_2\rule{0pt}{2em} \\
      -\frac{1}{2}w^T\bAxx w + C/(4\Delta x\Delta y)\rule{0pt}{2em}
    \end{array}
  \right]=
  \left[
    \begin{array}{c}
      \bNull \\ \bNull\rule{0pt}{2em} \\ 0\rule{0pt}{2em}
    \end{array}
  \right]
\end{equation}
for $w$ and $\phi$ with zero average and $\l$ using Newton's method.
The approach is very similar to that described in the previous
section, the resulting matrix is symmetric, has just one more row and
column than the matrix in~\pref{eq:newton_matrix}.

\subsection{Continuation}\label{sec:continuation}
To follow branches of solutions $(\l,w)$ of~\pref{eq:main_newton_p} we
adopt a continuation method described in~\cite{Keller}. We introduce a
parameter $s\in\R$ by adding a constraint---pseudo-arclength
normalization (in the $(\l,\Mod{w}_X)$-plane). For a given value of
$s$ we are solving
\begin{equation}
  \label{eq:main_cont}
  \bar\G(w,\phi,\l)=\left[
  \begin{array}{l}
    \G_1 \\ \G_2\rule{0pt}{2em} \\ \G_3\rule{0pt}{2em}
  \end{array}
  \right]=
  \left[
    \begin{array}{l}
      \bAbih w- \l\bAxx w+ \bAxx\phi -2[w,\phi]_1 \\
      -\bAbih\phi + \bAxx w - [w,w]_2\rule{0pt}{2em} \\
      \theta \langle \dot w_0,w-w_0\rangle_X + (1-\theta)
      \dot\l_0 (\l-\l_0) - (s-s_0)\rule{0pt}{2em}
    \end{array}
  \right]=
  \left[
    \begin{array}{c}
      \bNull \\ \bNull\rule{0pt}{2em} \\ 0\rule{0pt}{2em}
    \end{array}
  \right]
\end{equation}
for $w$, $\phi$ with zero average, and the load $\l$, where the value
of $\theta\in(0,1)$ is fixed (e.g., $\theta=\frac{1}{2}$). We assume
that we are given a value $s_0$, an initial point $(\l_0,w_0)$ on the
branch, and an approximate direction $(\dot\l_0,\dot w_0)$ of the
branch at this point (an approximation of the derivative
$\frac{d}{ds}(\l(s),w(s))|_{s=s_0}$).

System \pref{eq:main_cont} is solved for a discrete set of values of
$s$ in some interval $(s_0,s_1)$ by Newton's method. Then, a new
initial point on the branch is defined by setting $w_0=w(s_1)$,
$\l_0=\l(s_1)$, $s_0=s_1$, a new direction $(\dot\l_0,\dot w_0)$ at
this point is computed and the process is repeated. The matrix we need
to invert in Newton's method is
\begin{equation}
  \label{eq:cont_matrix}
  \bar \G'(w,\phi,\l)=
  \left[
    \begin{array}{cc}
      \G'(w,\phi,\l) & g \\
      h^T & d 
    \end{array}
  \right],\quad
  g=\left[
    \begin{array}{c}
      -\bAxx w \\ \bNull
    \end{array}
  \right],\ 
  h=\left[
    \begin{array}{c}
      4\theta \bAprod \dot w_0\,\Delta x\Delta y \\ \bNull
    \end{array}
  \right],\ 
  d= (1-\theta)\dot\l_0,
\end{equation}
where $\bAprod=\bAbih+\bAxx\bAbihinv\bAxx$.

Solving a linear system with this matrix amounts to solving system
\pref{eq:newt_lin_syst} for two right-hand sides. For a given
$u\in\R^{2MN}$ with $[\bOne^T\ \bOne^T]u=0$ and a given $\eta\in\R$ we
want to solve
\begin{equation}\label{eq:cont_lin_syst}
\left[
  \begin{array}{cc}
    \G'(w,\phi,\l) & g \\ h^T & d
  \end{array}
\right]
\left[
  \begin{array}{c}
    v \\ \zeta
  \end{array}
\right]=
\left[
  \begin{array}{c}
    u \\ \eta
  \end{array}
\right]\,.
\end{equation}
for $v$ with $[\bOne^T\ \bOne^T]v=0$ and $\zeta$. System
\pref{eq:cont_lin_syst} is solved in the following steps:
\[
\begin{array}{ll}
  {\rm solve}: & \G'(w,\phi,\l)v_1=g\,, \\
  & \G'(w,\phi,\l)v_2=u\,,
\end{array}\qquad
\zeta = \frac{\eta-h^Tv_2}{d-h^Tv_1}\,,\qquad
v=v_2-\zeta v_1\,.
\]

\begin{remark}
Note that in this implementation we simply follow a solution of the
equation, there is no guarantee that this remains a local minimum, a
MP-solution or constrained MP-solution (cf.~Fig.~\ref{fig:cont_100}).
\end{remark}

\begin{remark}
  Newton's method and continuation have been implemented only using
  a one-sided discretization of the mixed derivative $\partial_{xy}$
  and only on the domain $\Omq$ assuming symmetry~\pref{eq:sym}.
  The alternative discretization of $\partial_{xy}$ described in
  Sec.~\ref{sec:d_xy} uses the fast cosine/sine transform. The
  resulting matrix $\bAxyt$ is not sparse and therefore we would
  obtain a dense block $\bG_{12}^{}$ in system~\pref{eq:newt_lin_syst}
  which would prevent us from using a sparse solver.
  
  On the full domain $\Omega$ we assume periodicity of $w$ and $\phi$
  in the $y$-direction. Hence for a discretization with a small step
  $\Delta y$ the matrix we invert when solving \pref{eq:newt_lin_syst}
  would become close to singular. The shift in the $y$ direction is
  prevented by assuming the symmetry $w(x,y)=w(x,-y)$,
  $\phi(x,y)=\phi(x,-y)$.
\end{remark}
  
\section{Numerical solutions}
We fix the size of the domain and the size of the space step for the
following numerical computations: $a=b=100$, $\Delta x=\Delta y=0.5$.
We obtain solutions using the variational techniques SDM, MPA, CSDM,
and CMPA.
  Table~\ref{tab:methods} provides a summary of which discretization
  was used in which algorithm.

\begin{table}[htbp]
  \centering
  \renewcommand{\arraystretch}{1.2}
  \begin{tabular}{l|p{15mm}|p{15mm}|p{15mm}|p{15mm}|}
    \cline{2-5}
    & \multicolumn{2}{|c|}{mixed derivative $\partial_{xy}$} &
    \multicolumn{2}{|c|}{computational domain} \tabularnewline
    \cline{2-5}
    & \centering one-sided & \centering Fourier & \centering full &
    \centering $1/4$ \tabularnewline
    \hline
    \multicolumn{1}{|l|}{variational methods} &
    \centering $\surd$ & \centering $\surd$ & \centering $\surd$ &
    \centering $\surd$ \tabularnewline
    \hline
    \multicolumn{1}{|l|}{Newton/continuation} &
    \centering $\surd$ & & & \centering $\surd$ \tabularnewline
    \hline
  \end{tabular}
  \vspace{3mm}
  \caption{Summary of which spatial discretization was used in the
    different numerical methods.}
  \label{tab:methods}
\end{table}

\subsection{A mountain-pass solution on the full domain $\Omega$}\label{sec:num_mp}

\begin{figure}[htbp]
\begin{center}
\color[rgb]{.5,.5,.5}
\fbox{\includegraphics[height=35mm]{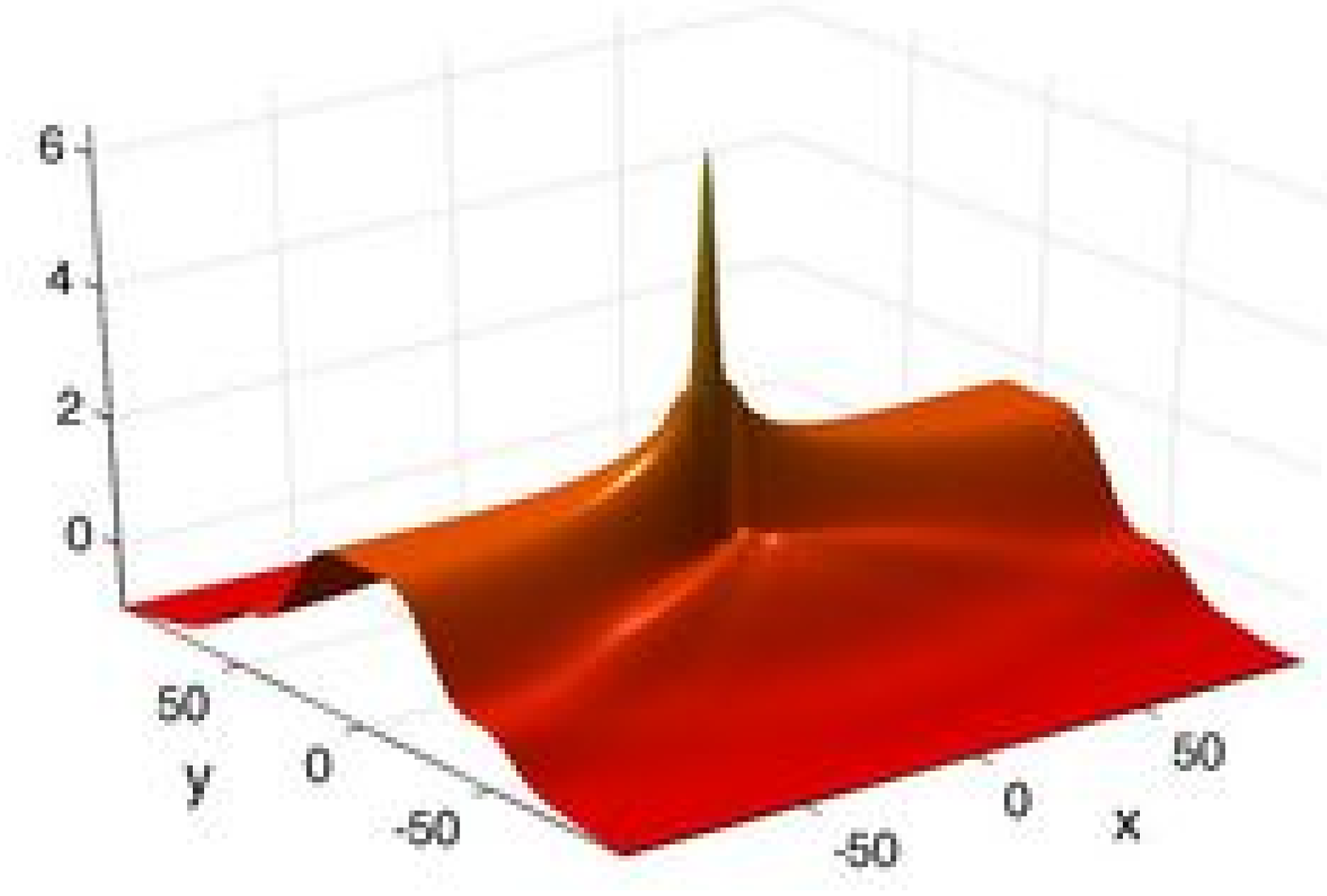}
  \includegraphics[height=35mm]{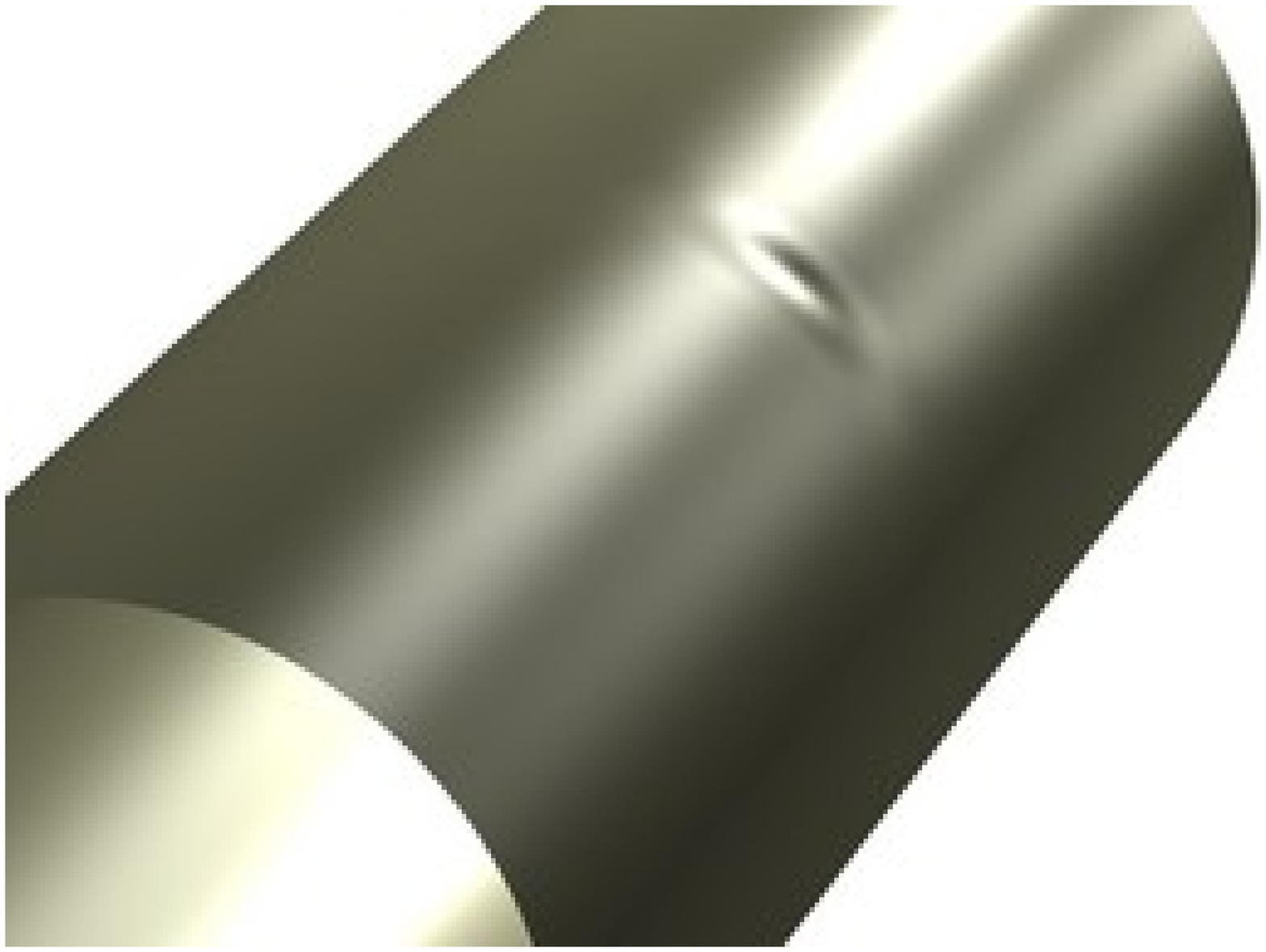}}
\color{black}
\end{center}
\caption{Mountain-pass solution for $\l=1.4$ found using the MPA on the full
  domain $\Omega$ with $\partial_{xy}$ discretized using the fast
  Fourier transform.}
\label{fig:num_wMP}
\end{figure}

The first numerical experiments are done on the full domain $\Omega$, without symmetry restrictions, and with the unbiased (Fourier) discretization of $\partial_{xy}$ (Sec.~\ref{sec:d_xy}). For a fixed load $\l=1.4$ we computed a mountain-pass solution using 
the MPA (Sec.~\ref{sec:MPA}). As end points were taken $w_1=0$ and a second point $w_2$ obtained by the SDM (here the initial point for the SDM was chosen to have a single peak centered
at $x=y=0)$. The graph of the
solution  $\wMP$ is shown in Fig.~\ref{fig:num_wMP}
(left). The figure on the right shows $\wMP$ rendered on a cylinder 
and we see it has the form of a single dimple. The value of shortening 
for this solution is $S(\wMP)=14.93529$.

Alternatively, if we apply the CSDM with $S=14.93529$ and use a function
with a single peak in the center of the domain as the initial
condition $w_0$ we also obtain the same solution $\wMP$, this time as 
a local minimizer of $E$ under constrained $S$.

We remark that although the MPA and the CSDM have a local character we have 
not found any numerical mountain-pass solution with the total
potential $F_\l$ smaller than $F_\l(\wMP)$ for $\l=1.4$. Similarly,
using the CSDM we have not found any solution with energy $E$ smaller than
$E(\wMP)$ under the constraint $S=14.93529$.
We briefly discuss the physical relevance of this solution
in~Sec.~\ref{sec:conc}.

\subsection{Solutions under symmetry restrictions}
The solution $\wMP$ of Fig.~\ref{fig:num_wMP} satisfies the symmetry
property~\pref{eq:sym}. In the computations described below we
enforced this symmetry and worked on the quarter domain $\Omq$, thus
reducing the complexity of the problem. In order to improve the
variational methods by combining them with Newton's method we also
discretized the mixed derivative $\partial_{xy}$ using left-sided
finite differences. The influence of this choice on the numerical
solution is described in~Sec.~\ref{sec:bias}.

\subsubsection{Constrained steepest descent method}
\begin{landscape}
  \begin{figure}[htbp]
    \centering
    \setlength{\unitlength}{1mm}
    \begin{picture}(181,124)
      \color[rgb]{.5,.5,.5}
      \put(0,64){\frparbcenter{\includegraphics[width=4.233cm]{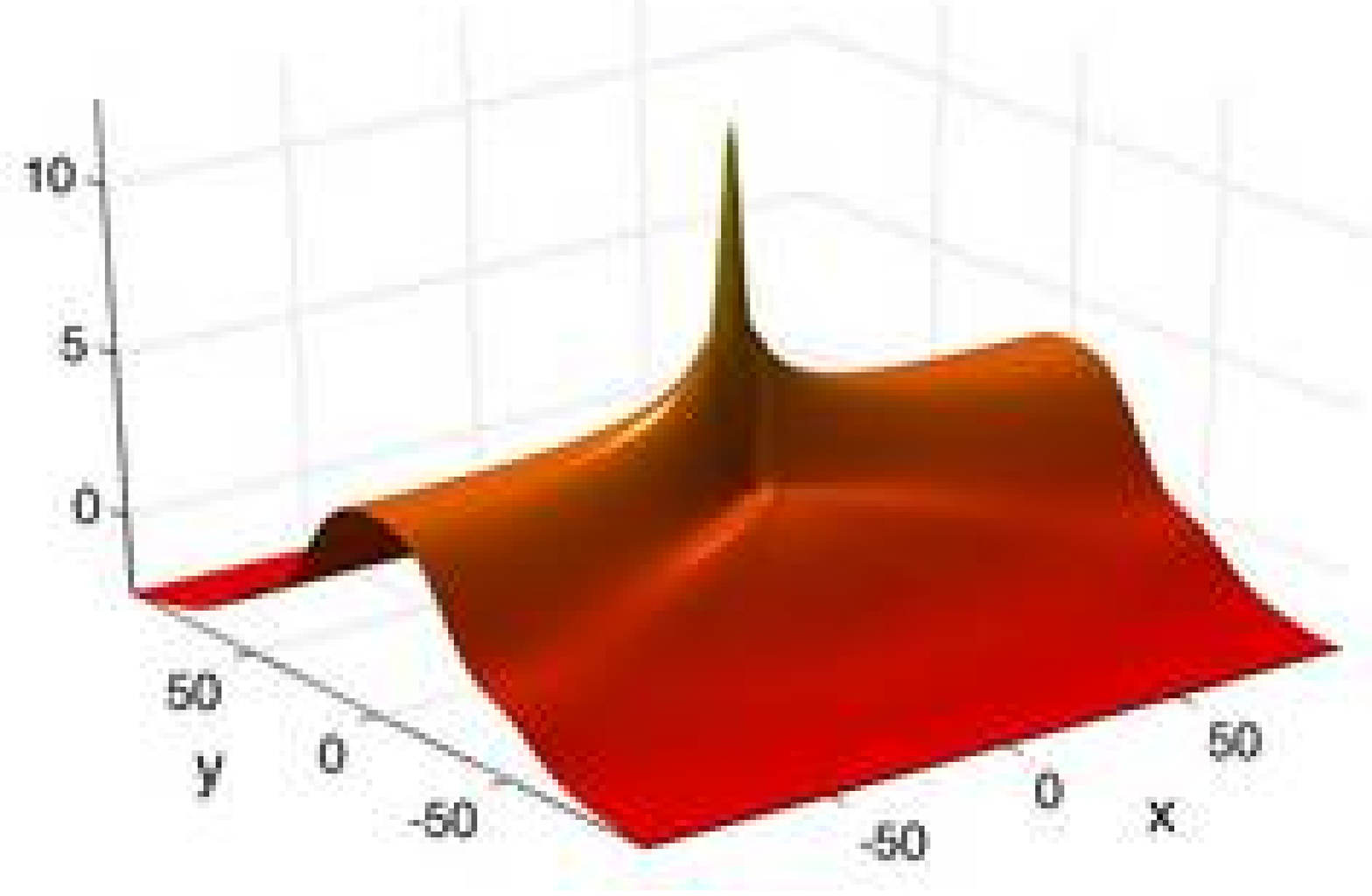}\\{\includegraphics[width=4.064cm]{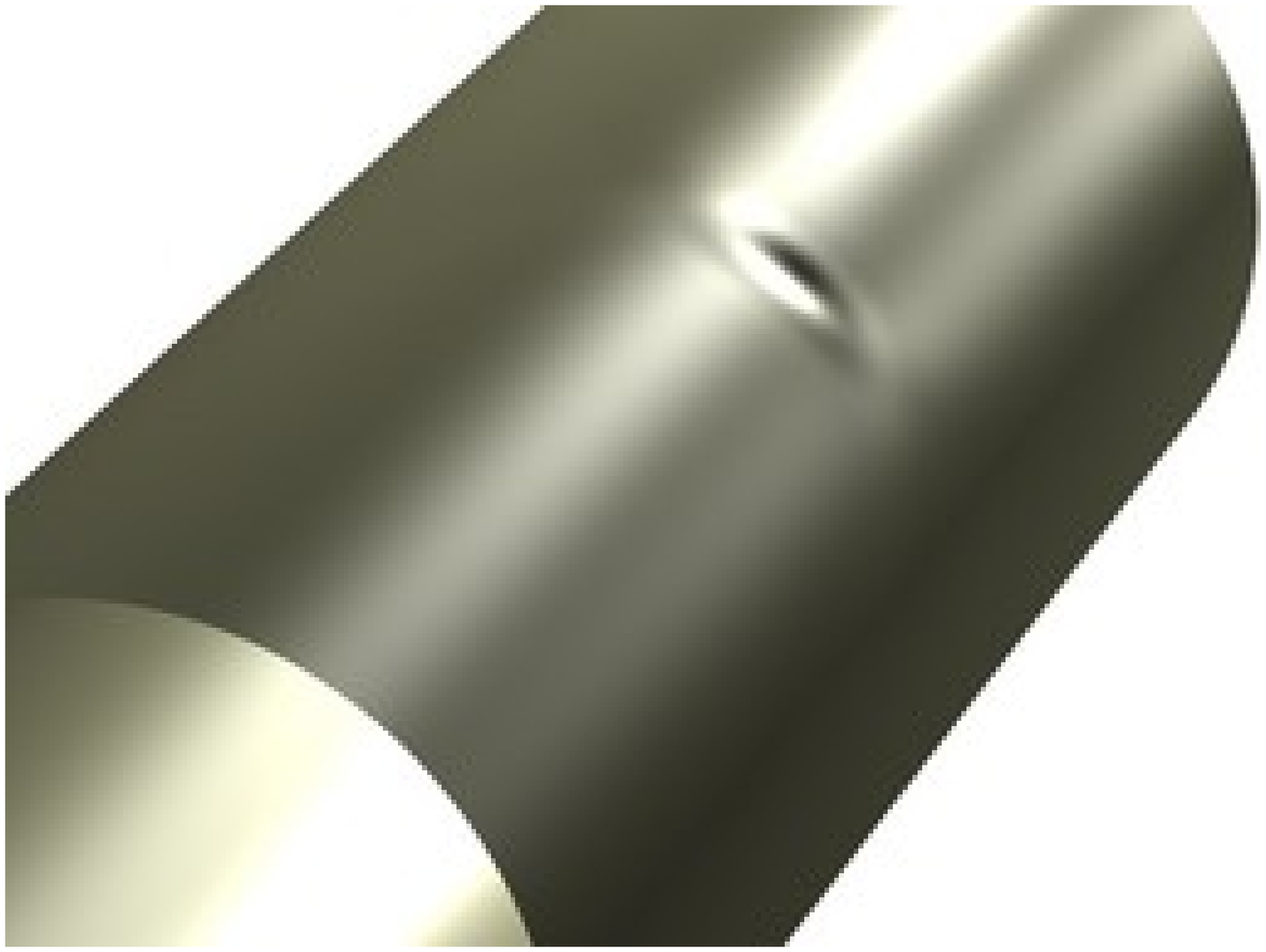}}}}
      \put(46,64){\frparbcenter{\includegraphics[width=4.233cm]{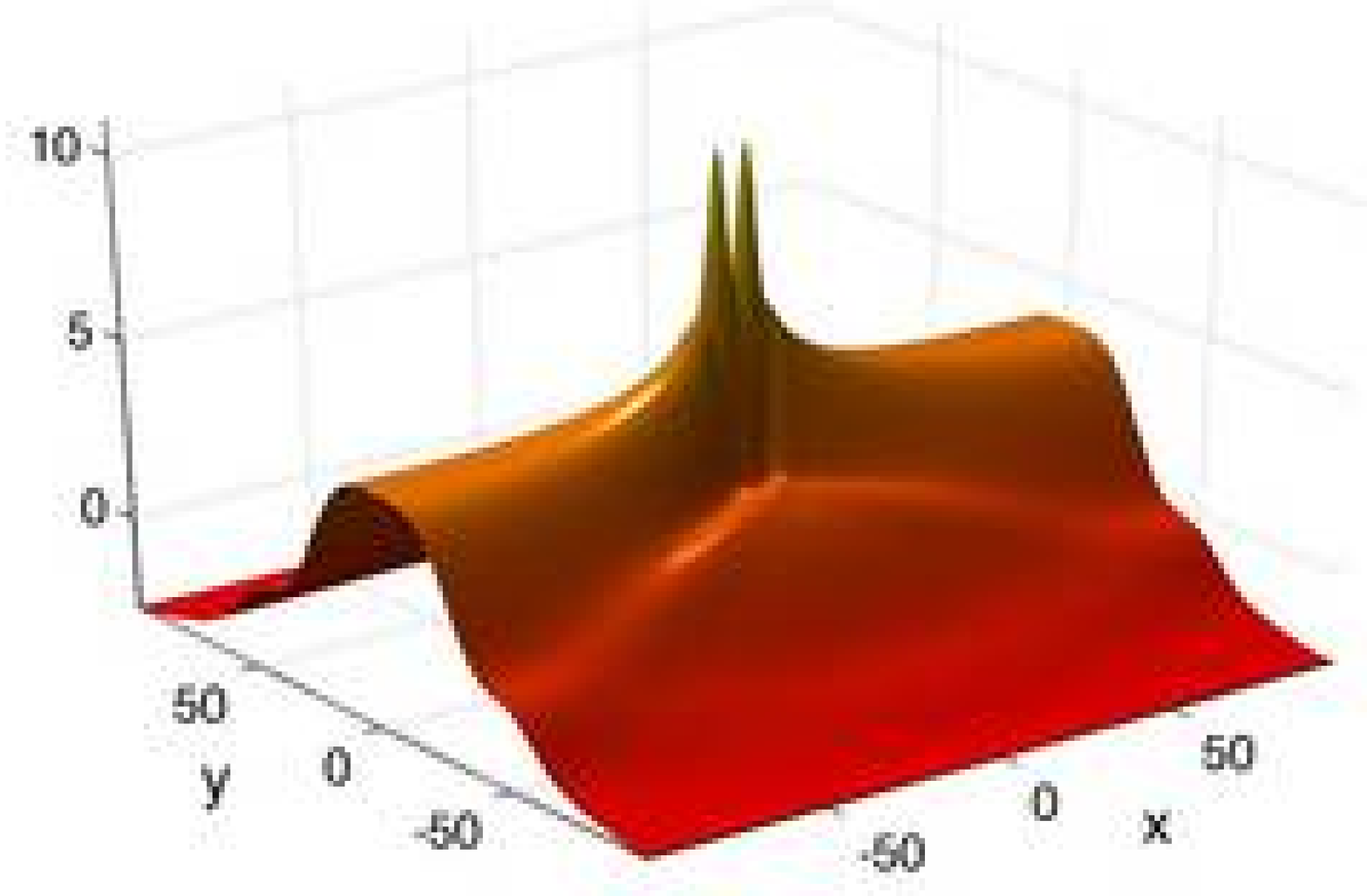}\\{\includegraphics[width=4.064cm]{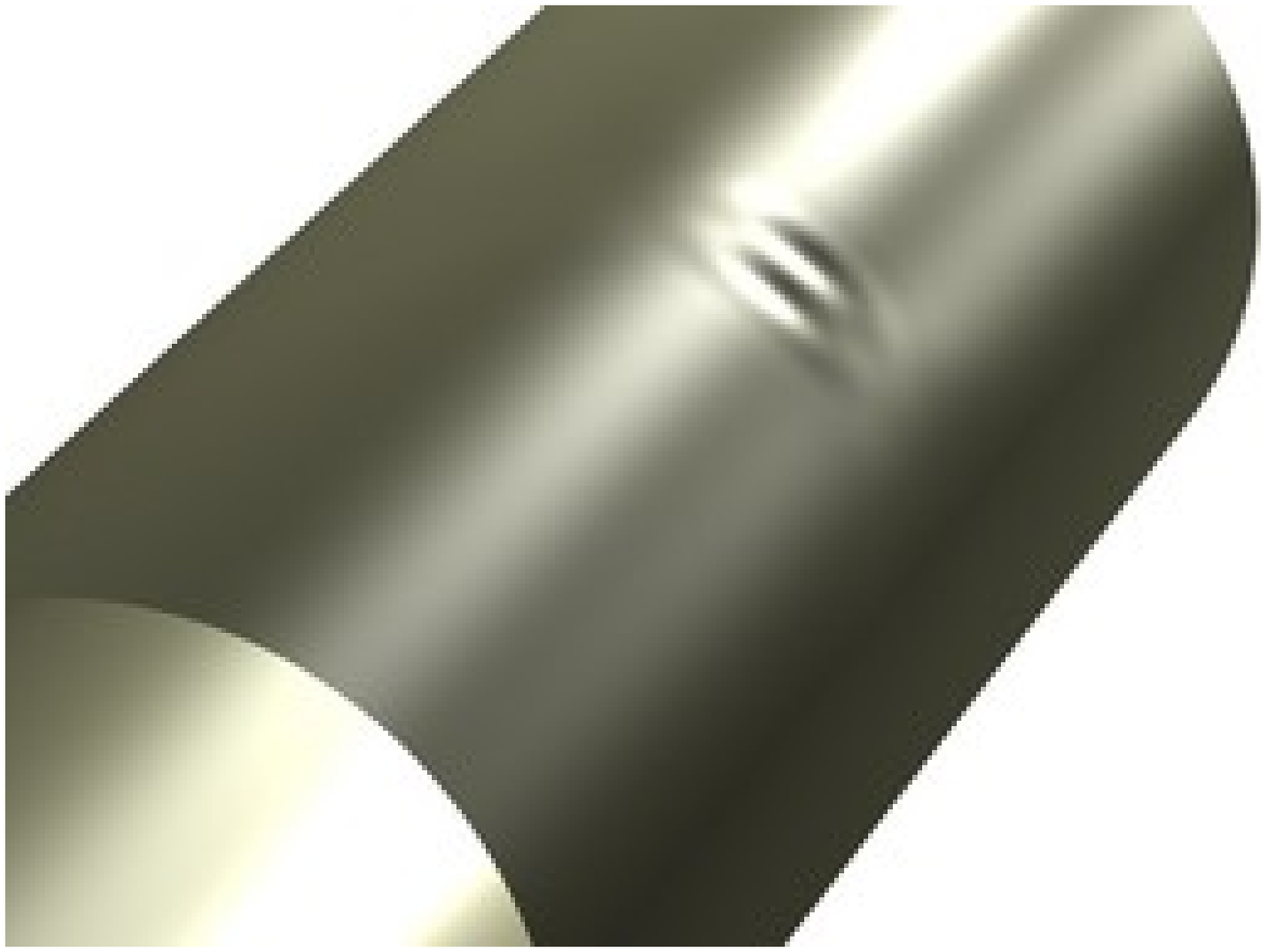}}}}
      \put(92,64){\frparbcenter{\includegraphics[width=4.233cm]{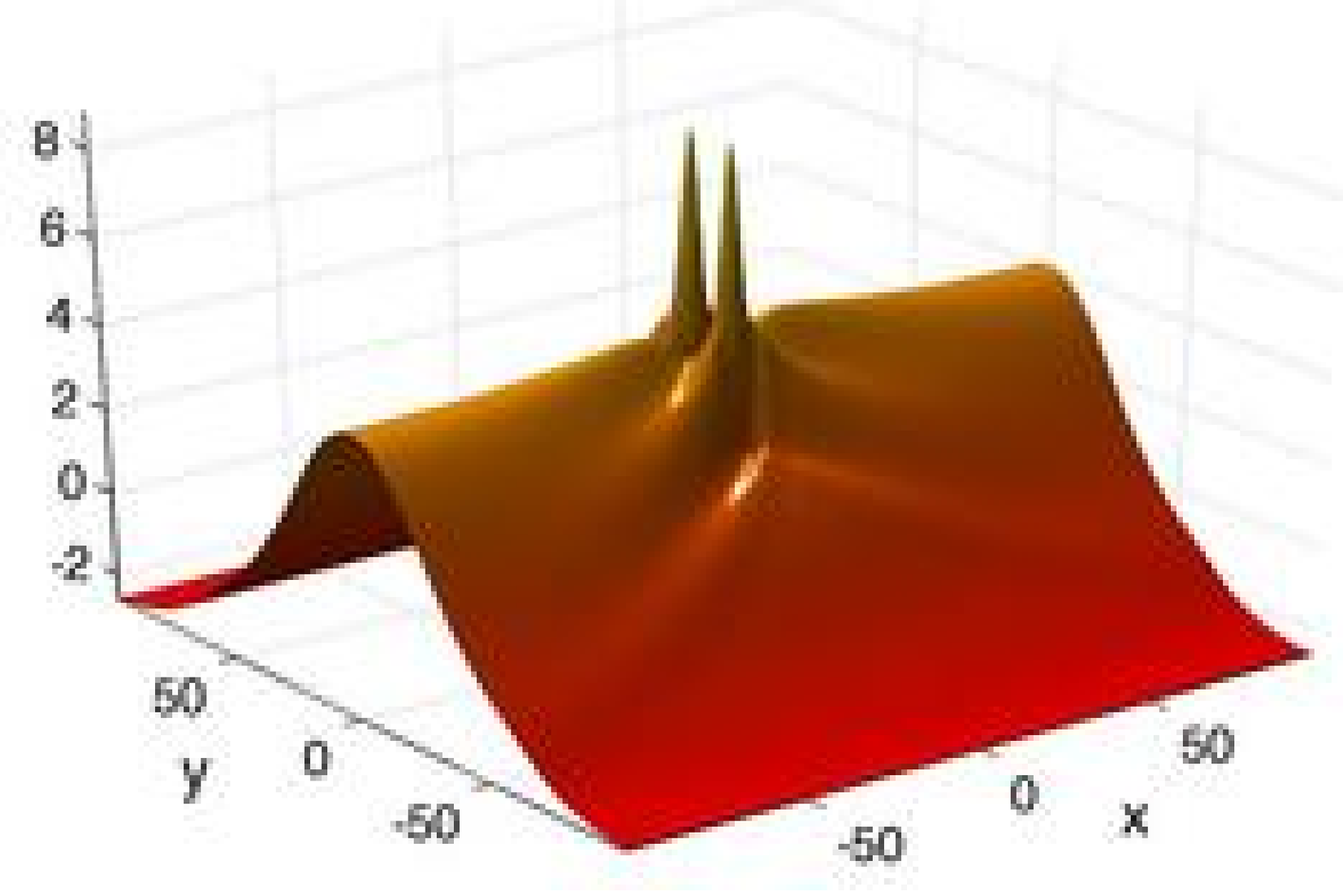}\\{\includegraphics[width=4.064cm]{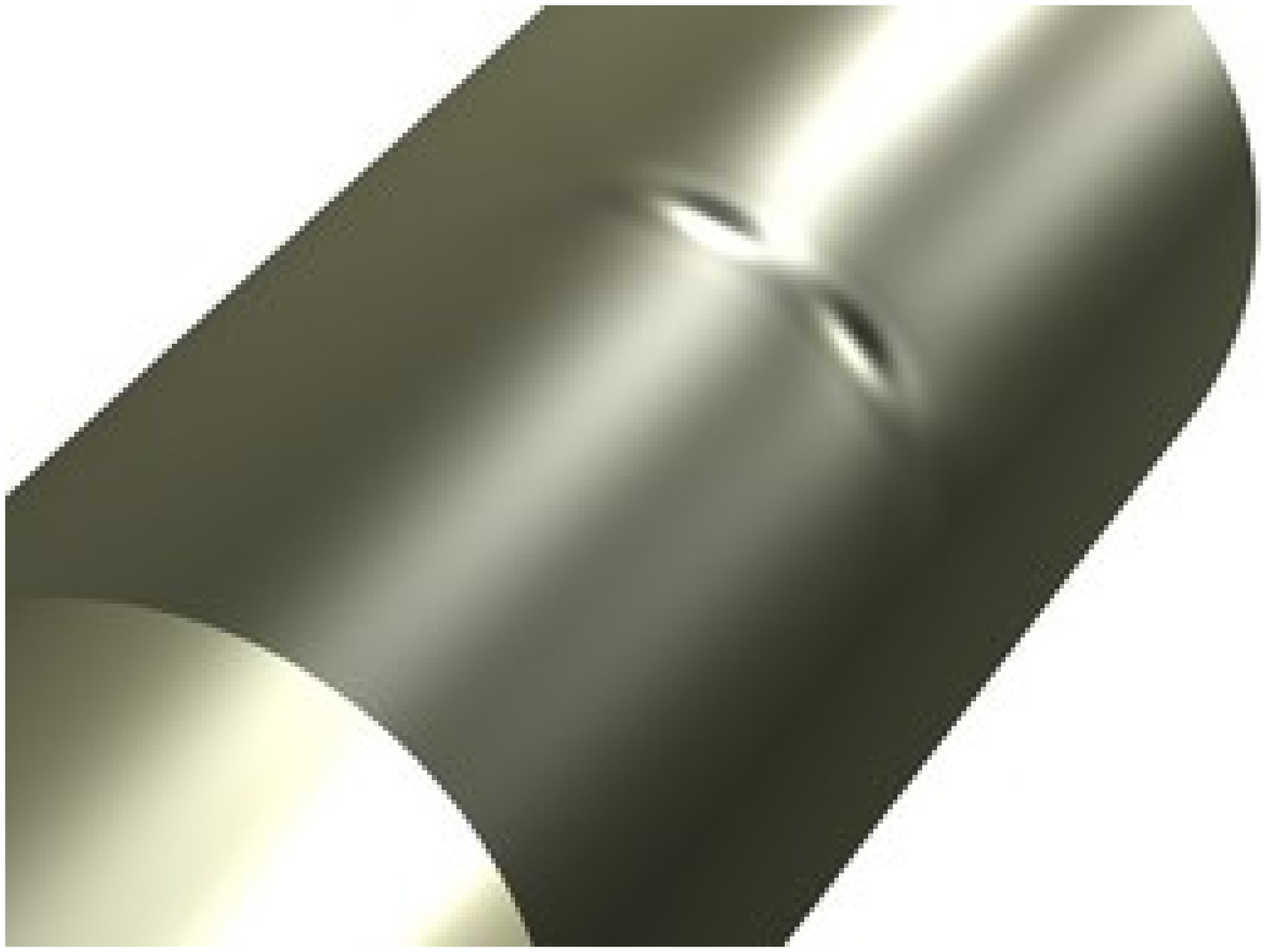}}}}
      \put(138,64){\frparbcenter{\includegraphics[width=4.233cm]{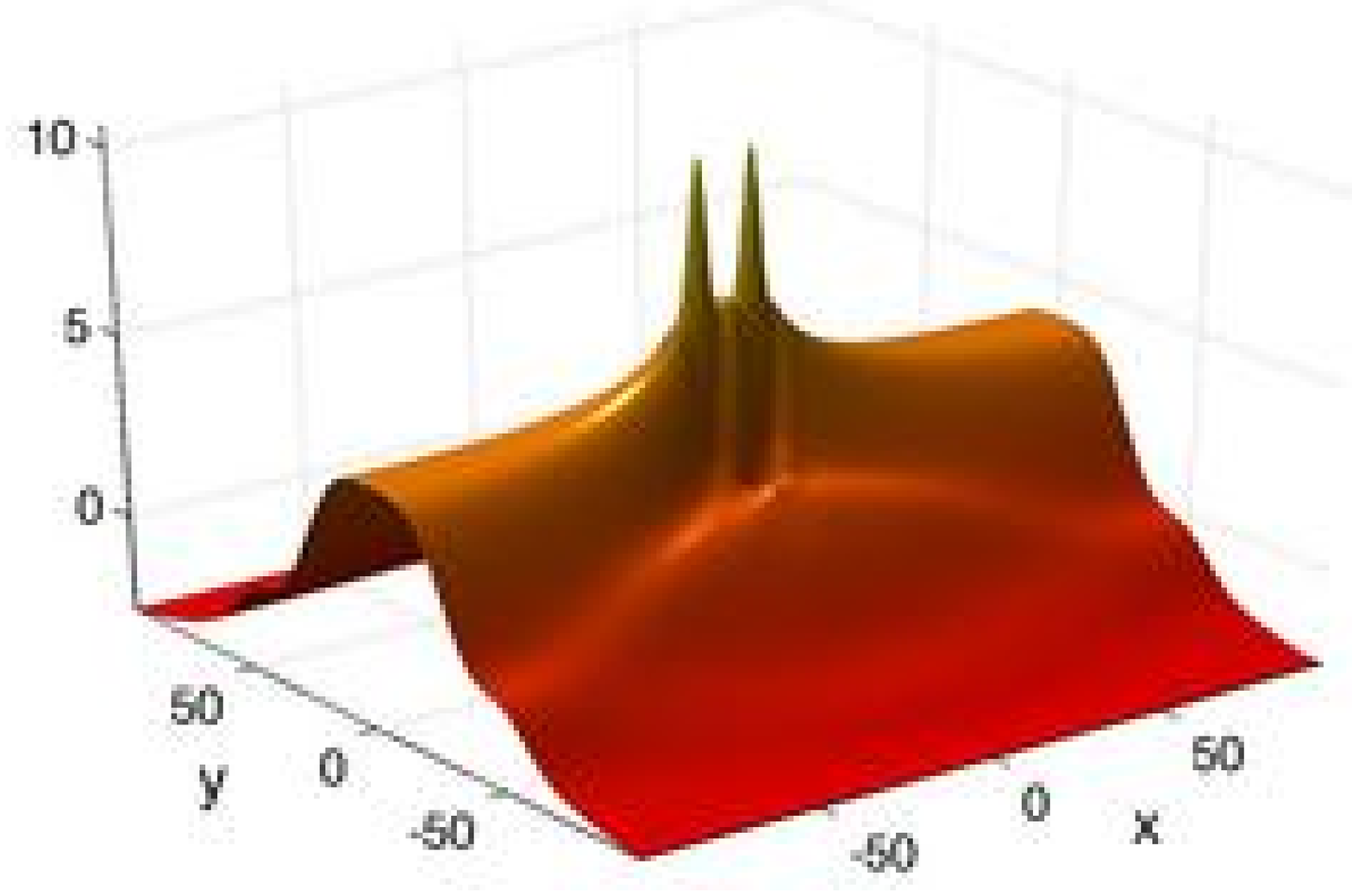}\\{\includegraphics[width=4.064cm]{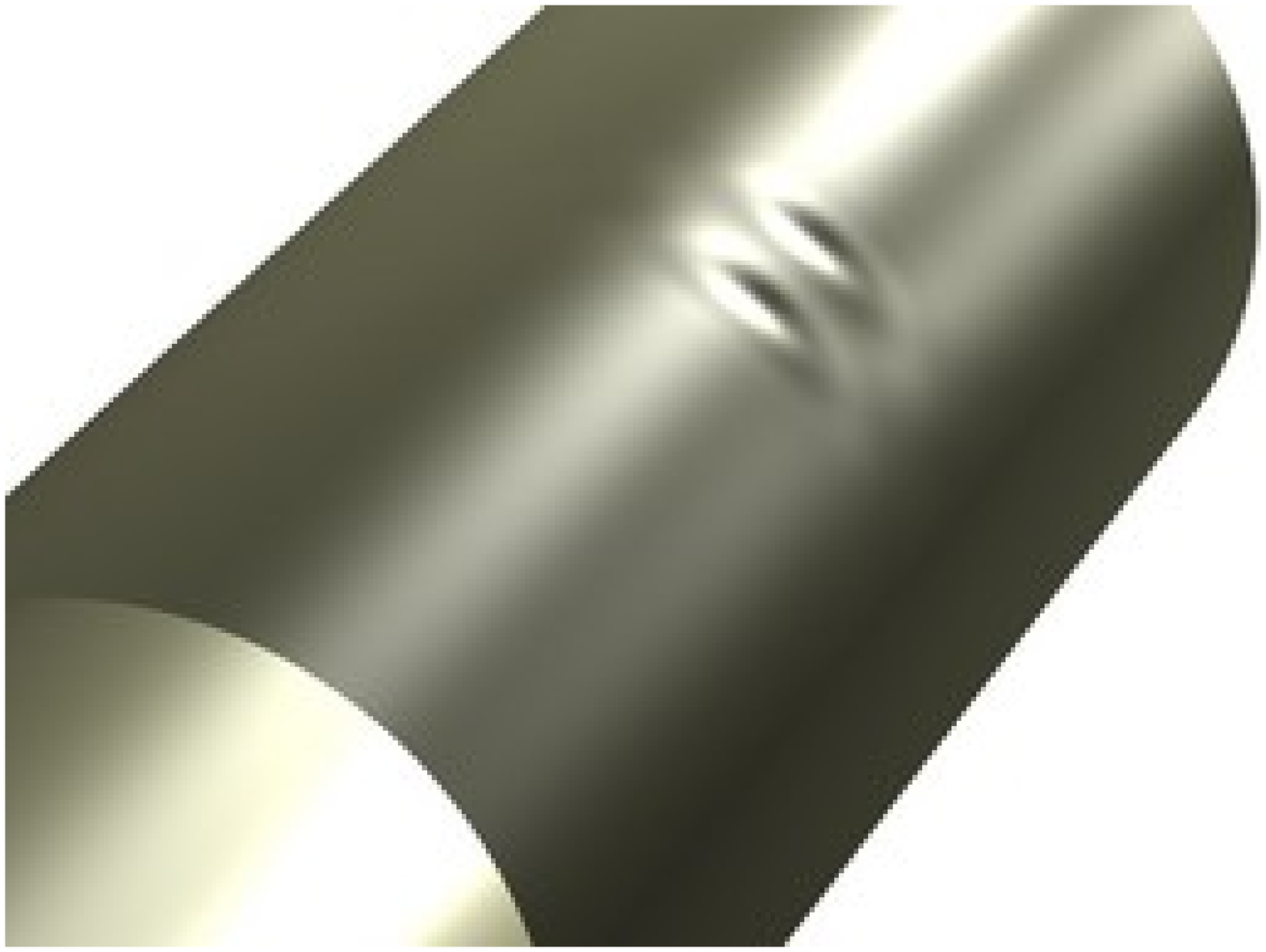}}}}
      \put(0,0){\frparbcenter{\includegraphics[width=4.233cm]{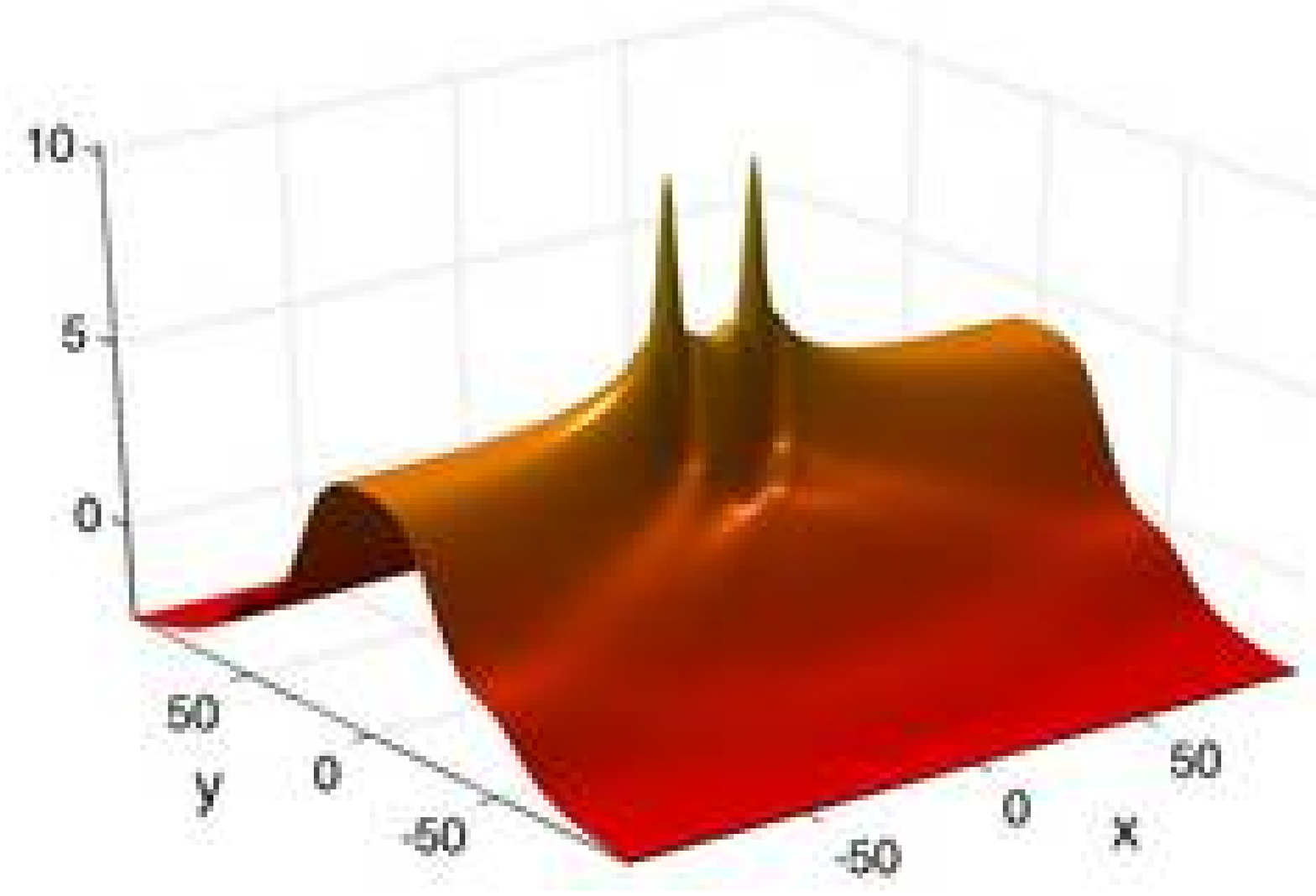}\\{\includegraphics[width=4.064cm]{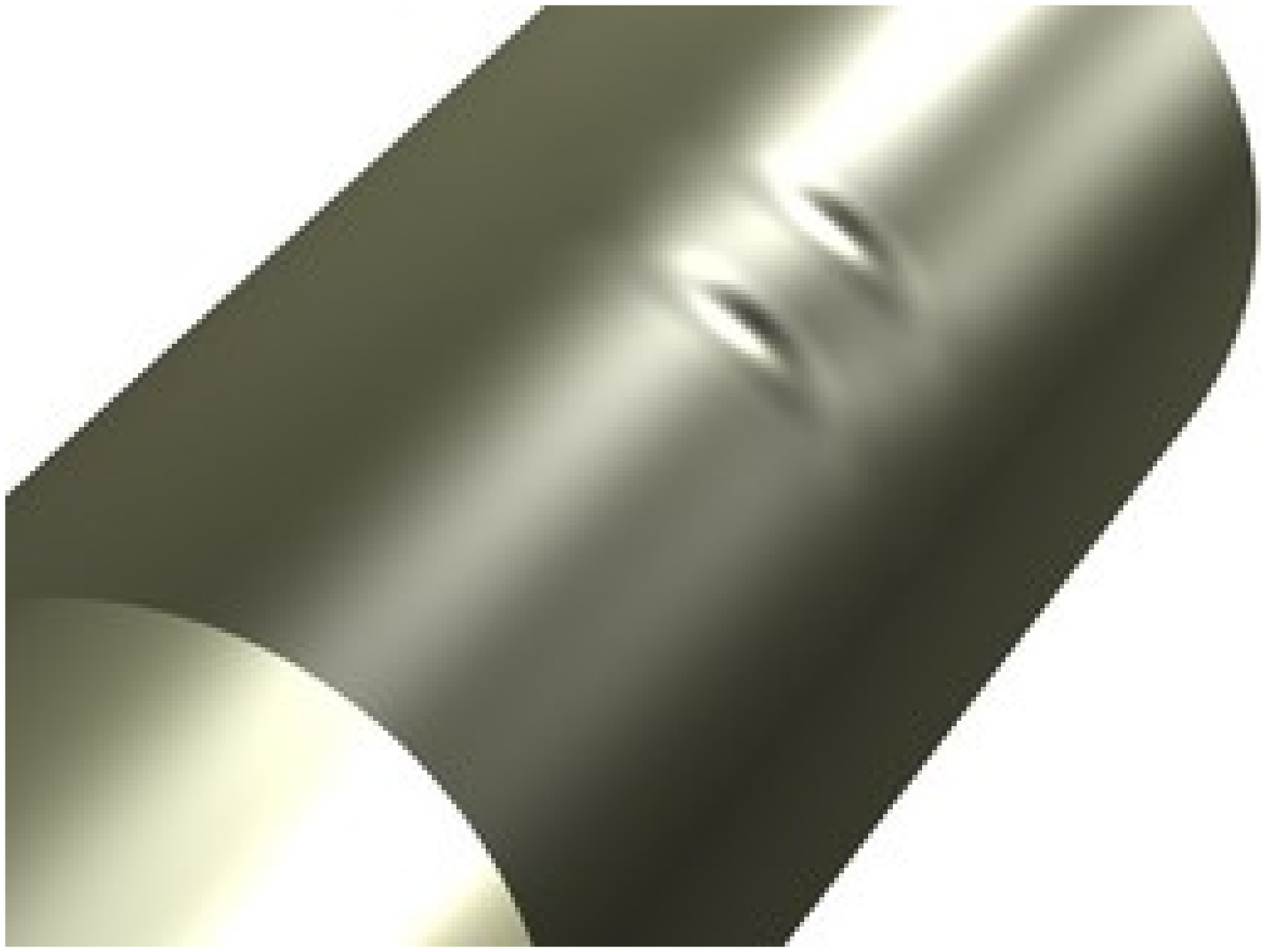}}}}
      \put(46,0){\frparbcenter{\includegraphics[width=4.233cm]{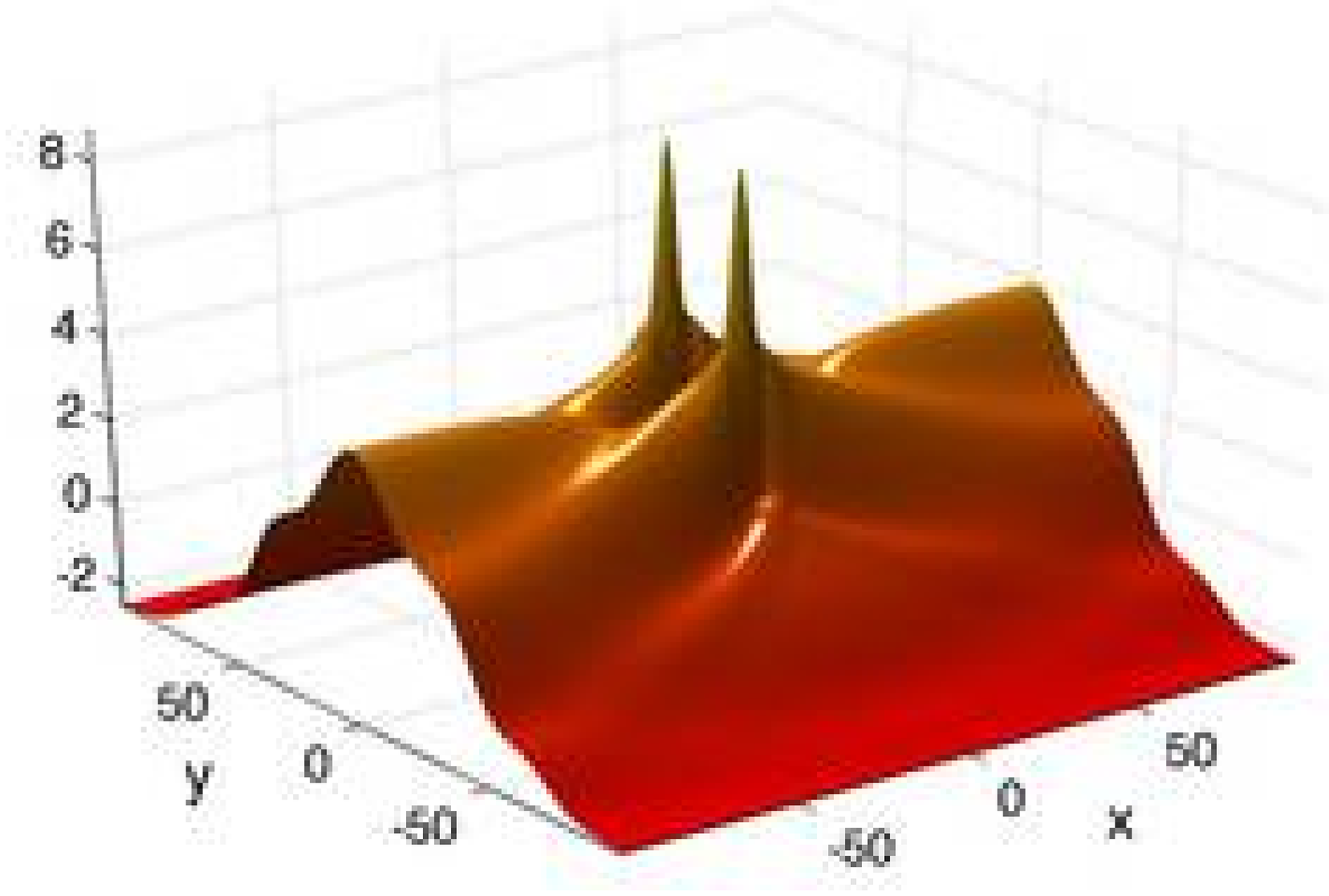}\\{\includegraphics[width=4.064cm]{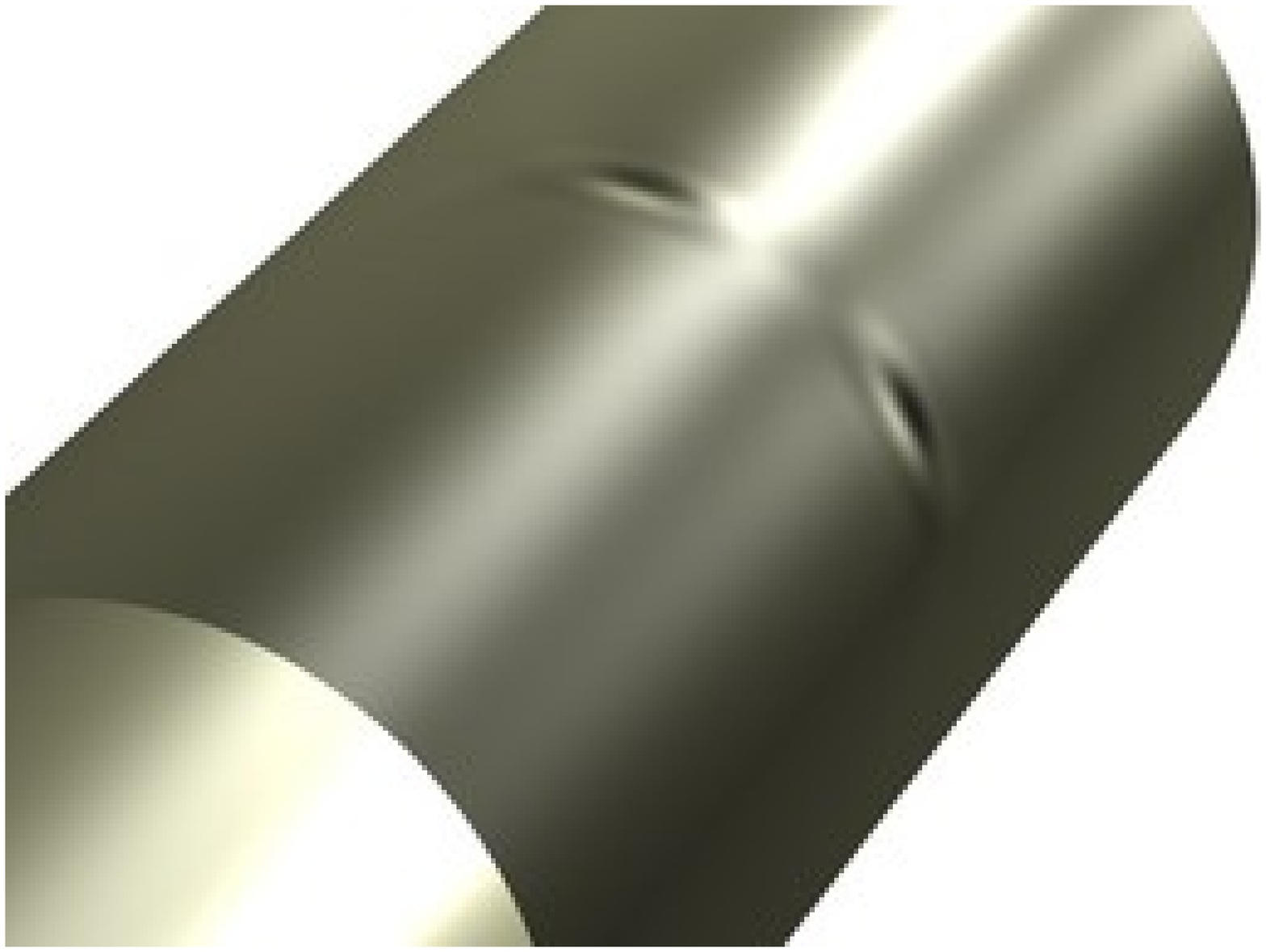}}}}
      \put(92,0){\frparbcenter{\includegraphics[width=4.233cm]{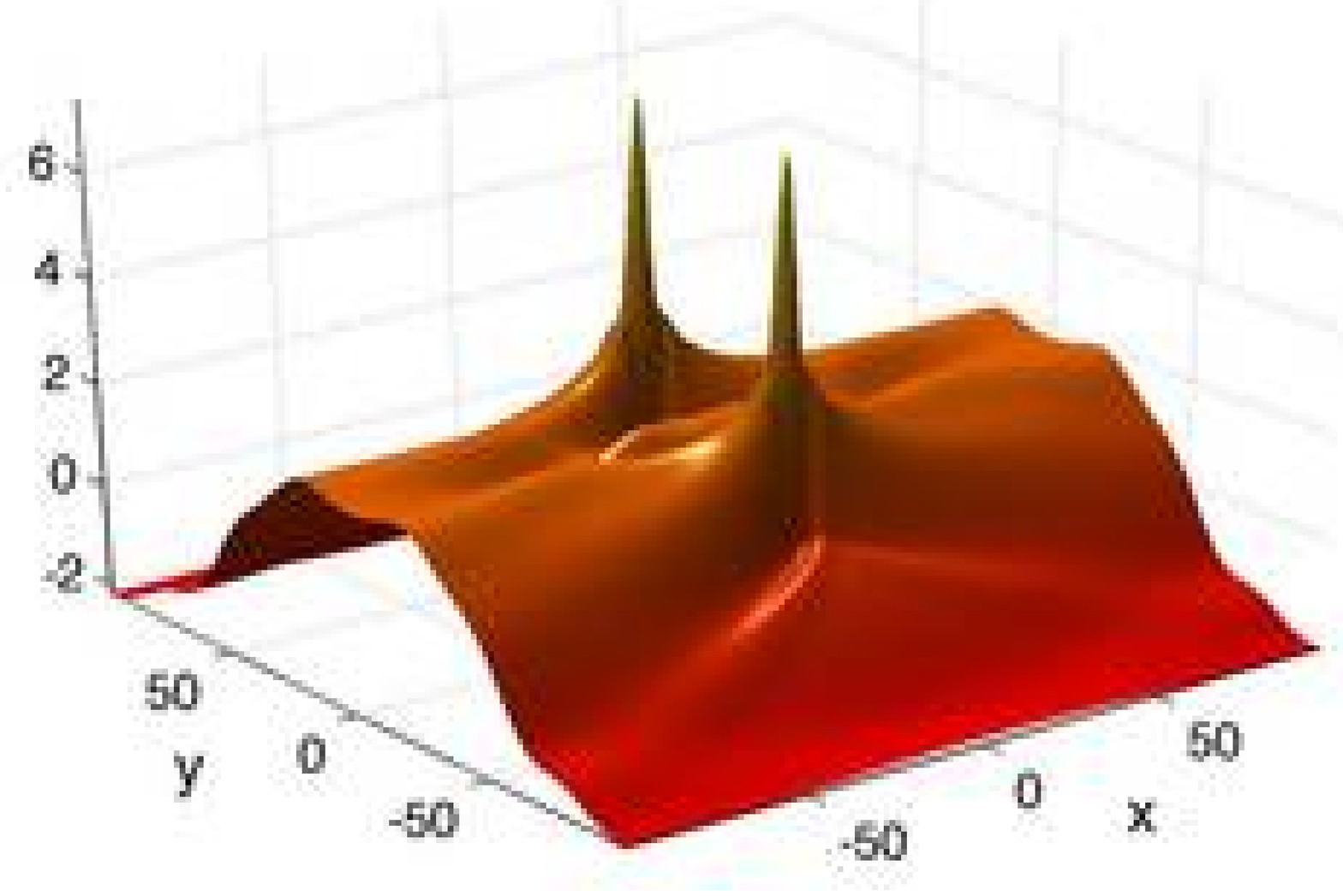}\\{\includegraphics[width=4.064cm]{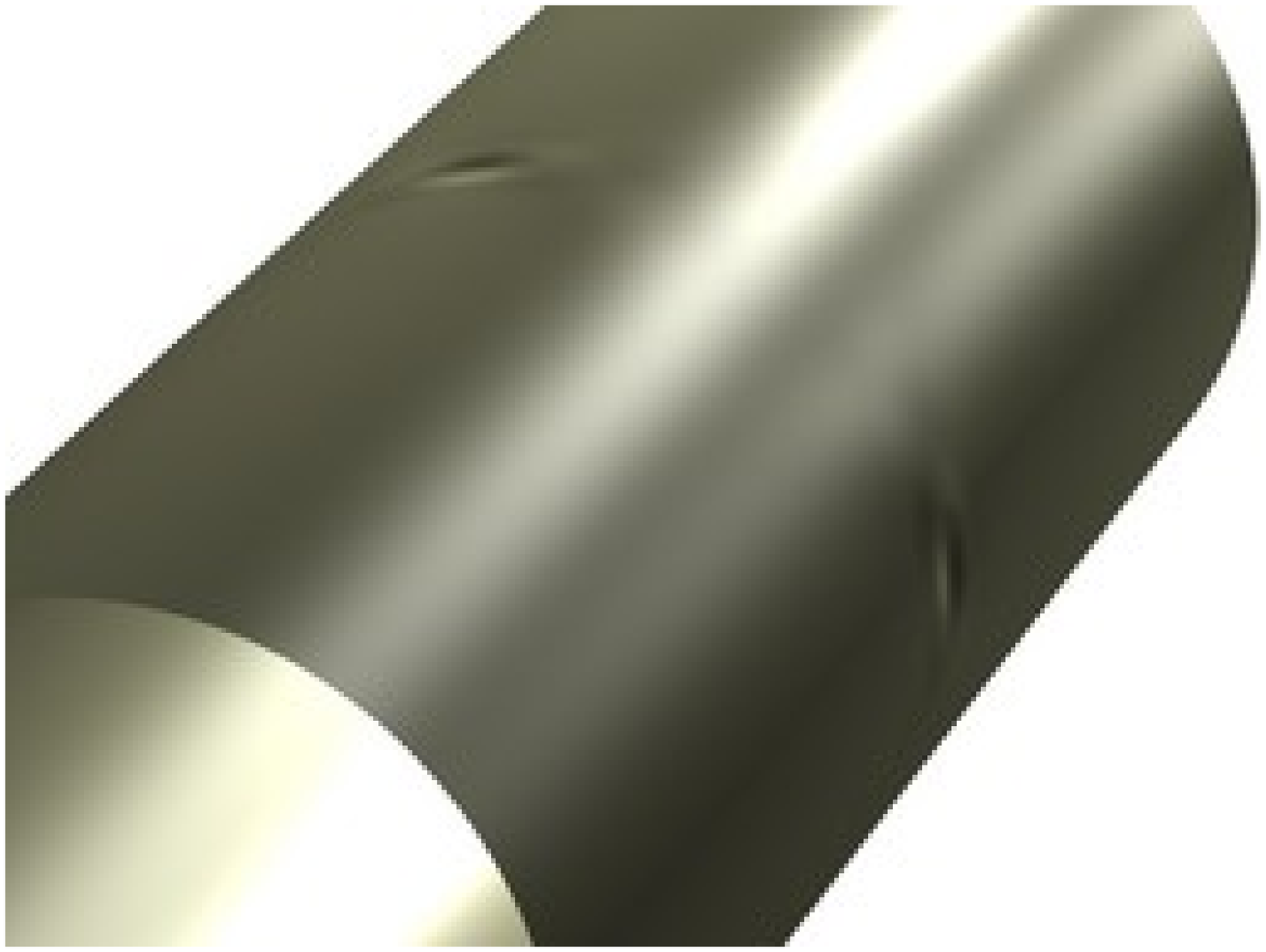}}}}
      \color{black}
      \put(31,65){\makebox(11,5)[rb]{(1.1)}}
      \put(77,65){\makebox(11,5)[rb]{(1.2)}}
      \put(123,65){\makebox(11,5)[rb]{(1.3)}}
      \put(169,65){\makebox(11,5)[rb]{(1.4)}}
      \put(31,1){\makebox(11,5)[rb]{(1.5)}}
      \put(77,1){\makebox(11,5)[rb]{(1.6)}}
      \put(123,1){\makebox(11,5)[rb]{(1.7)}}
    \end{picture}
  \caption{Numerical solutions found using the CSDM with axial end
    shortening $S=40$. More details are in Table~\ref{tab:sol_csdm}.}
  \label{fig:sol_csdm}
\end{figure}
\end{landscape}

We first fixed $S=40$ and used the CSDM to obtain constrained local
minimizers of $E$ described in~Table~\ref{tab:sol_csdm}. They are
ordered according to the increasing value of stored energy $E$. Their
graphs and renderings on a cylinder are shown
in~Fig.~\ref{fig:sol_csdm}. 
Solution 1.1 is similar to the single
dimple solution $\wMP$ described above 
and according to Table~\ref{tab:sol_csdm} it has, indeed,  the
smallest value of $E$.

\begin{table}[htbp]
  \centering
  \begin{tabular}{|c||c|c|c|c|c|}
    \hline
    {\bf CSDM} & $\l$ & $S$ & $E$ & $F_\l$ & same shape as {\bf MPA}\tabularnewline
    \hline
    1.1 & 1.108121 & 40 & 56.85636 & 12.53151 & 2.1 \tabularnewline 
    \hline
    1.2 & 1.299143 & 40 & 62.76150 & 10.79577 & 2.2 \tabularnewline 
    \hline
    1.3 & 1.316146 & 40 & 63.21646 & 10.57063 & 2.3 \tabularnewline 
    \hline
    1.4 & 1.311687 & 40 & 63.64083 & 11.17334 & 2.4 \tabularnewline 
    \hline
    1.5 & 1.309586 & 40 & 63.70623 & 11.32278 & 2.5 \tabularnewline 
    \hline
    1.6 & 1.328997 & 40 & 64.00875 & 10.84889 & 2.6 \tabularnewline 
    \hline
    1.7 & 1.344898 & 40 & 64.52244 & 10.72651 & 2.7 \tabularnewline 
    \hline
  \end{tabular}
  \vspace{3mm}
  \caption{Numerical solutions obtained by the CSDM on $\Omq$ with $\partial_{xy}$ discretized using left-sided finite differences. Graphs are shown in~Fig.~\ref{fig:sol_csdm}.}
  \label{tab:sol_csdm}
\end{table}

\subsubsection{Mountain-pass algorithm}
We then used the MPA for fixed $\l=1.4$ and various choices of $w_2$ to
obtain the local mountain-pass points of $F_\l$ described
in~Table~\ref{tab:sol_mpa}. They are ordered according to the
increasing value of the total potential $F_\l$. The shape of their
graph is very similar to that of the CSDM solutions discussed above
and depicted in Fig.~\ref{fig:sol_csdm} and we do not show
their graphs here. Solution 2.1 is again the single dimple solution and the
table shows that it has the smallest value of $F_\l$.  

\begin{table}[htbp]
  \centering
  \begin{tabular}{|c||c|c|c|c|c|}
    \hline
    {\bf MPA} & $\l$ & $S$ & $E$ & $F_\l$ & same shape as {\bf CSDM}\tabularnewline
    \hline
    2.1 & 1.4 & 17.73822 & 29.42997 & 4.596460 & 1.1 \tabularnewline 
    \hline
    2.2 & 1.4 & 29.85121 & 49.08882 & 7.297132 & 1.2 \tabularnewline 
    \hline
    2.3 & 1.4 & 31.28849 & 51.39952 & 7.595635 & 1.3 \tabularnewline 
    \hline
    2.4 & 1.4 & 31.41723 & 52.01893 & 8.034809 & 1.4 \tabularnewline 
    \hline
    2.5 & 1.4 & 31.22992 & 51.84074 & 8.118852 & 1.5 \tabularnewline 
    \hline
    2.6 & 1.4 & 32.77491 & 54.15818 & 8.273314 & 1.6 \tabularnewline 
    \hline
    2.7 & 1.4 & 34.19888 & 56.56472 & 8.686284 & 1.7 \tabularnewline 
    \hline
  \end{tabular}
  \vspace{3mm}
  \caption{Numerical solutions obtained by the MPA on $\Omq$ with $\partial_{xy}$ discretized using left-sided finite differences.}
  \label{tab:sol_mpa}
\end{table}

\subsubsection{Constrained mountain-pass algorithm}
\begin{landscape}
\begin{figure}[!ht]
\begin{center}
\setlength{\unitlength}{1mm}
\begin{picture}(227,124)
\color[rgb]{.5,.5,.5}
\put(0,64){\frparbcenter{\includegraphics[width=4.233cm]{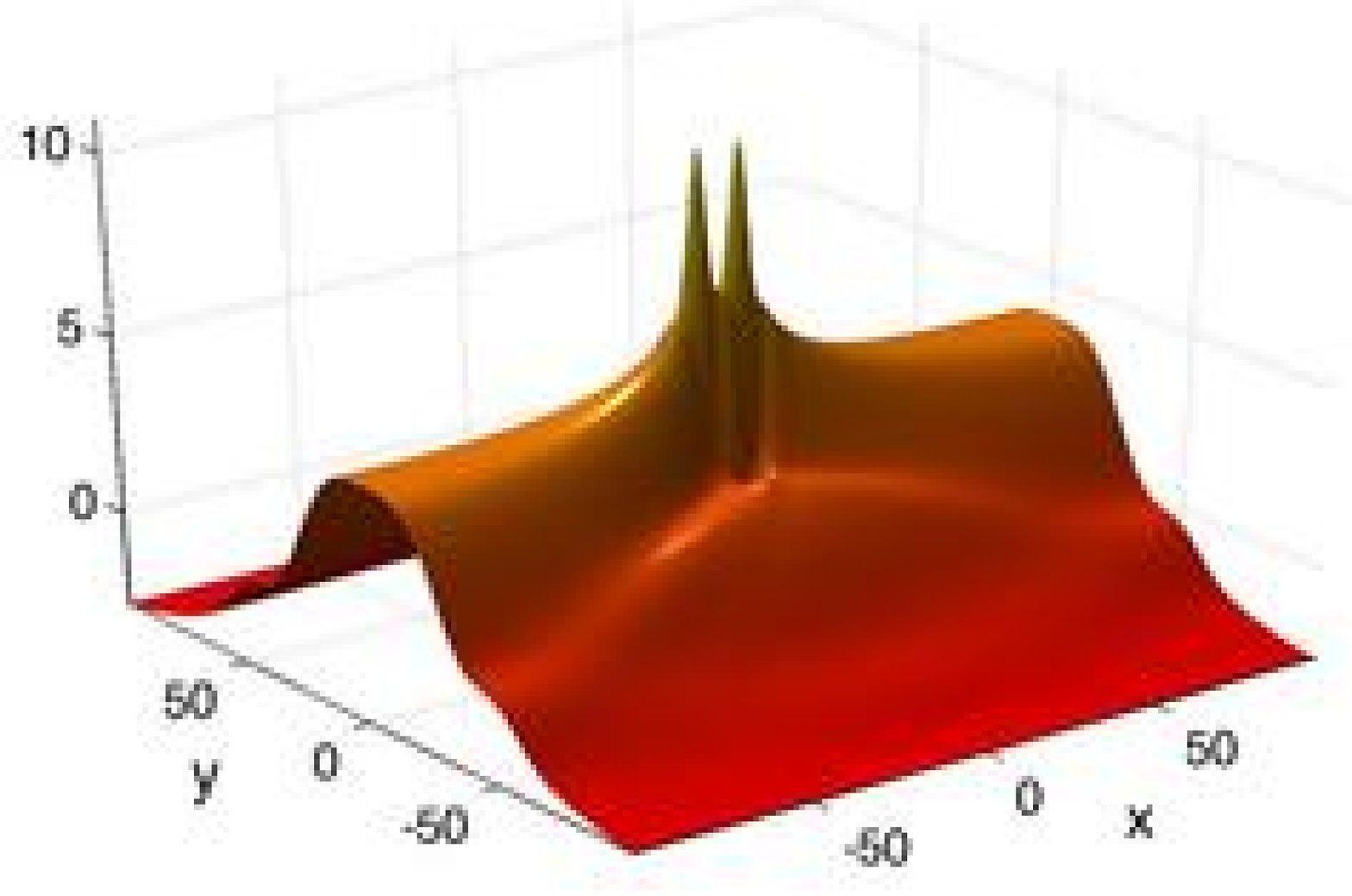}\\\includegraphics[width=4.064cm]{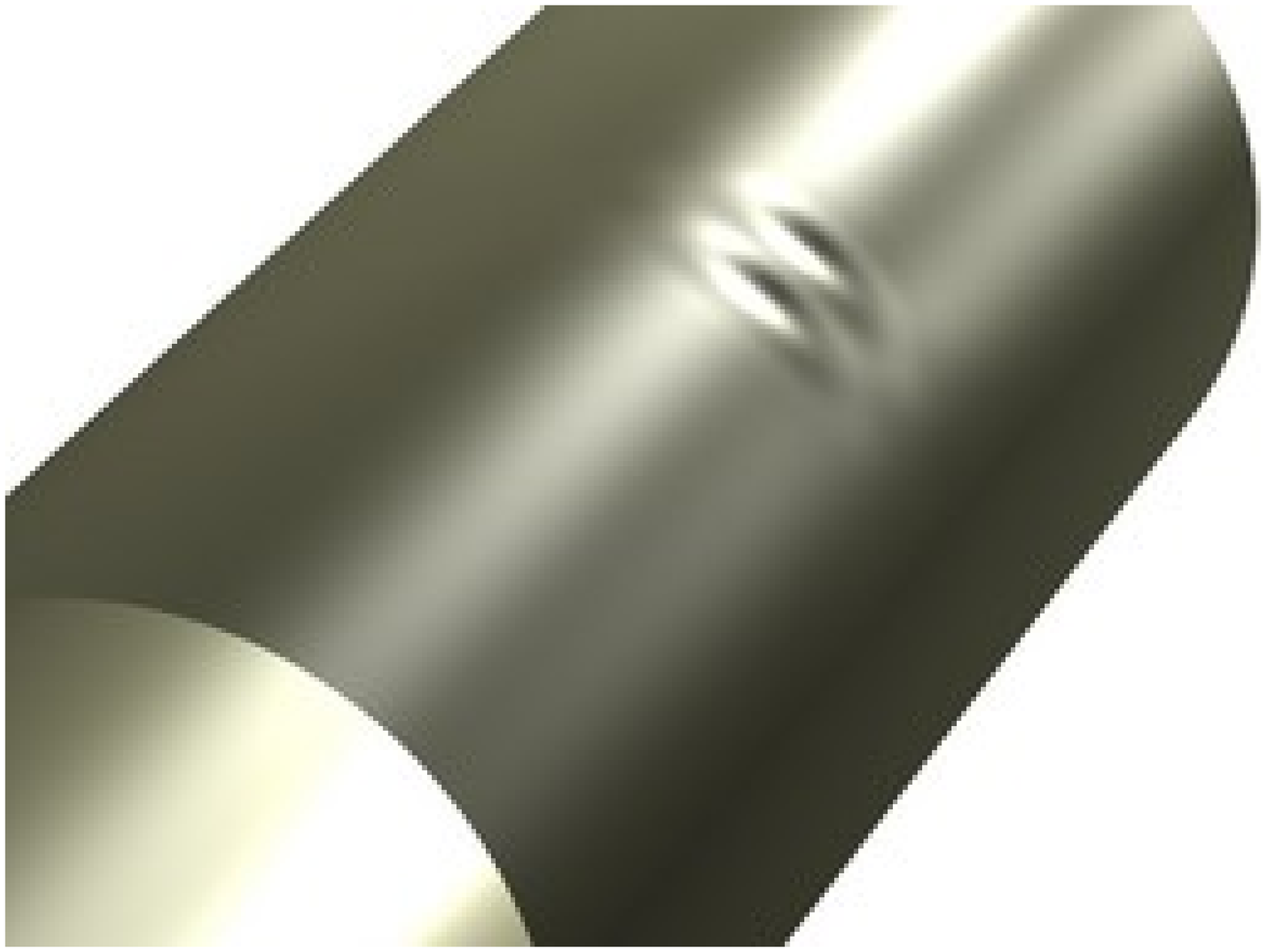}}}
\put(46,64){\frparbcenter{\includegraphics[width=4.233cm]{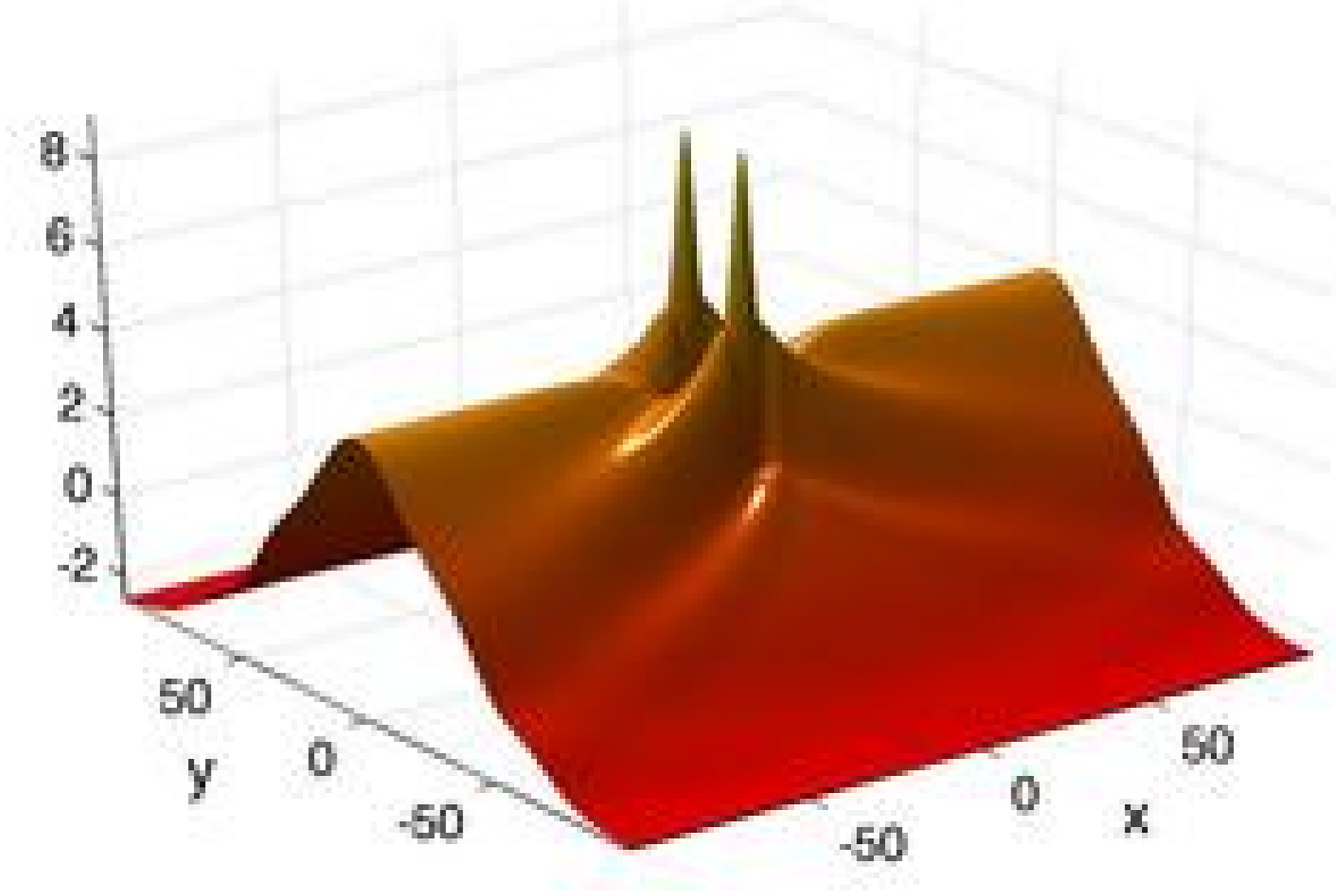}\\\includegraphics[width=4.064cm]{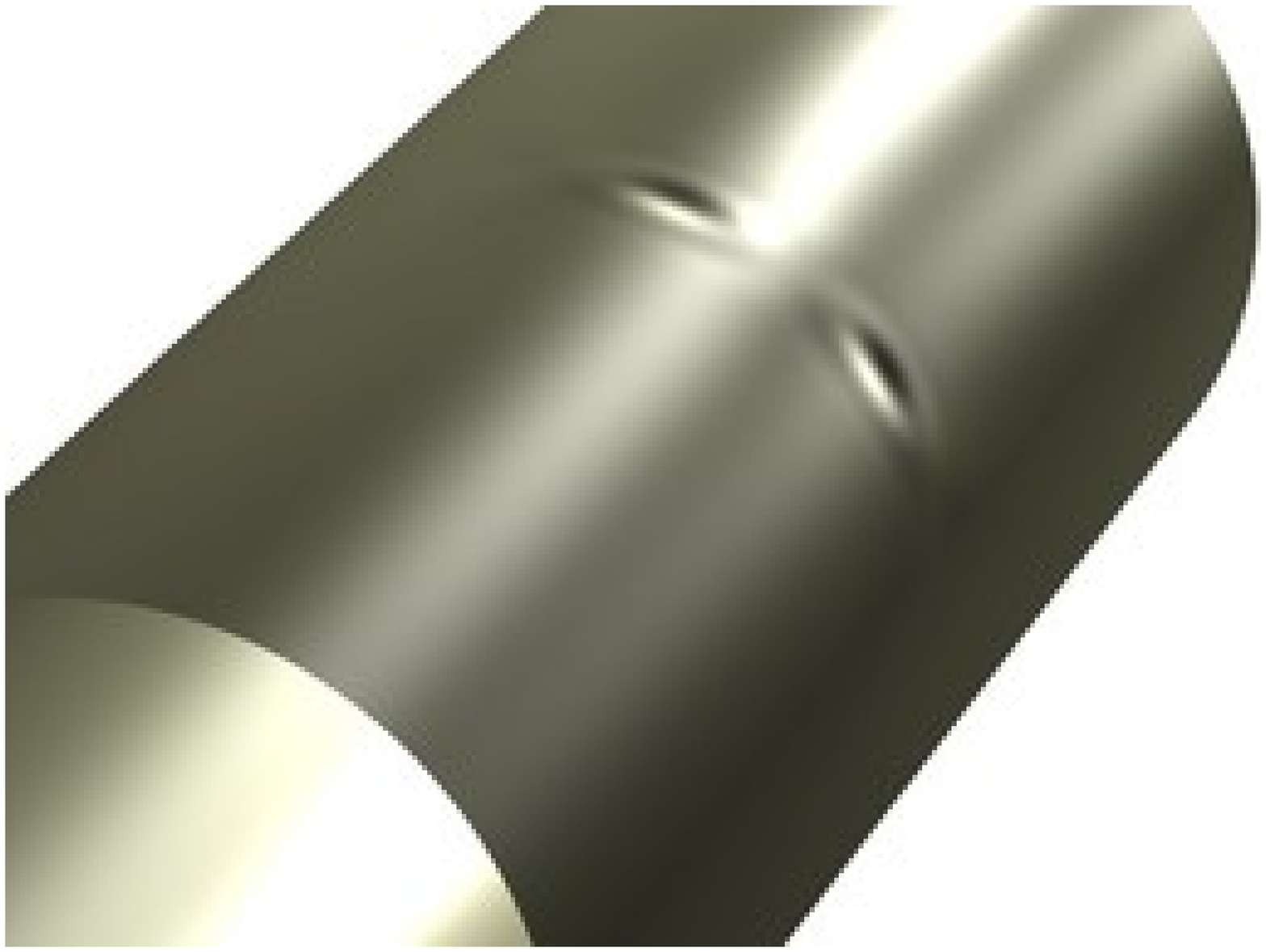}}}
\put(92,64){\frparbcenter{\includegraphics[width=4.233cm]{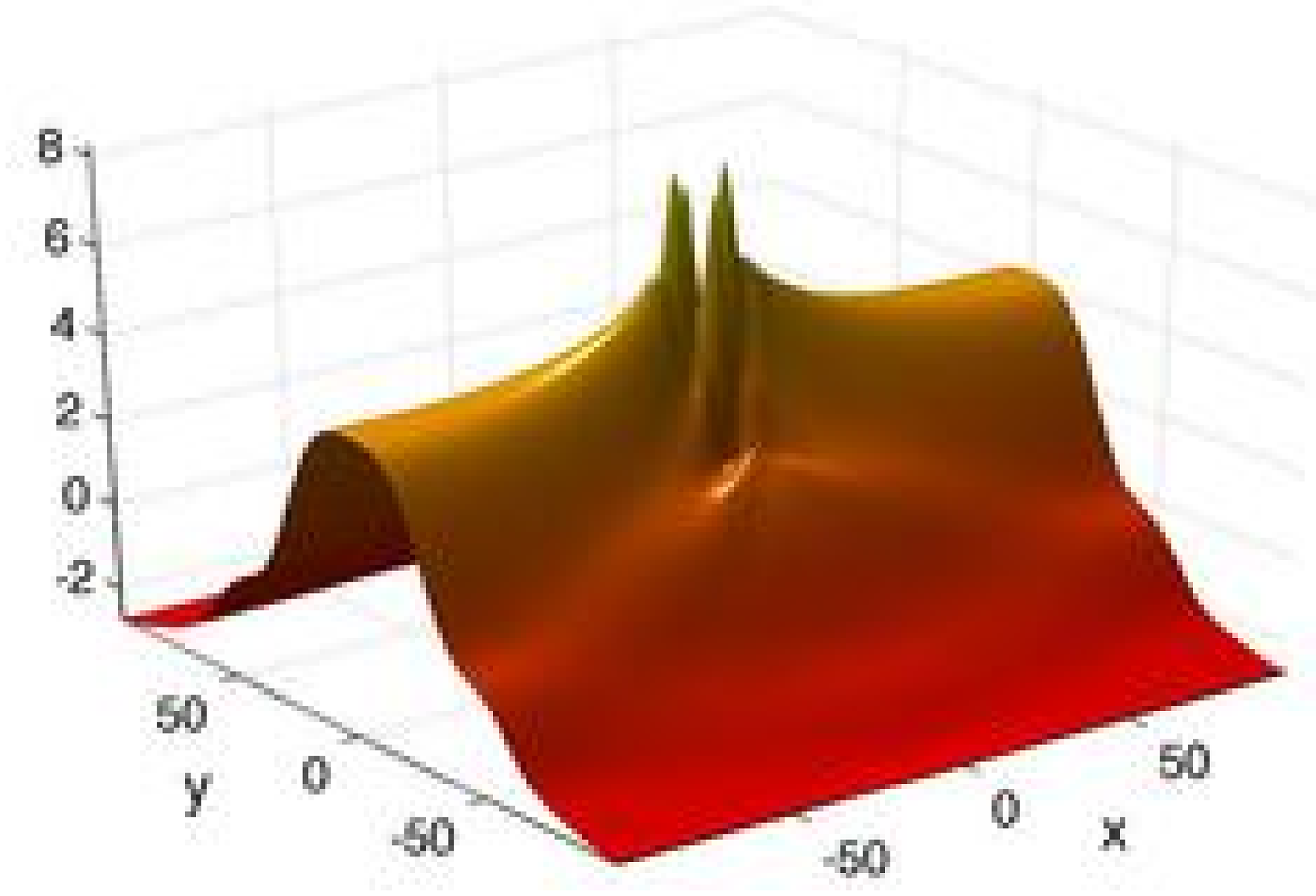}\\\includegraphics[width=4.064cm]{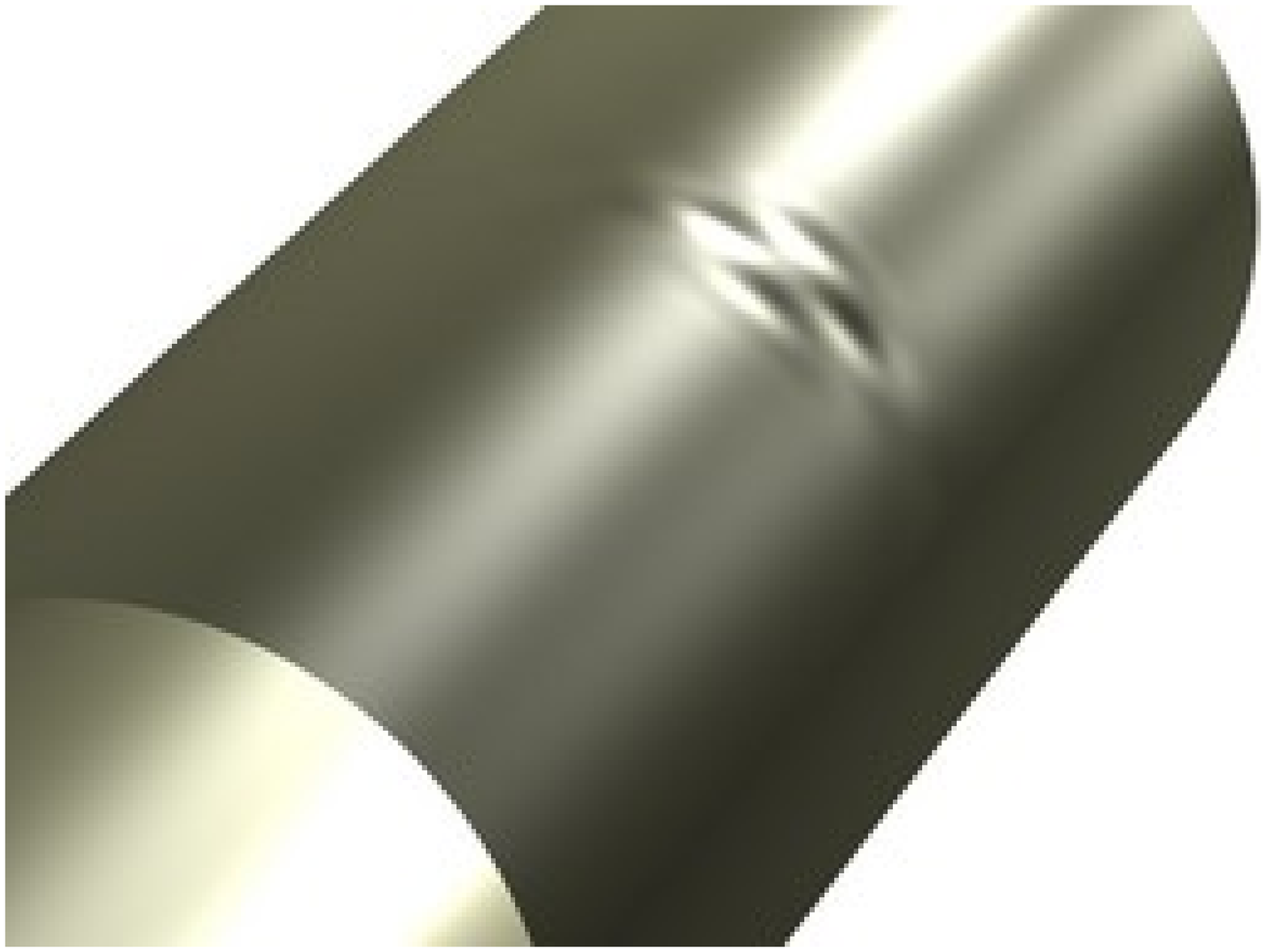}}}
\put(138,64){\frparbcenter{\includegraphics[width=4.233cm]{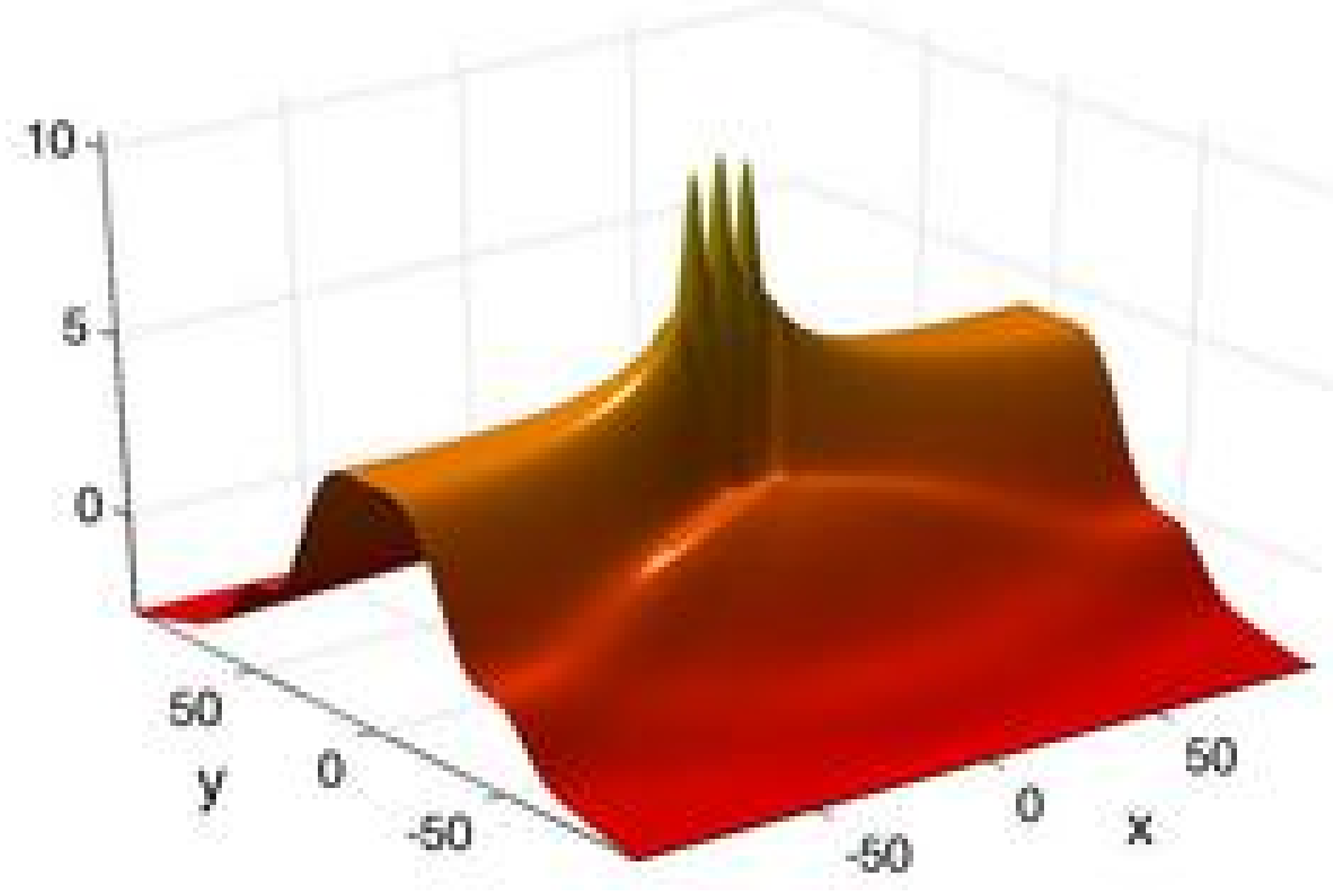}\\\includegraphics[width=4.064cm]{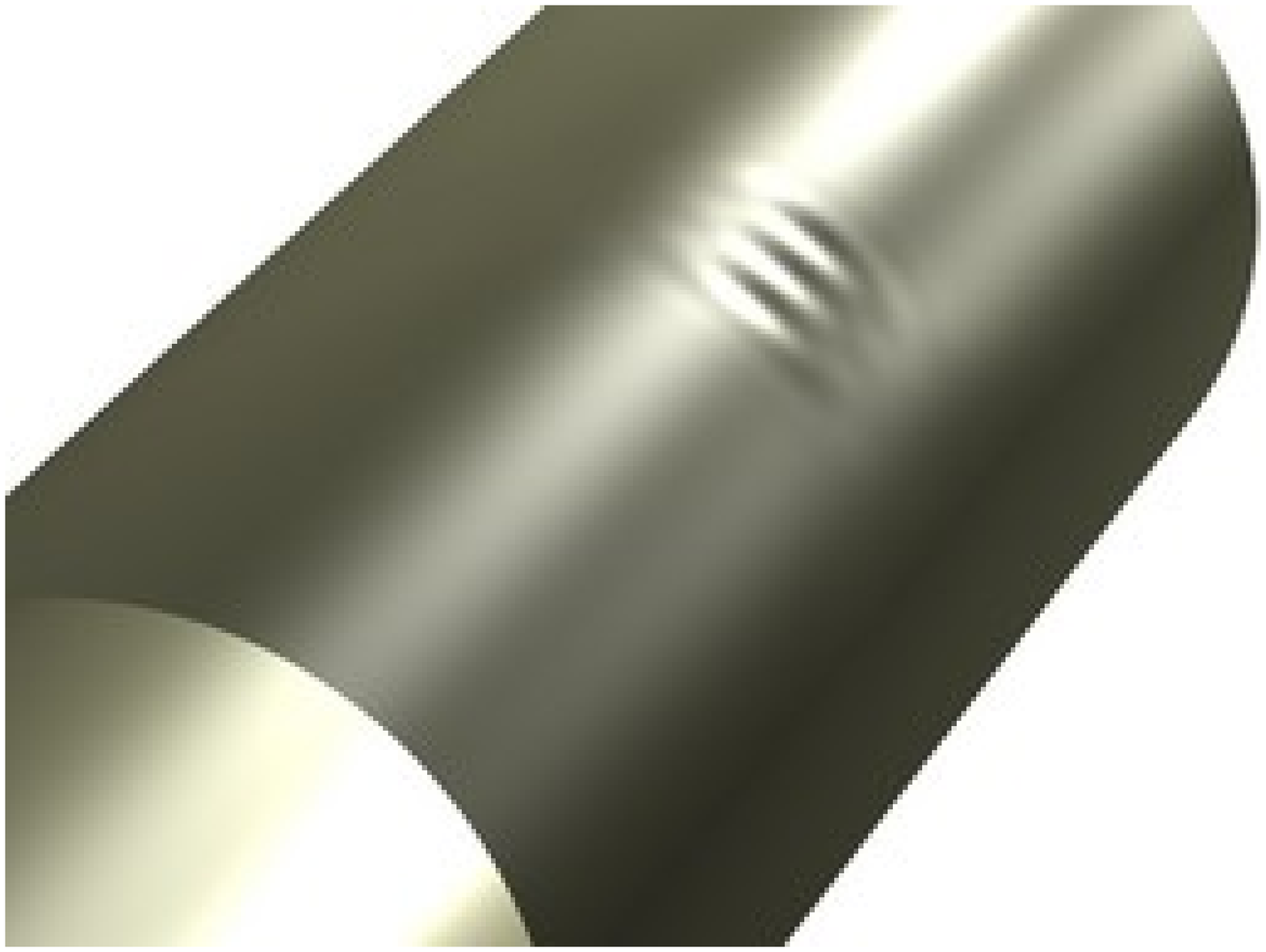}}}
\put(184,64){\frparbcenter{\includegraphics[width=4.233cm]{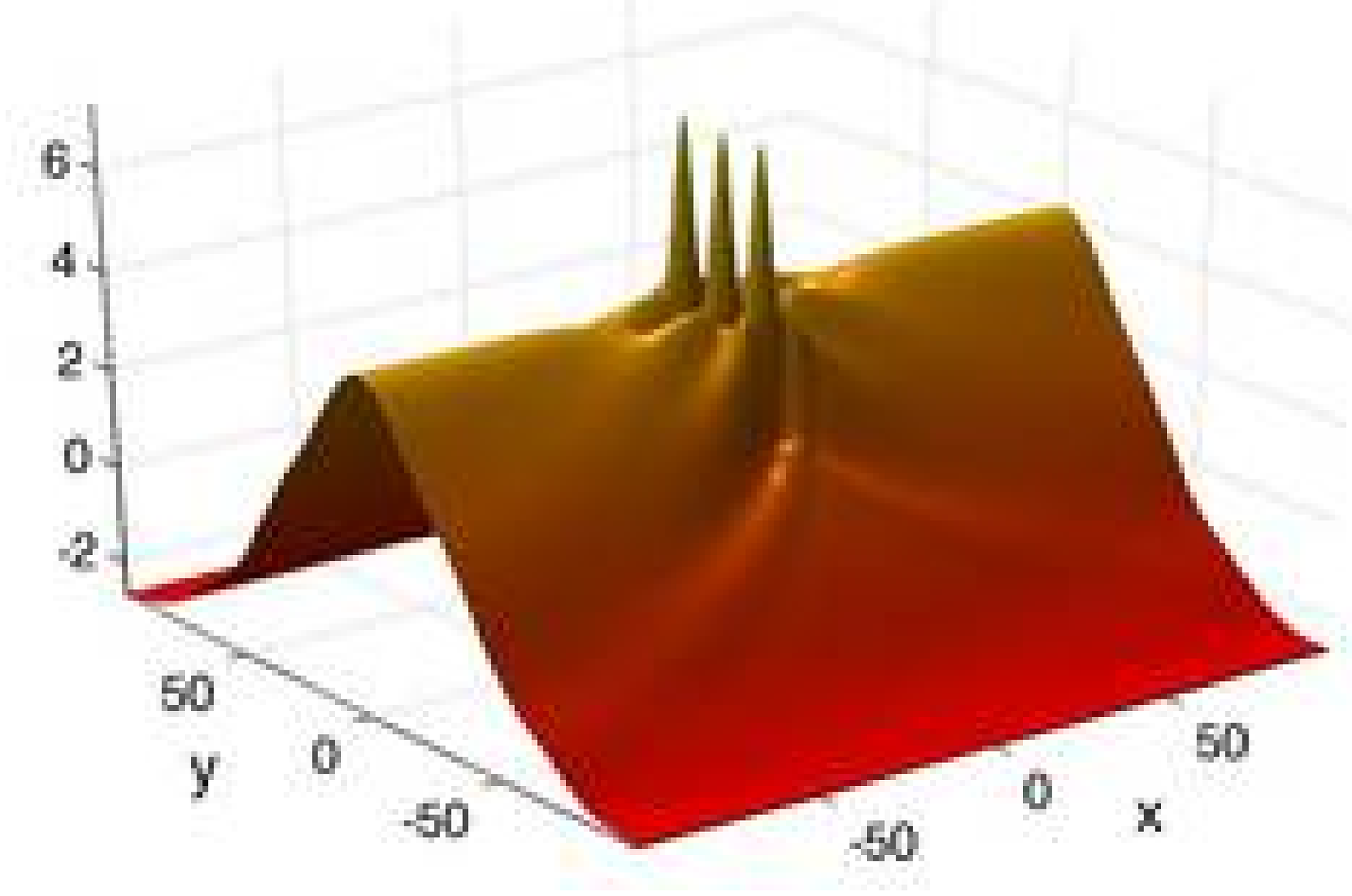}\\\includegraphics[width=4.064cm]{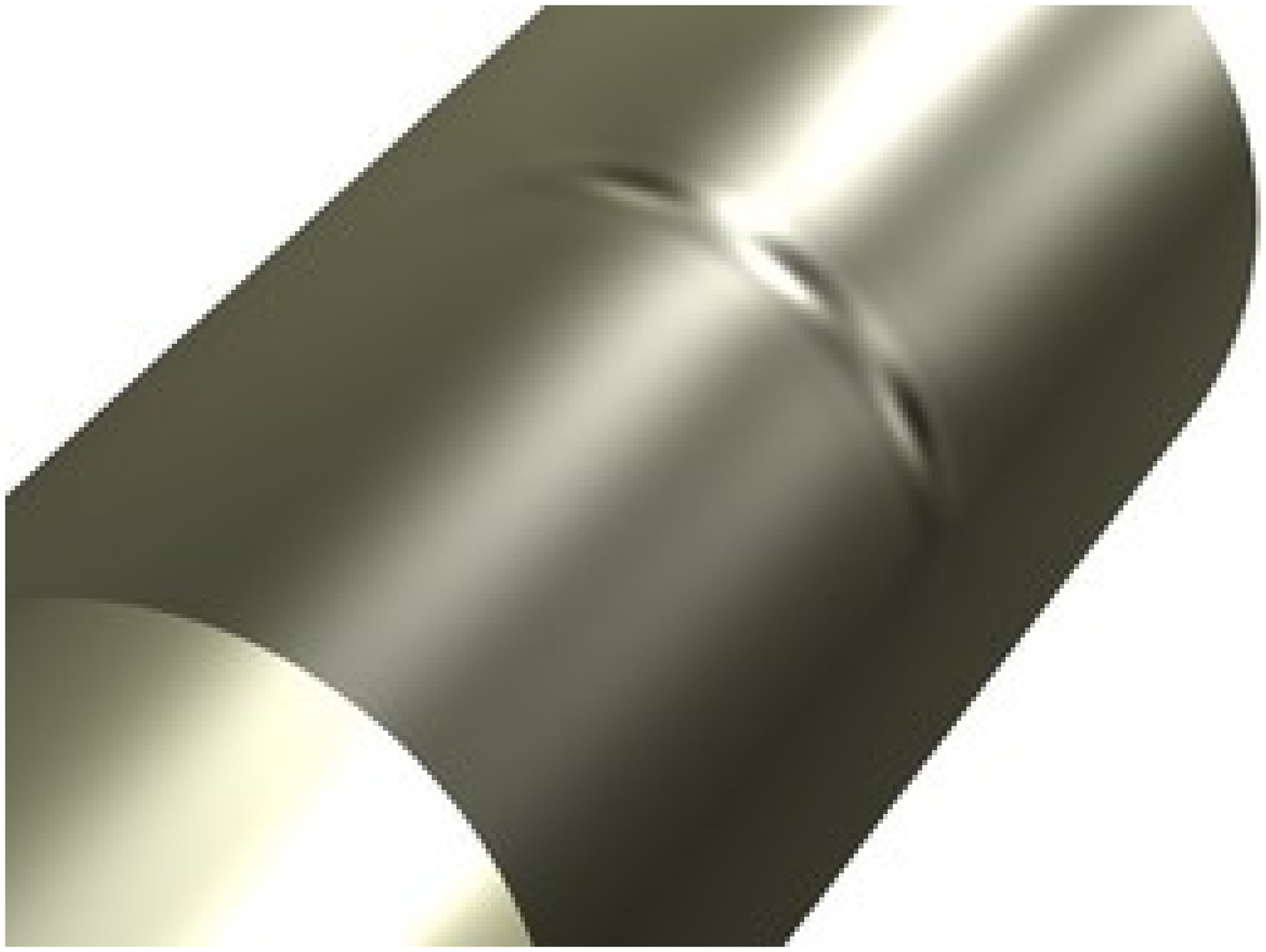}}}
\put(0,0){\frparbcenter{\includegraphics[width=4.233cm]{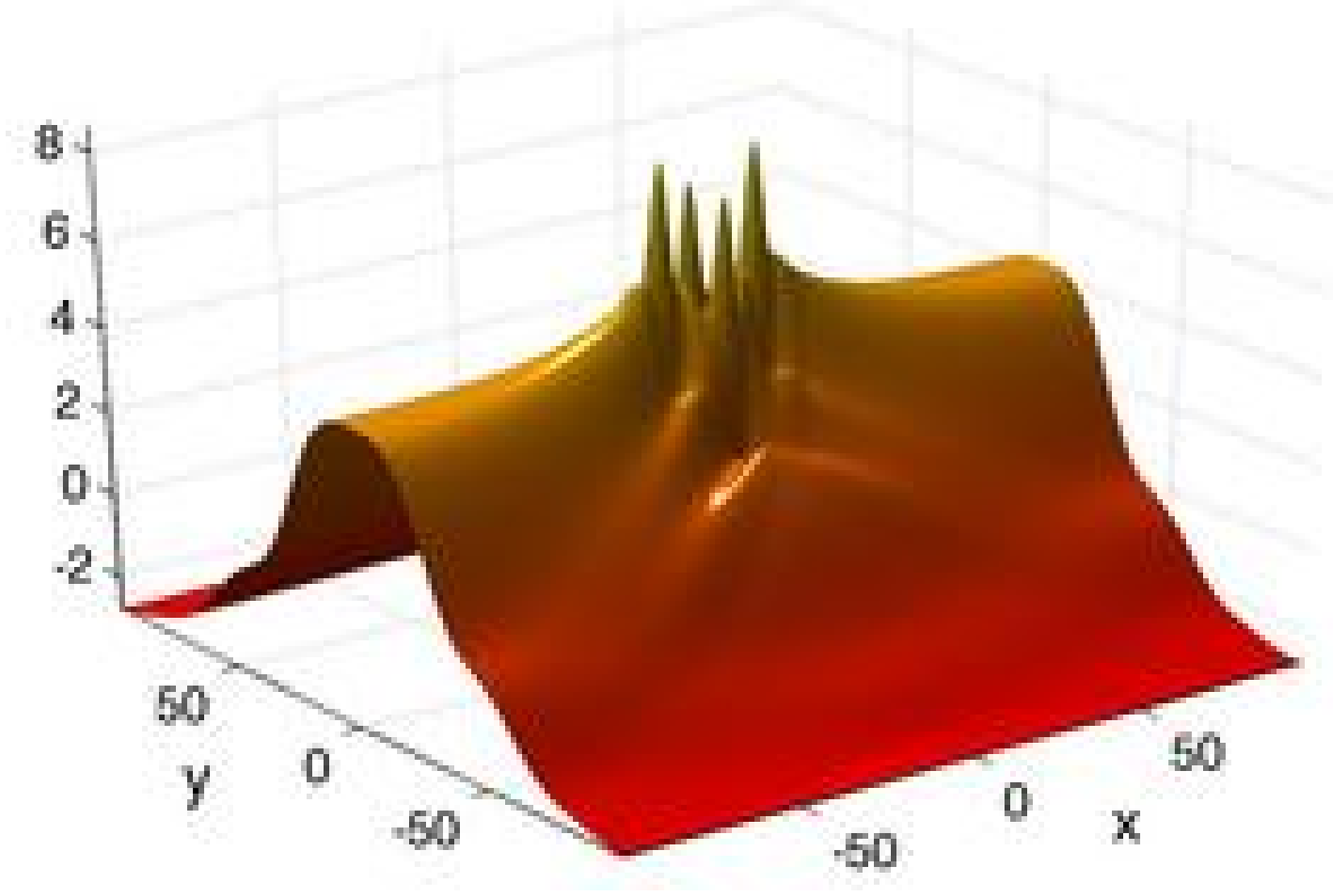}\\\includegraphics[width=4.064cm]{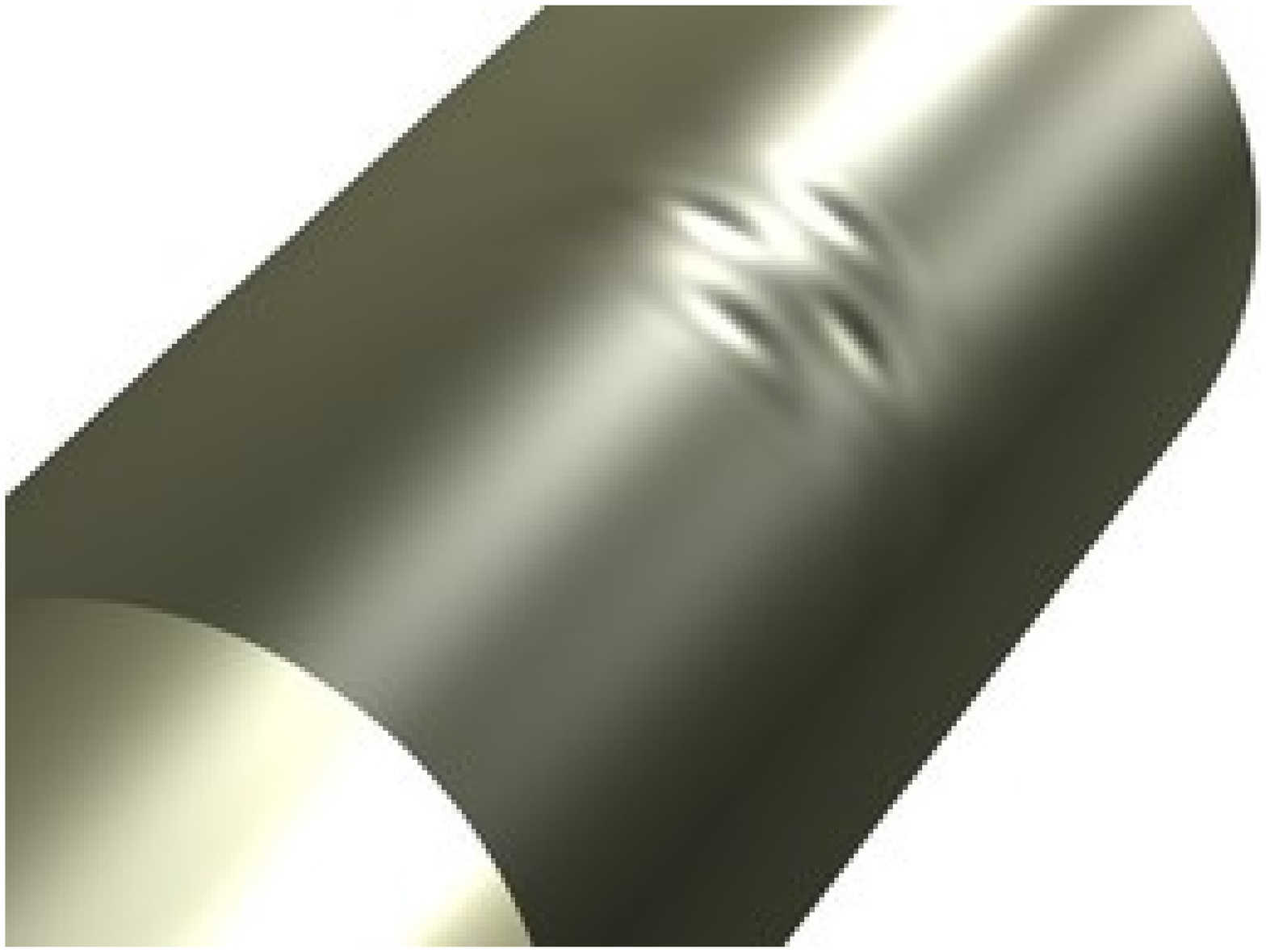}}}
\put(46,0){\frparbcenter{\includegraphics[width=4.233cm]{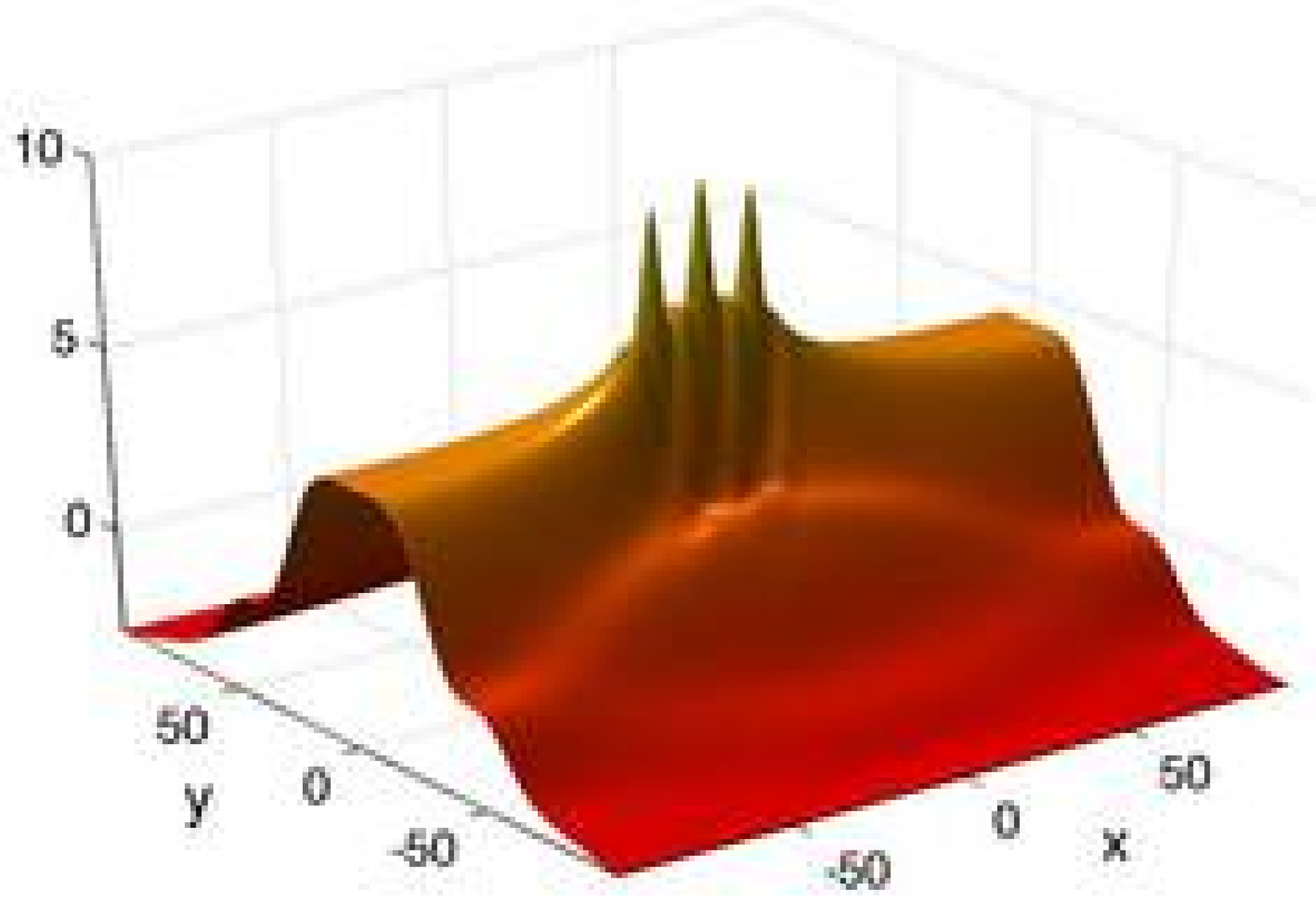}\\\includegraphics[width=4.064cm]{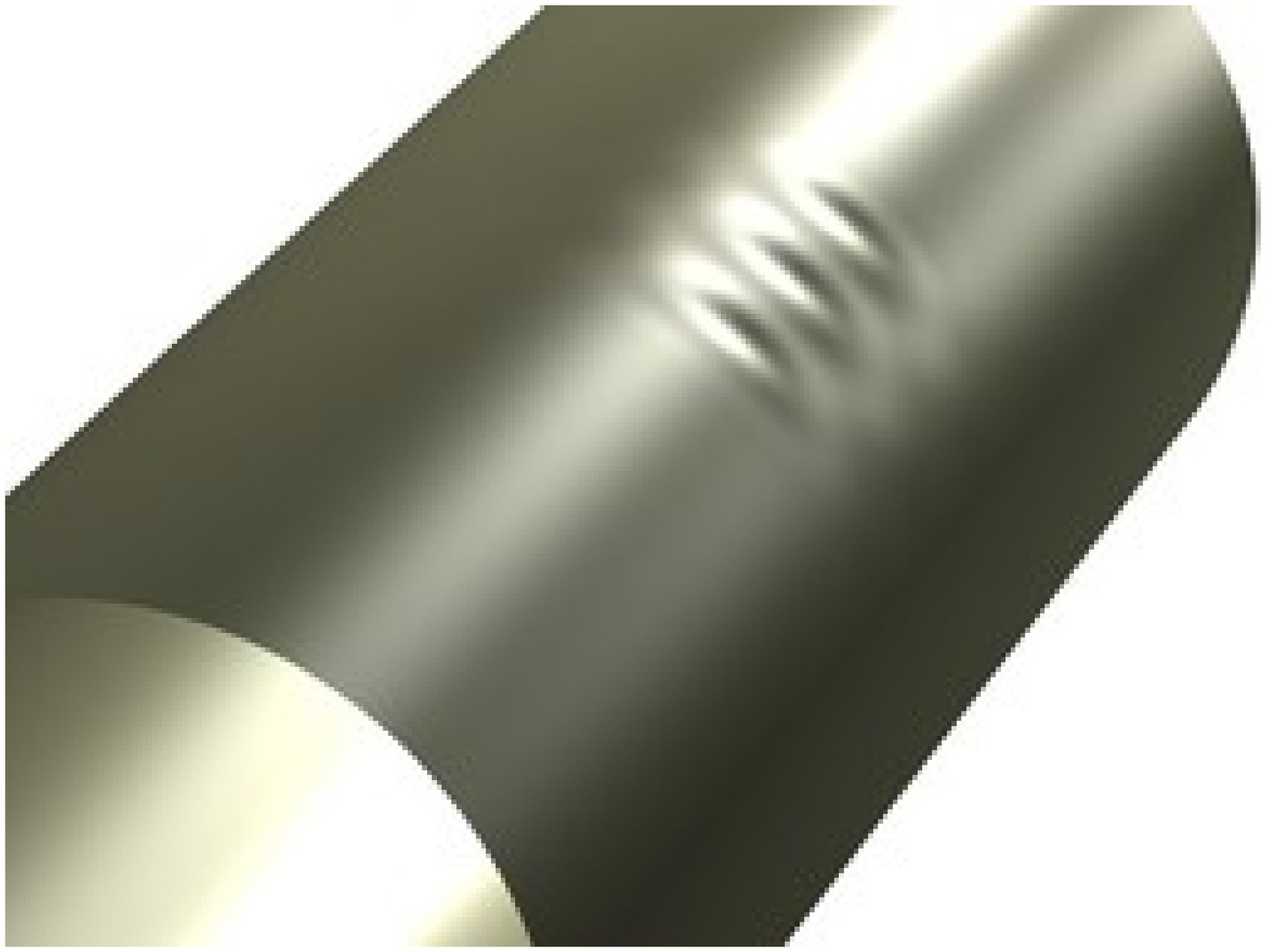}}}
\put(92,0){\frparbcenter{\includegraphics[width=4.233cm]{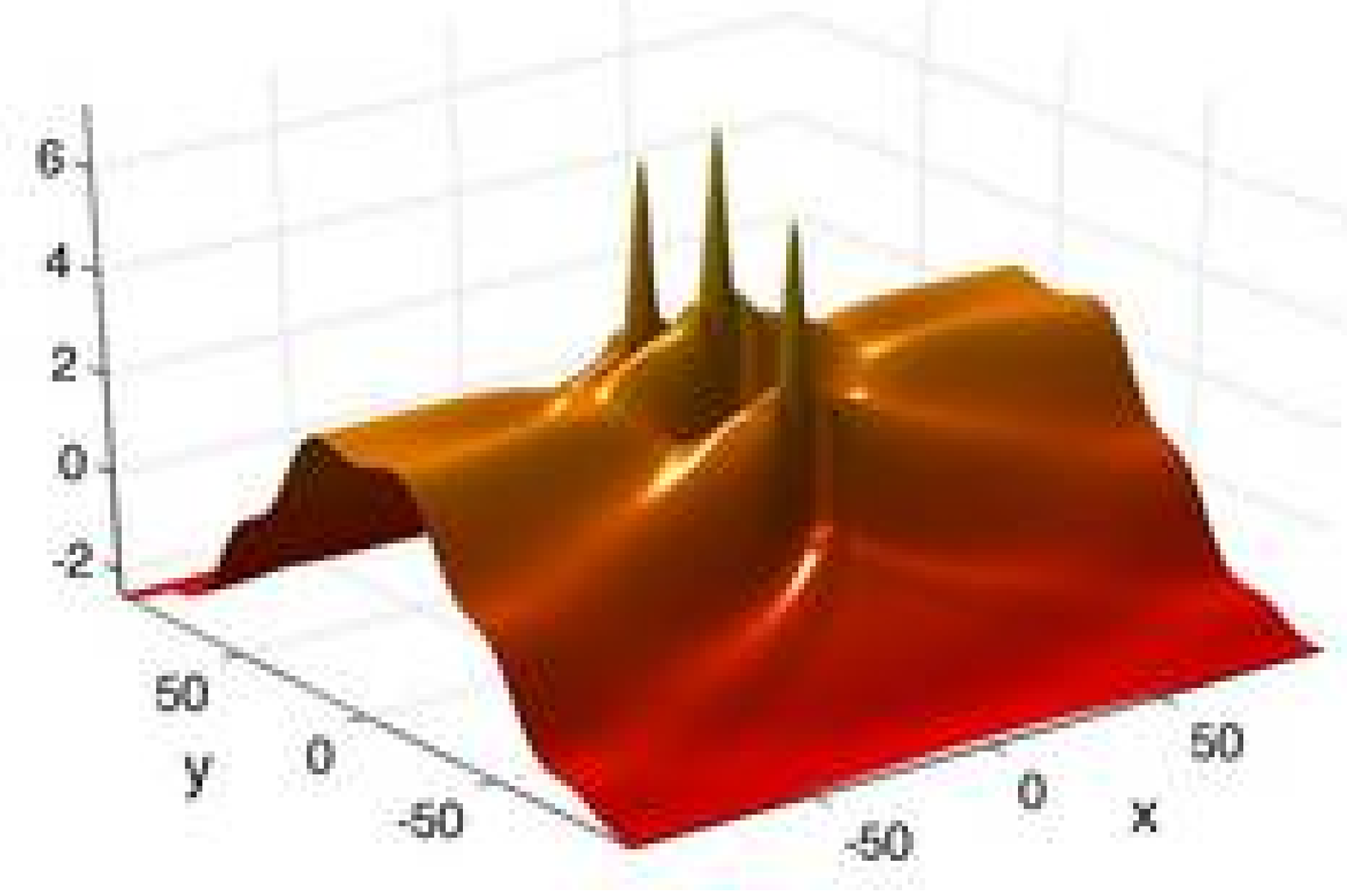}\\\includegraphics[width=4.064cm]{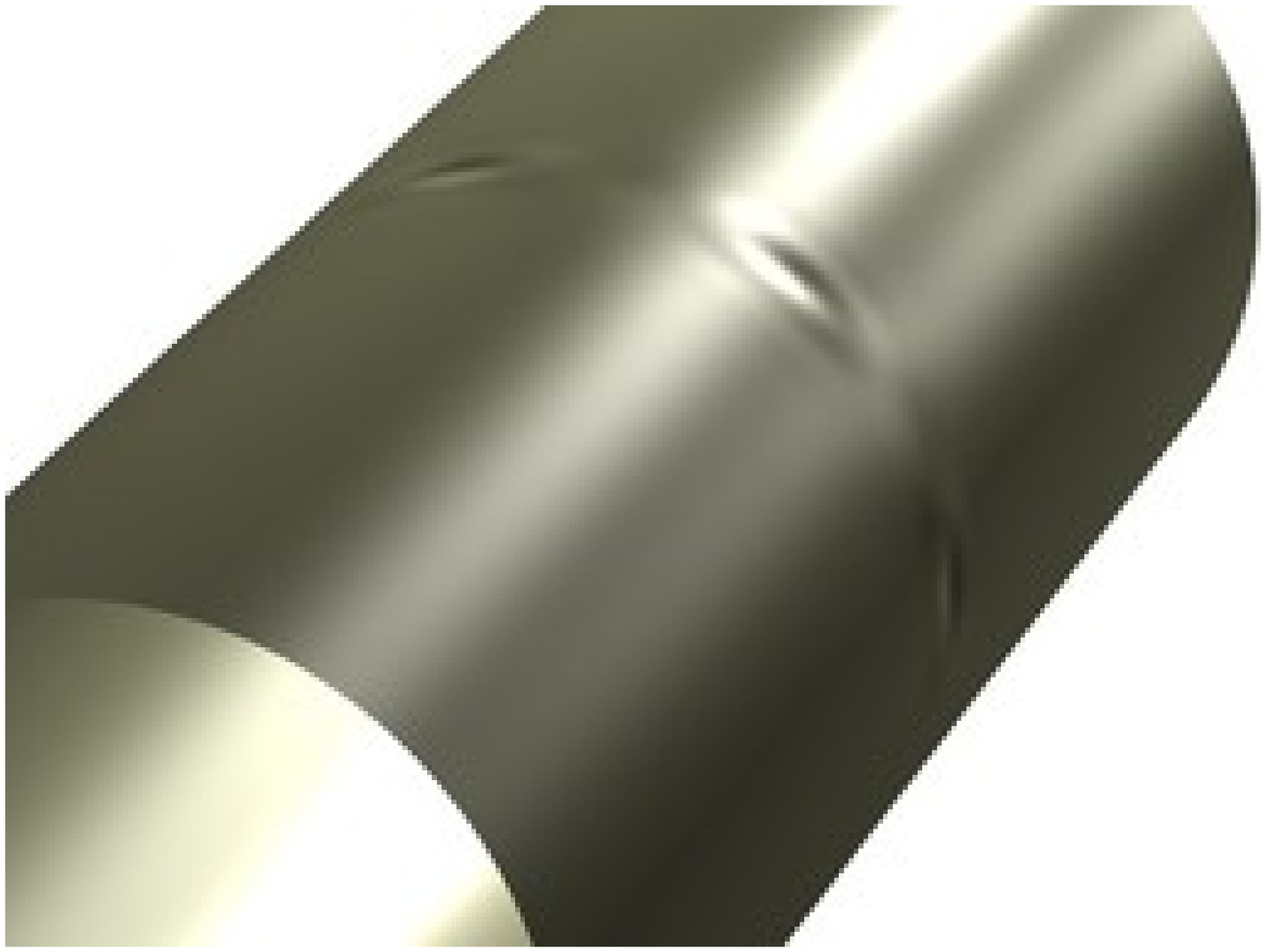}}}
\put(138,0){\frparbcenter{\includegraphics[width=4.233cm]{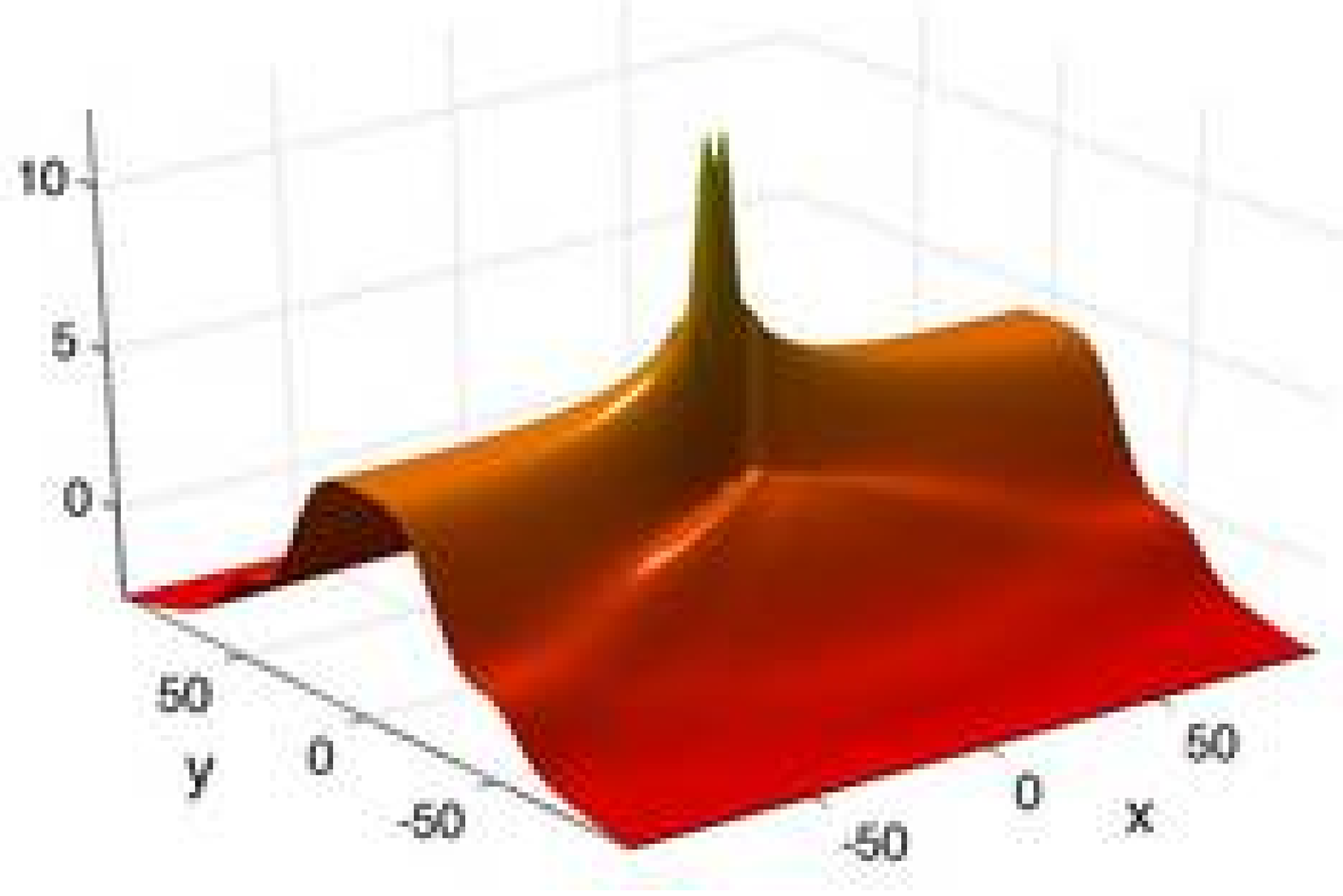}\\\includegraphics[width=4.064cm]{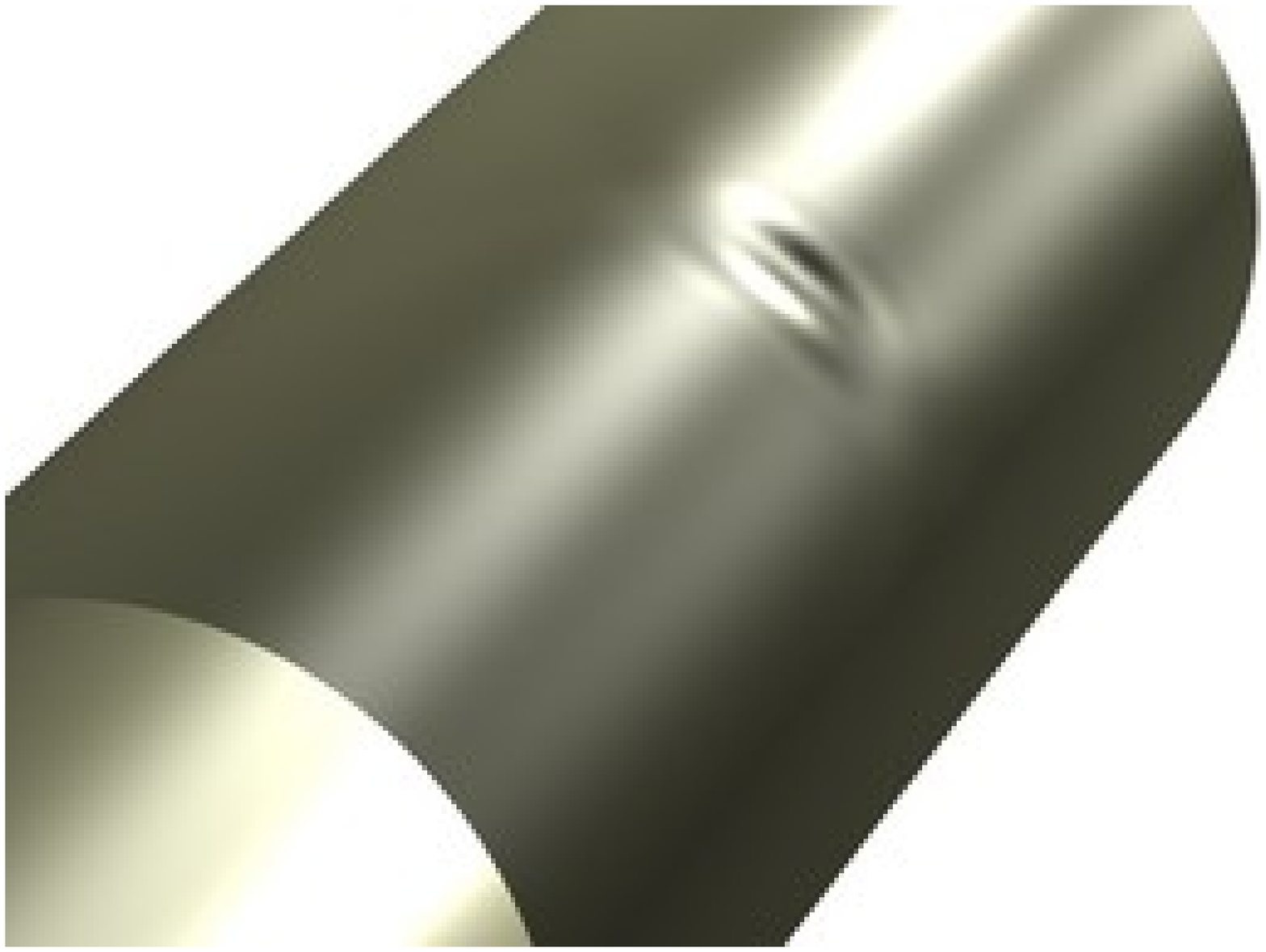}}}
\put(184,0){\frparbcenter{\includegraphics[width=4.233cm]{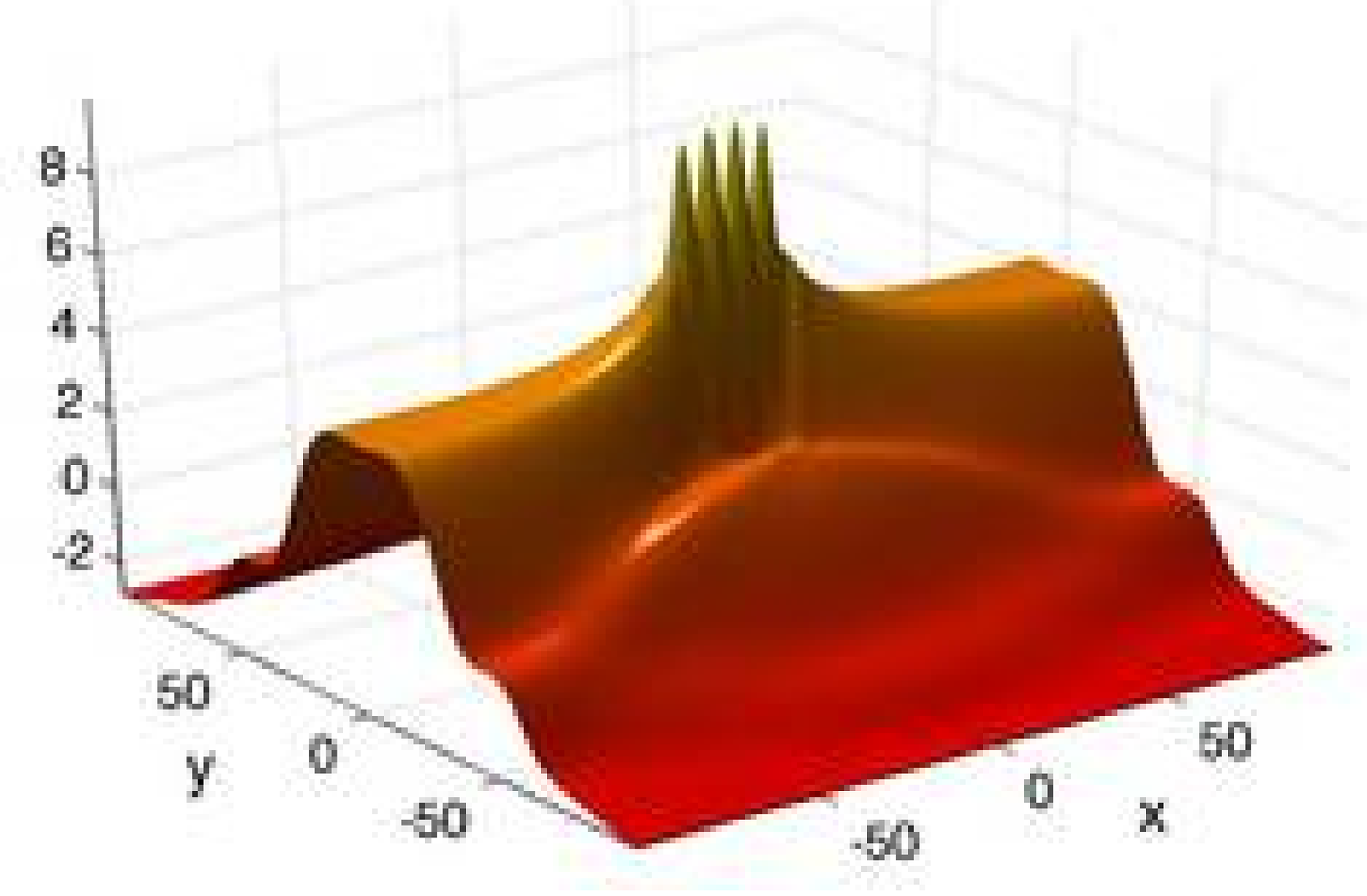}\\\includegraphics[width=4.064cm]{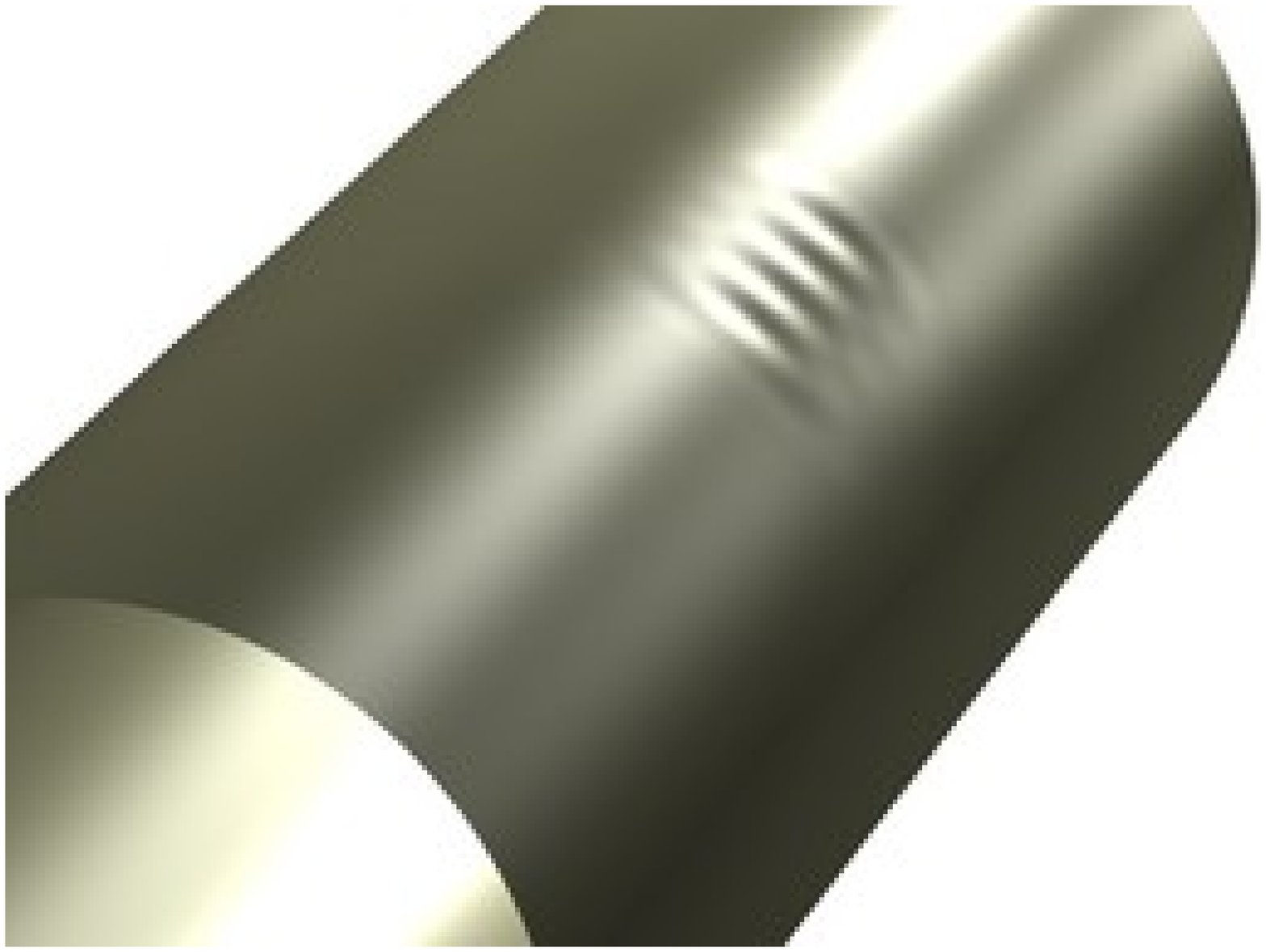}}}
\color{black}
\put(31,65){\makebox(11,5)[rb]{(3.1)}}
\put(77,65){\makebox(11,5)[rb]{(3.2)}}
\put(123,65){\makebox(11,5)[rb]{(3.3)}}
\put(169,65){\makebox(11,5)[rb]{(3.4)}}
\put(215,65){\makebox(11,5)[rb]{(3.5)}}
\put(31,1){\makebox(11,5)[rb]{(3.6)}}
\put(77,1){\makebox(11,5)[rb]{(3.7)}}
\put(123,1){\makebox(11,5)[rb]{(3.8)}}
\put(169,1){\makebox(11,5)[rb]{(3.9)}}
\put(215,1){\makebox(11,5)[rb]{(3.10)}}
\end{picture}
\end{center}
\caption{Numerical solutions found using the CMPA/Newton with $S=40$. More details are in Table~\ref{tab:sol_cmpa}.}
\label{fig:sol_cmpa1}
\end{figure}
\end{landscape}

\begin{landscape}
\begin{figure}[!ht]
\begin{center}
\setlength{\unitlength}{1mm}
\begin{picture}(227,124)
\color[rgb]{.5,.5,.5}
\put(0,64){\frparbcenter{\includegraphics[width=4.233cm]{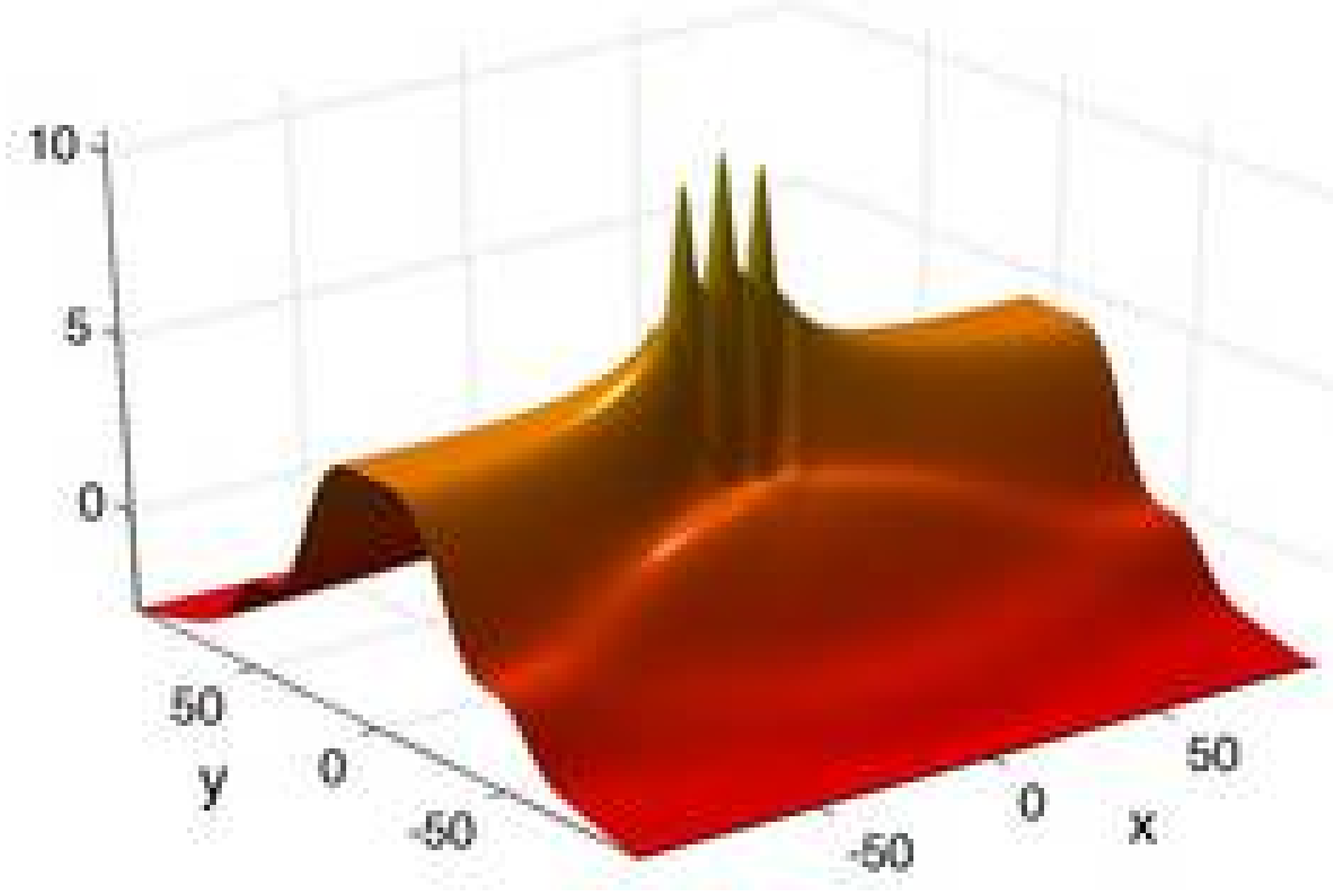}\\\includegraphics[width=4.064cm]{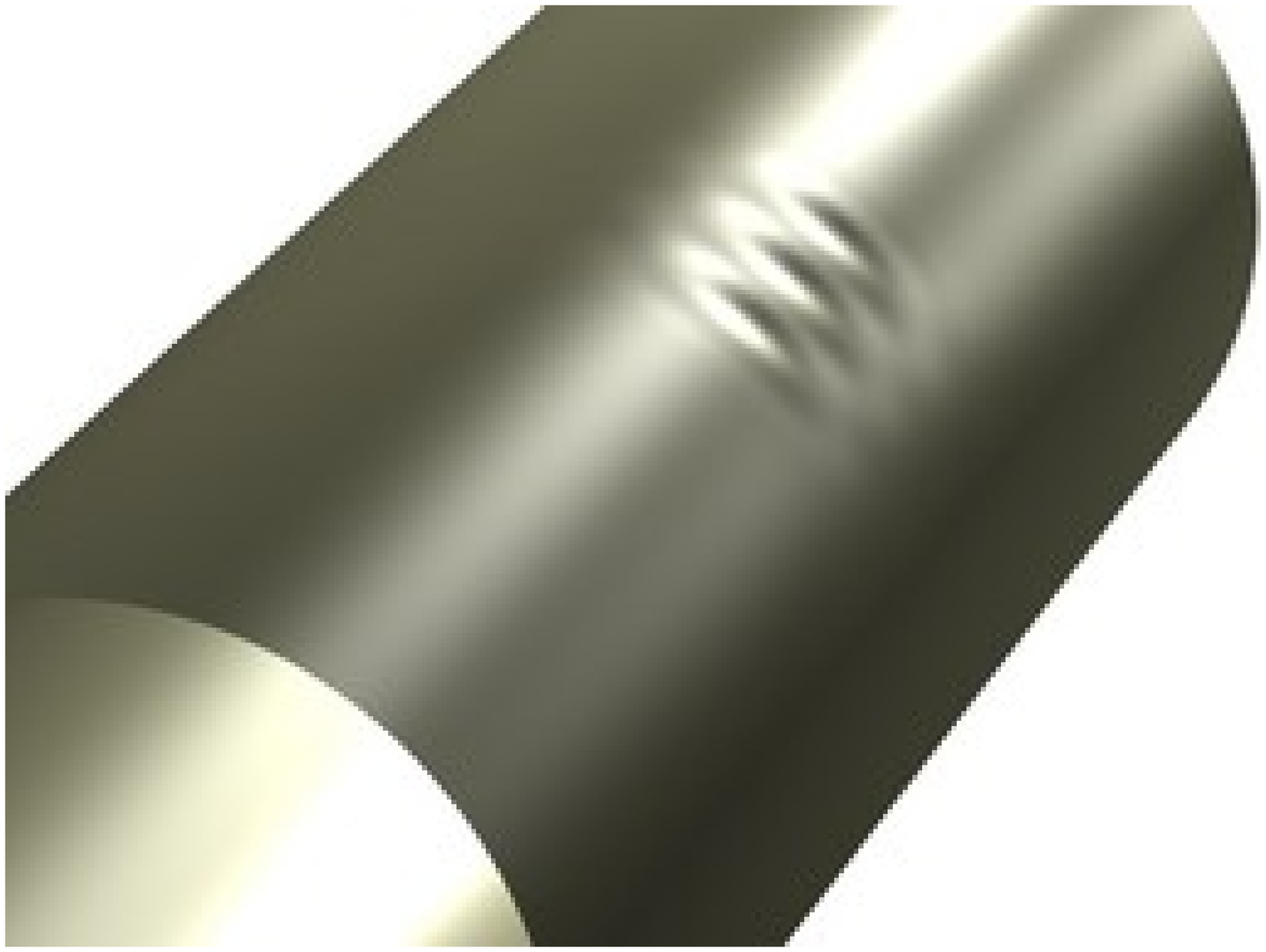}}}
\put(46,64){\frparbcenter{\includegraphics[width=4.233cm]{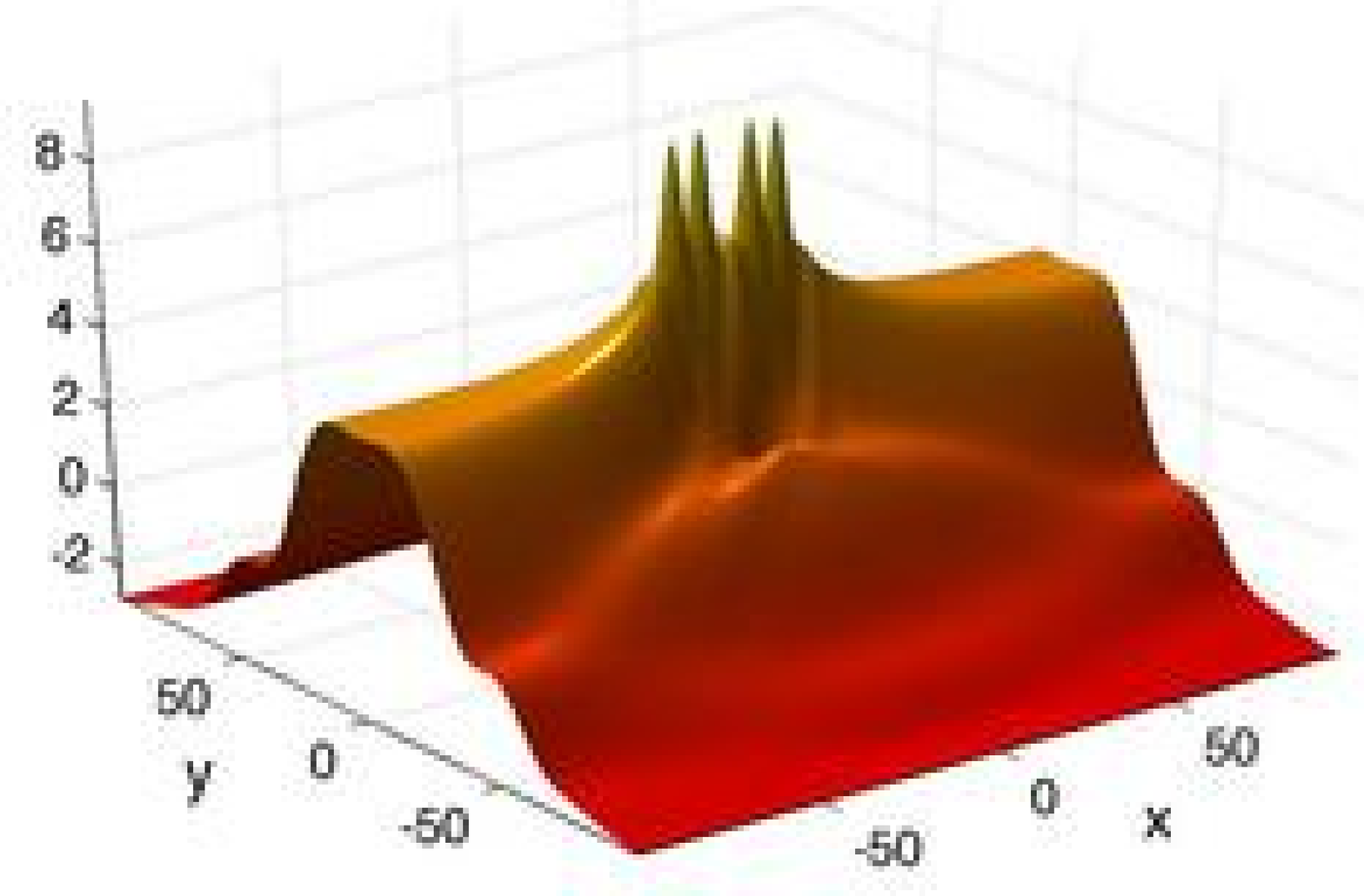}\\\includegraphics[width=4.064cm]{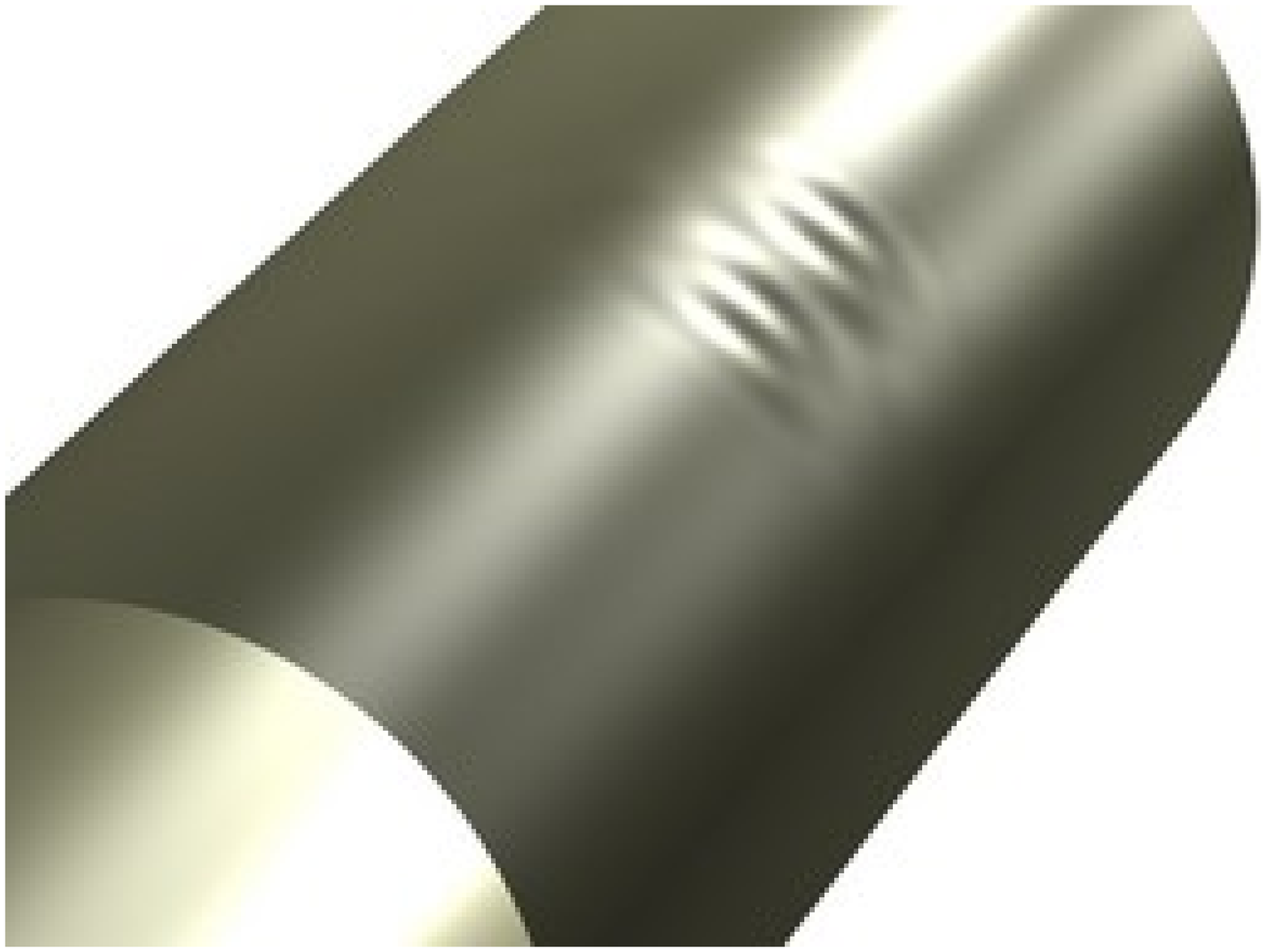}}}
\put(92,64){\frparbcenter{\includegraphics[width=4.233cm]{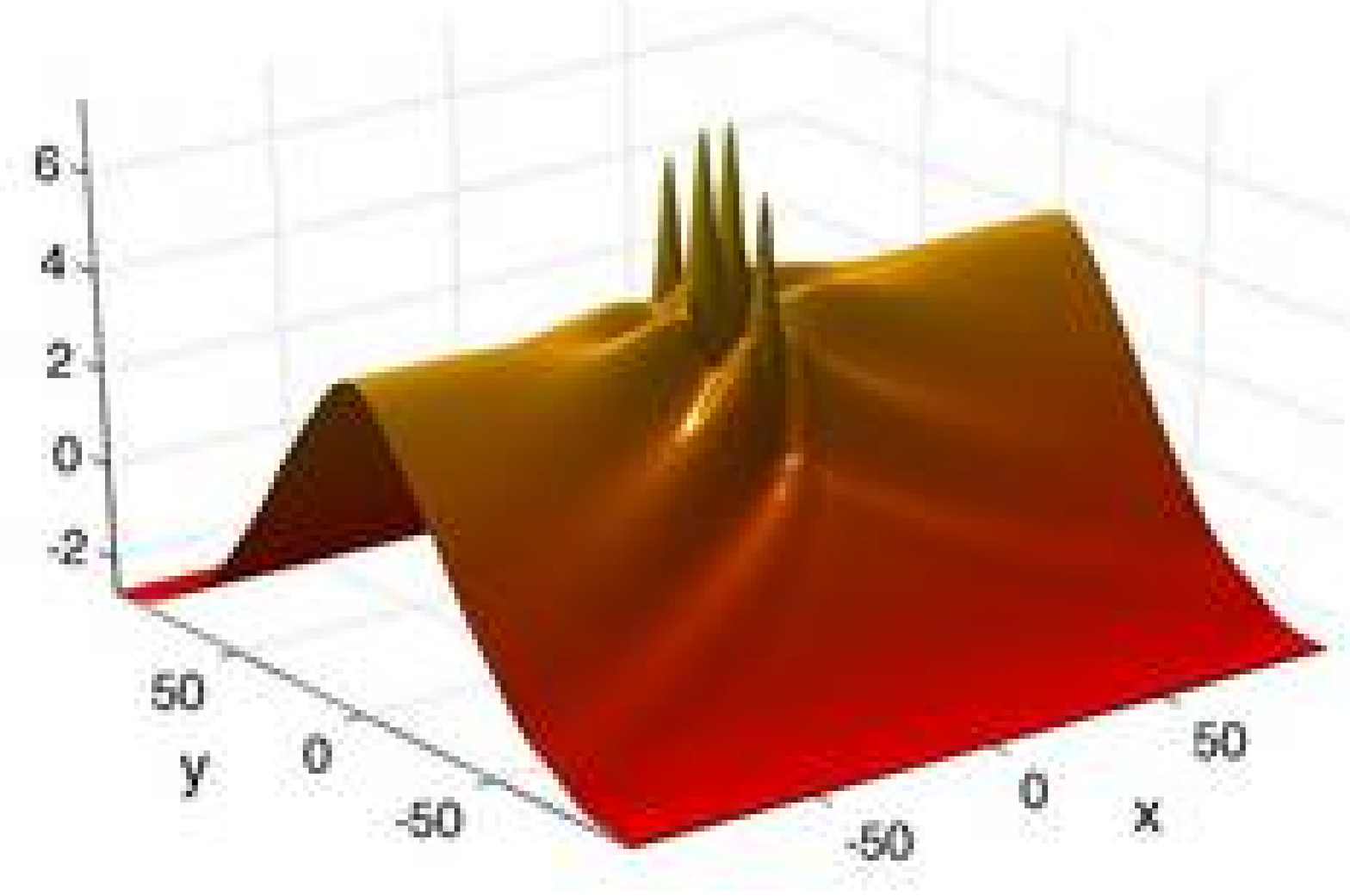}\\\includegraphics[width=4.064cm]{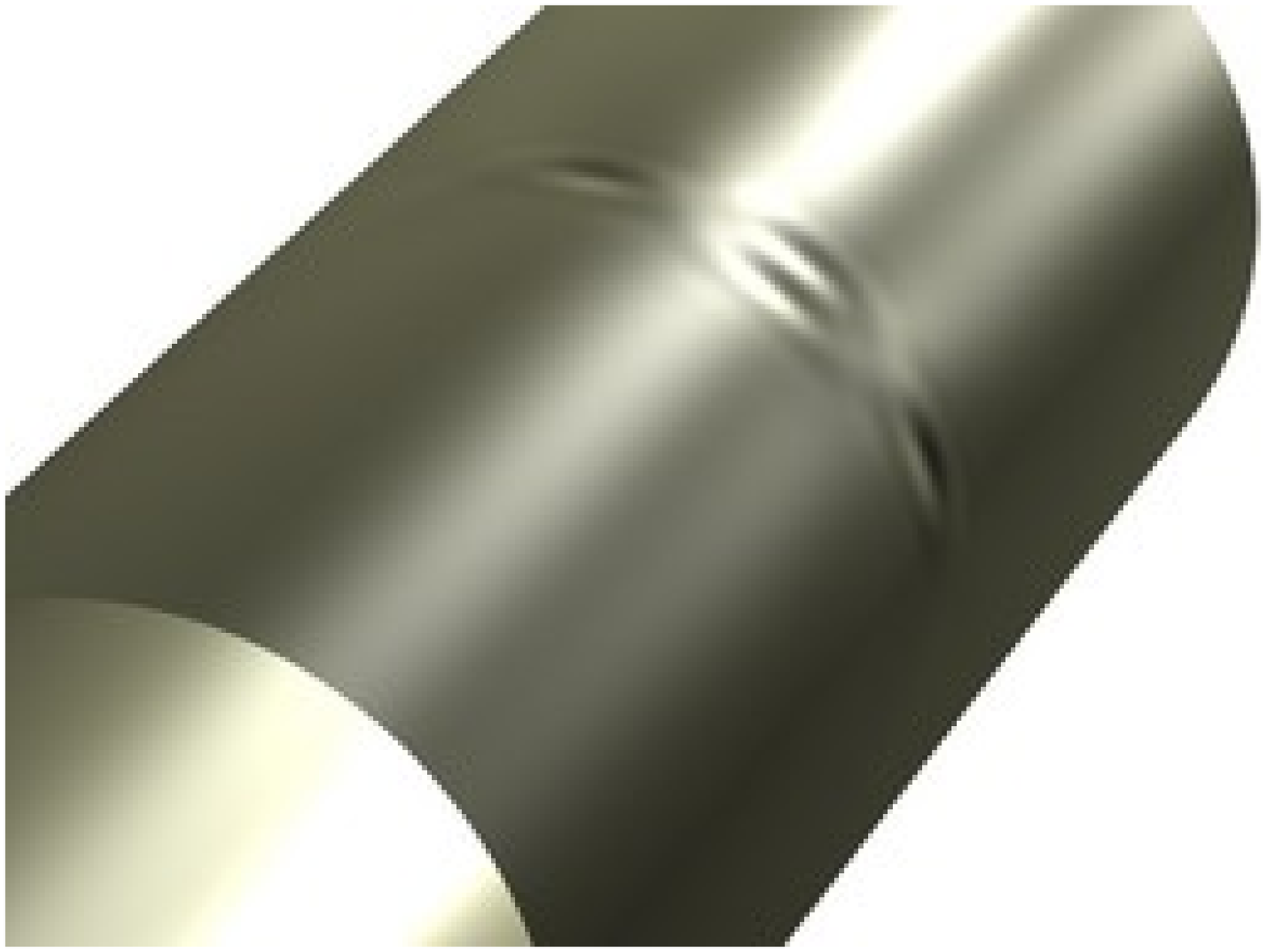}}}
\put(138,64){\frparbcenter{\includegraphics[width=4.233cm]{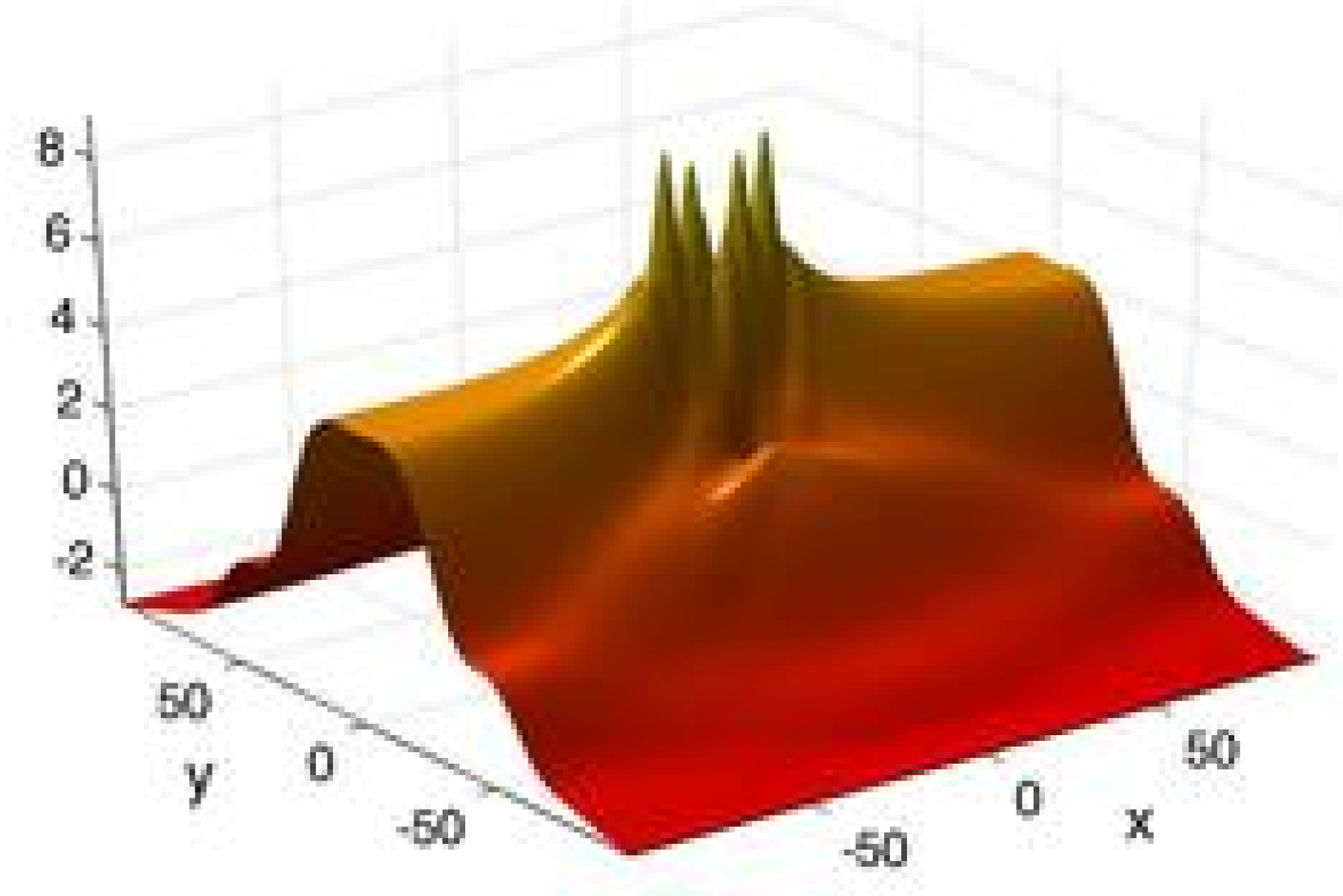}\\\includegraphics[width=4.064cm]{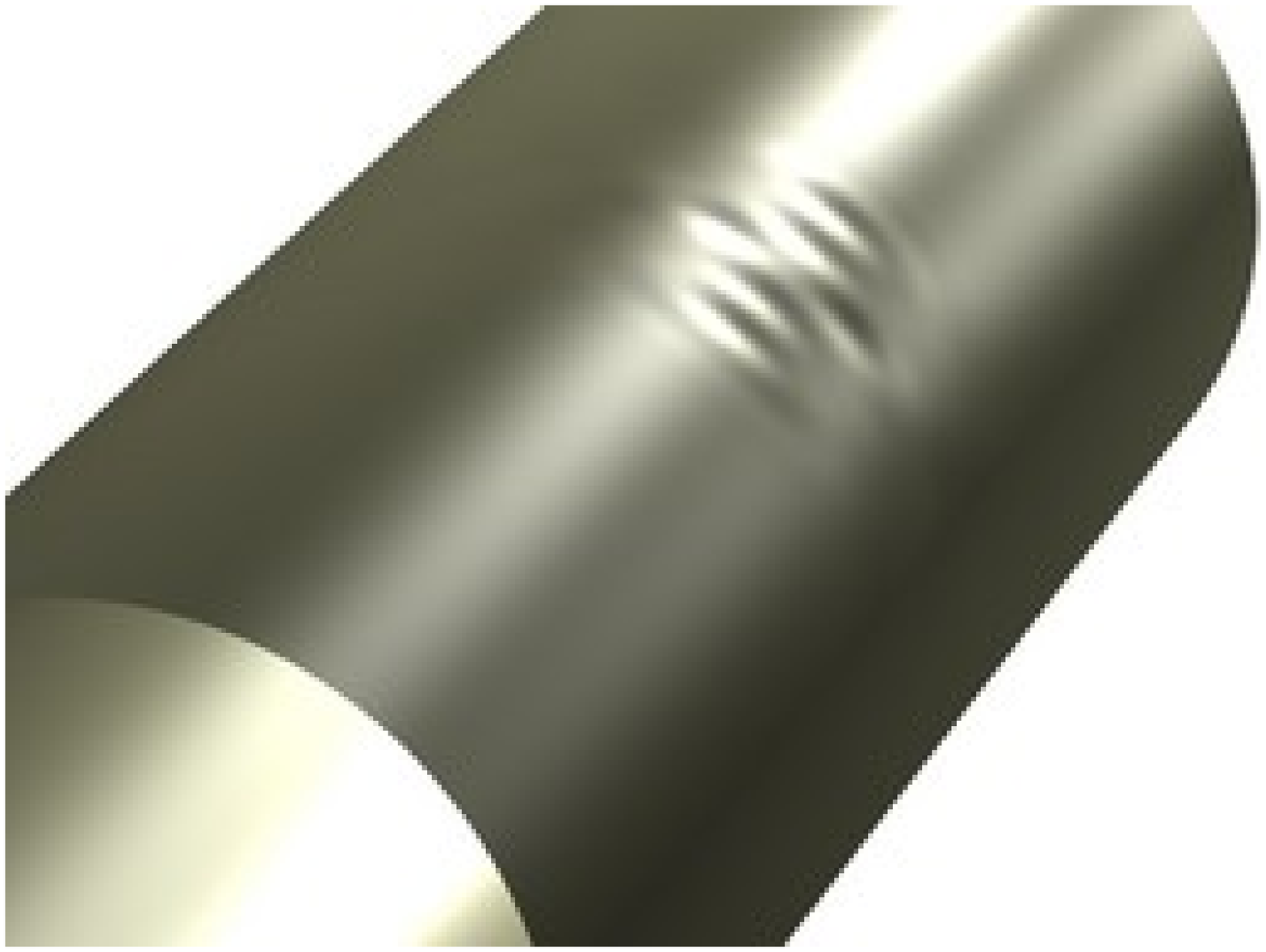}}}
\put(184,64){\frparbcenter{\includegraphics[width=4.233cm]{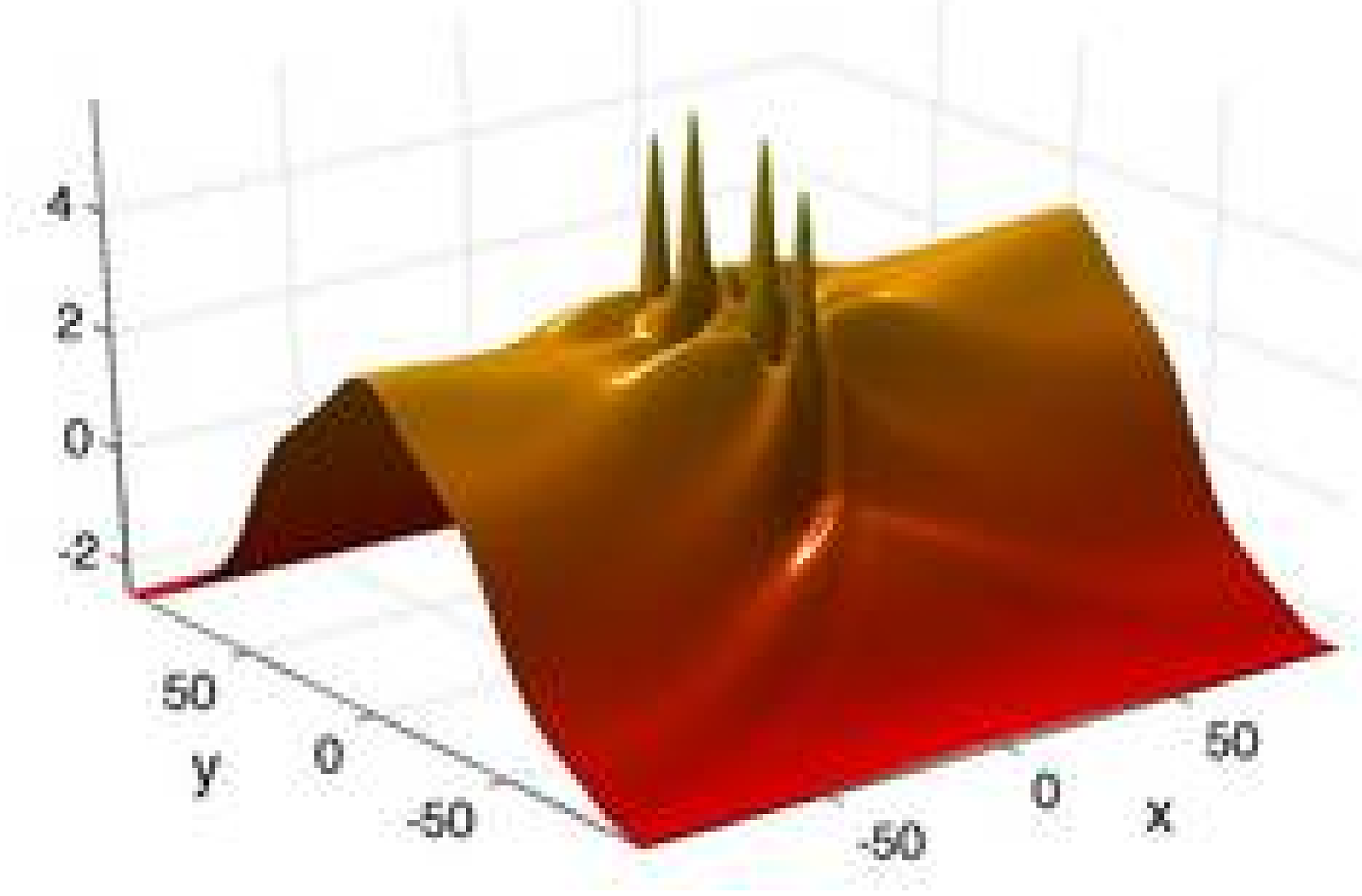}\\\includegraphics[width=4.064cm]{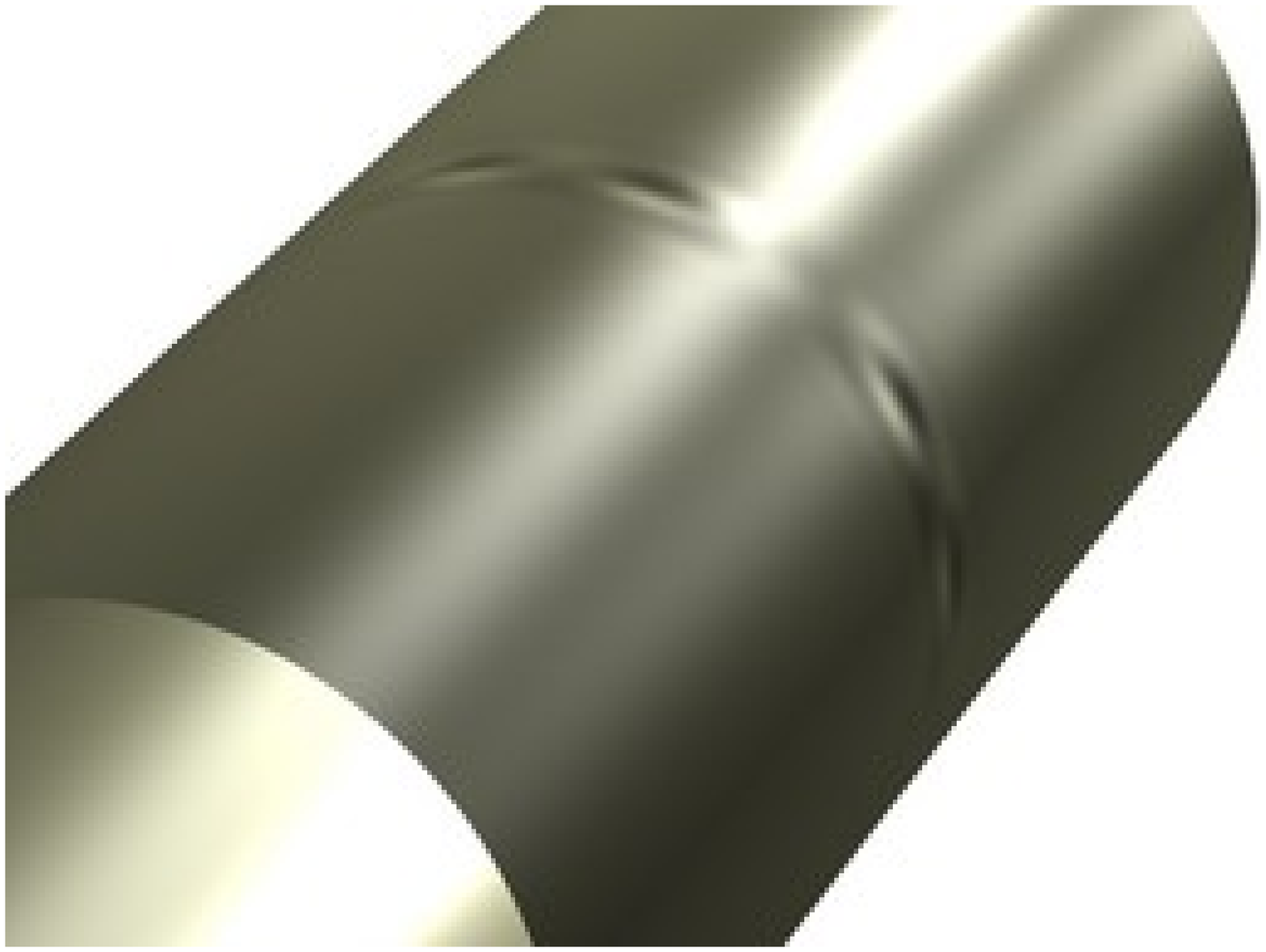}}}
\put(0,0){\frparbcenter{\includegraphics[width=4.233cm]{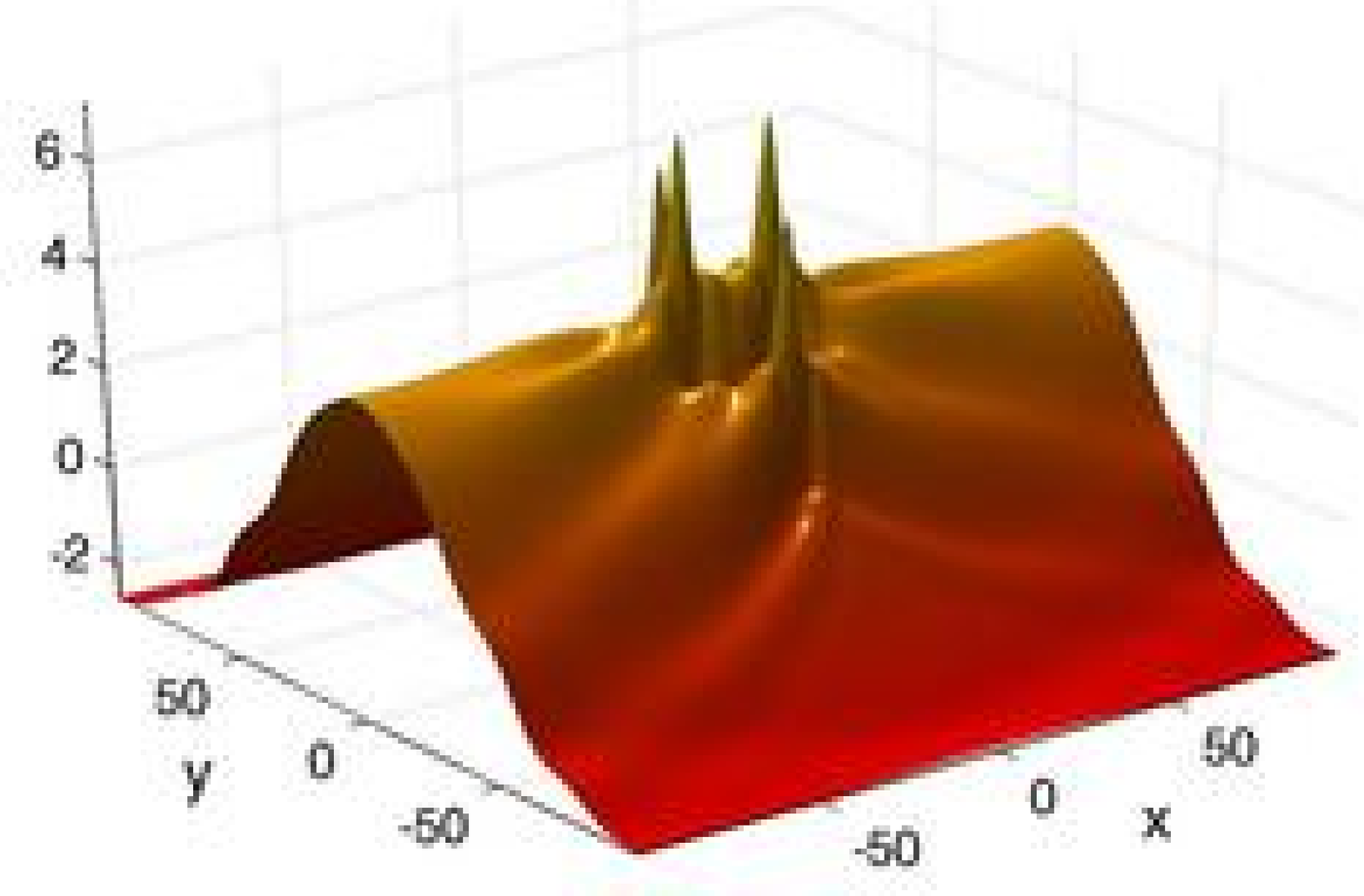}\\\includegraphics[width=4.064cm]{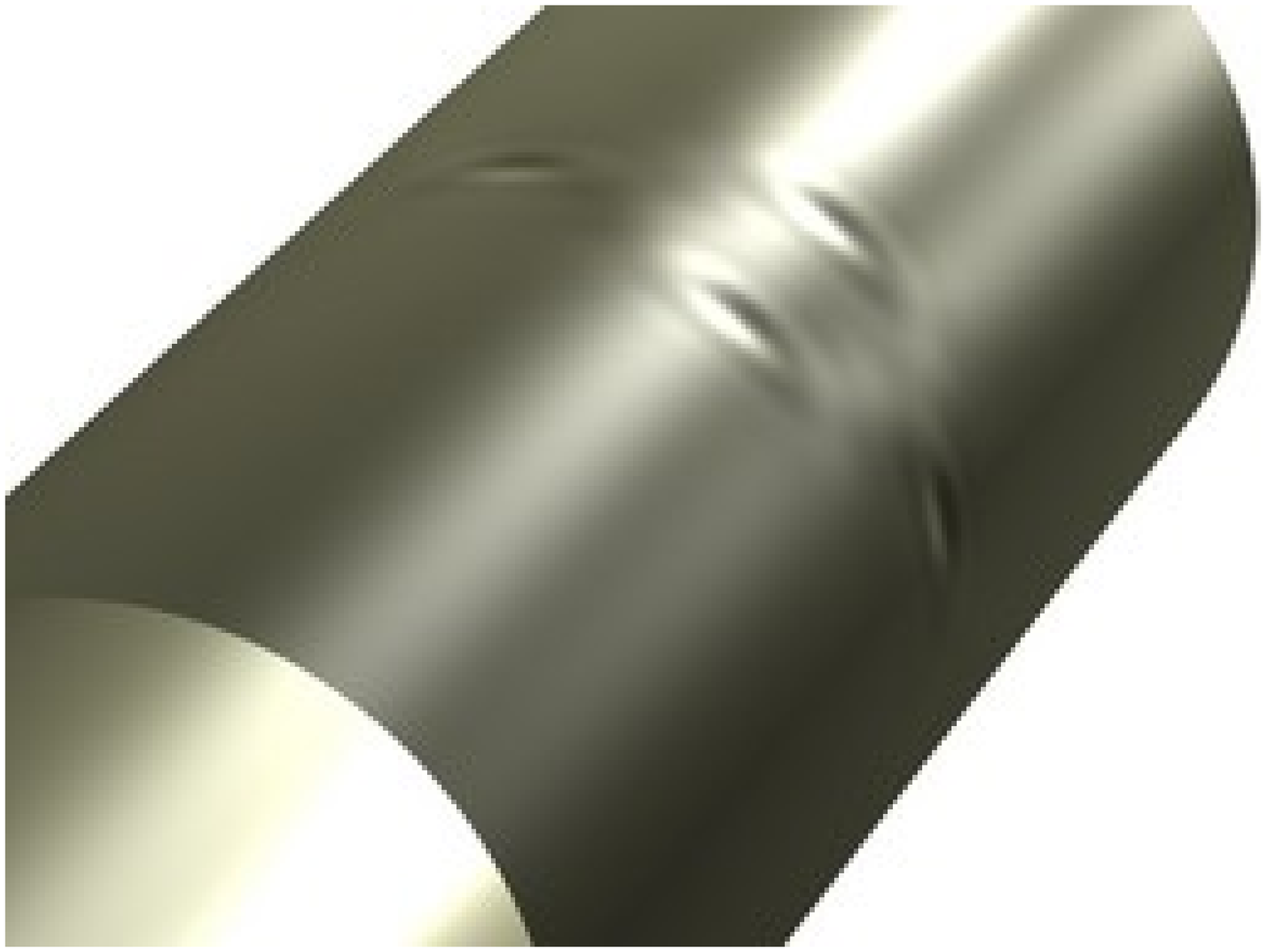}}}
\put(46,0){\frparbcenter{\includegraphics[width=4.233cm]{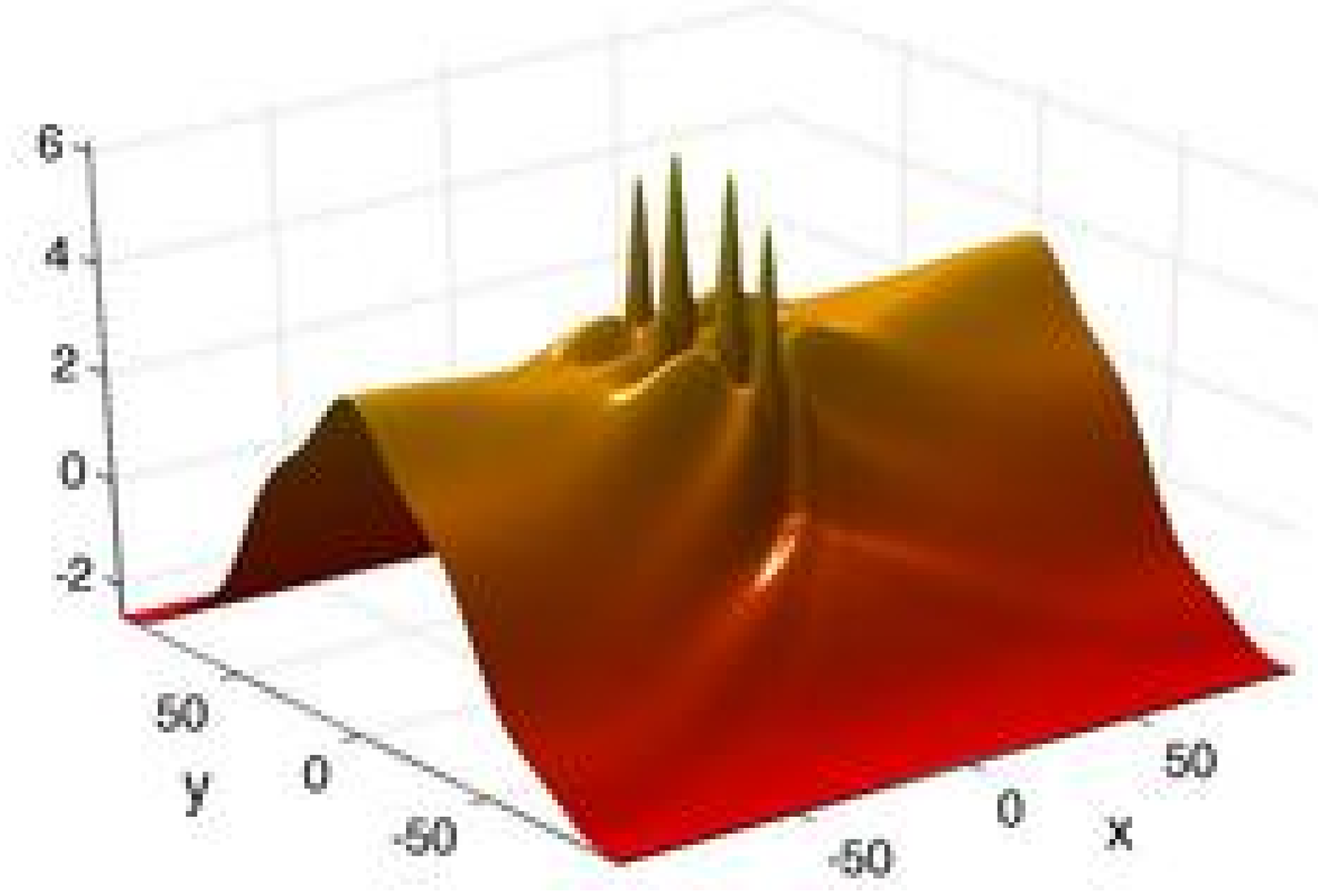}\\\includegraphics[width=4.064cm]{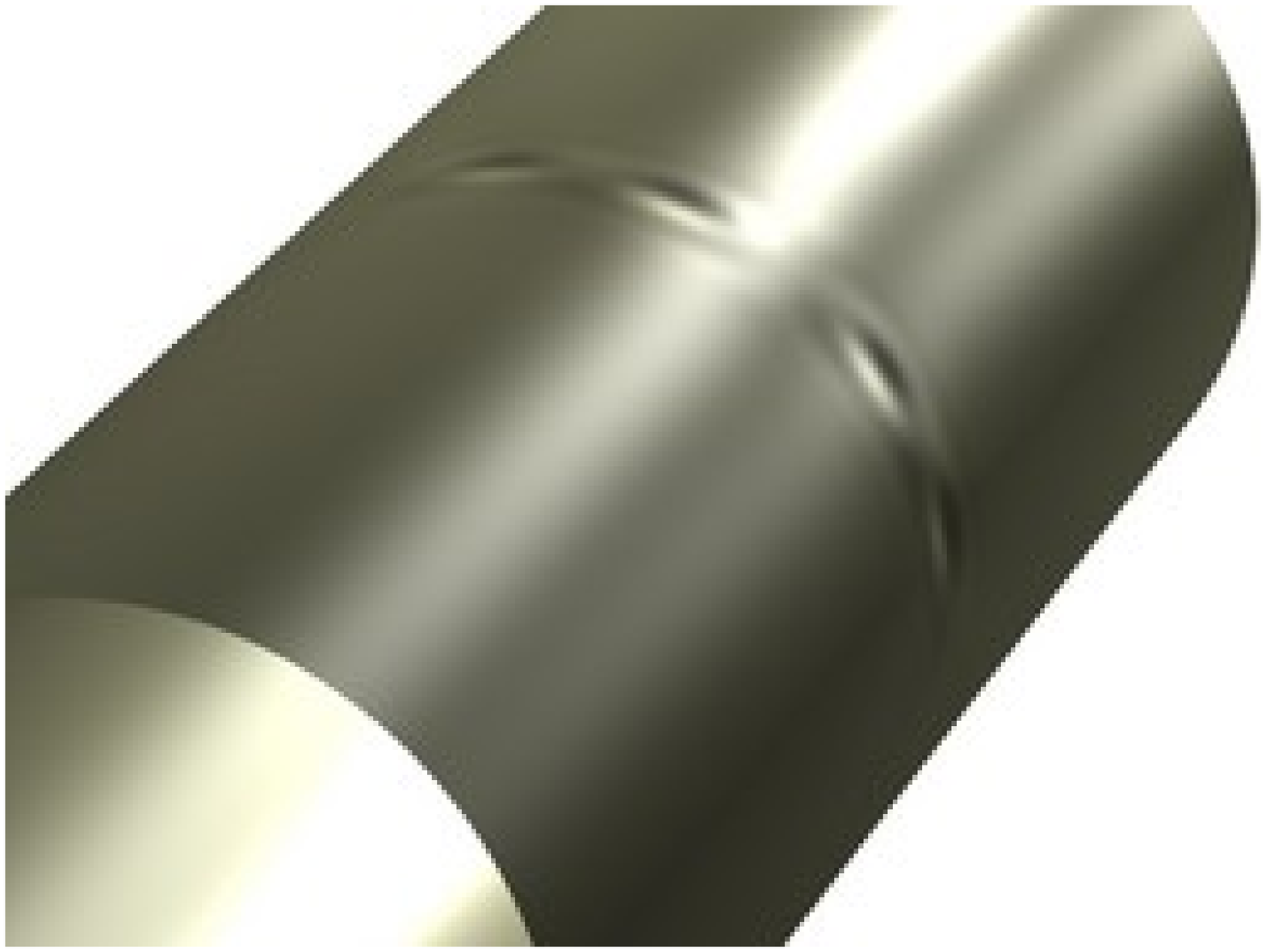}}}
\put(92,0){\frparbcenter{\includegraphics[width=4.233cm]{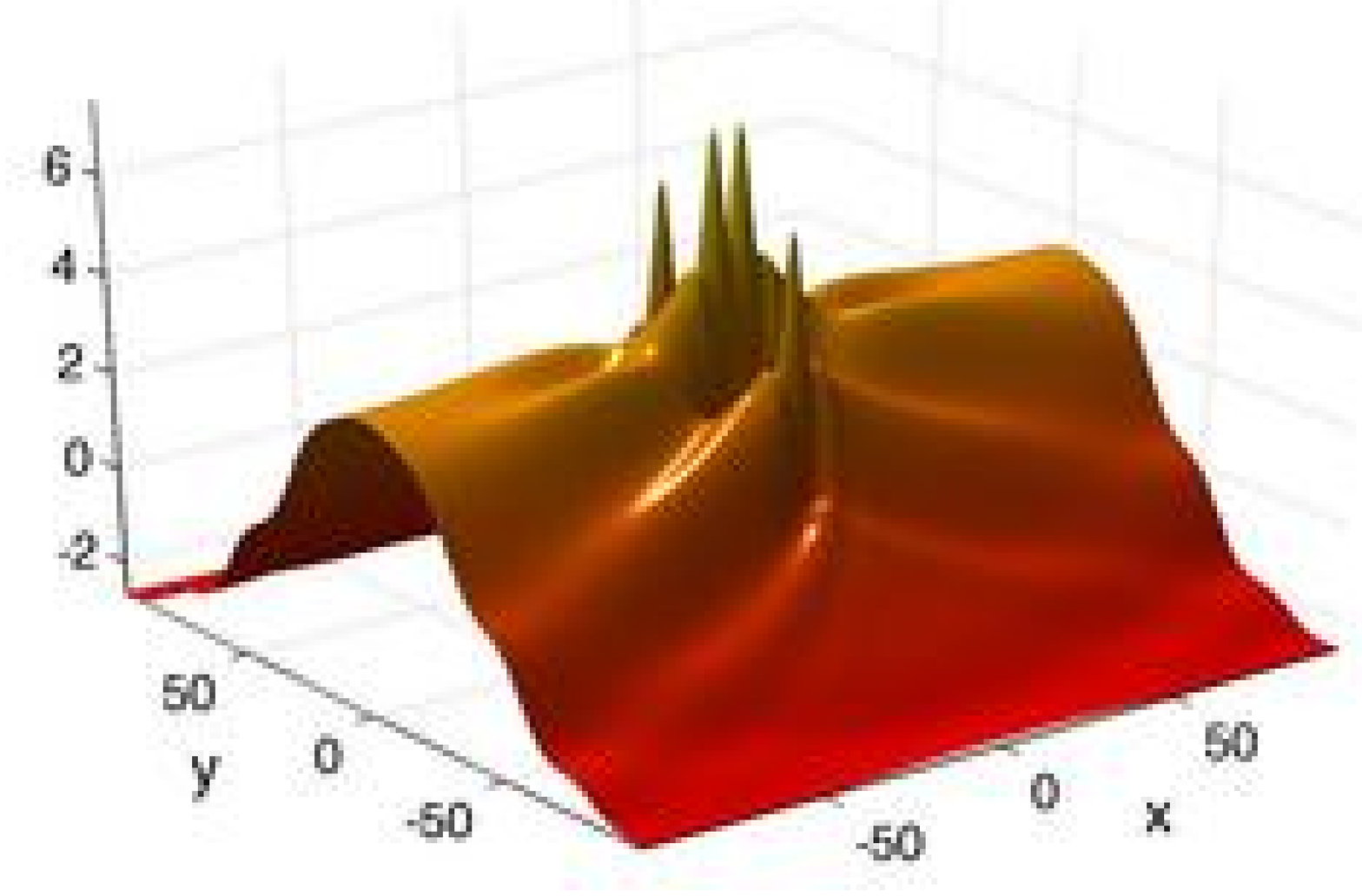}\\\includegraphics[width=4.064cm]{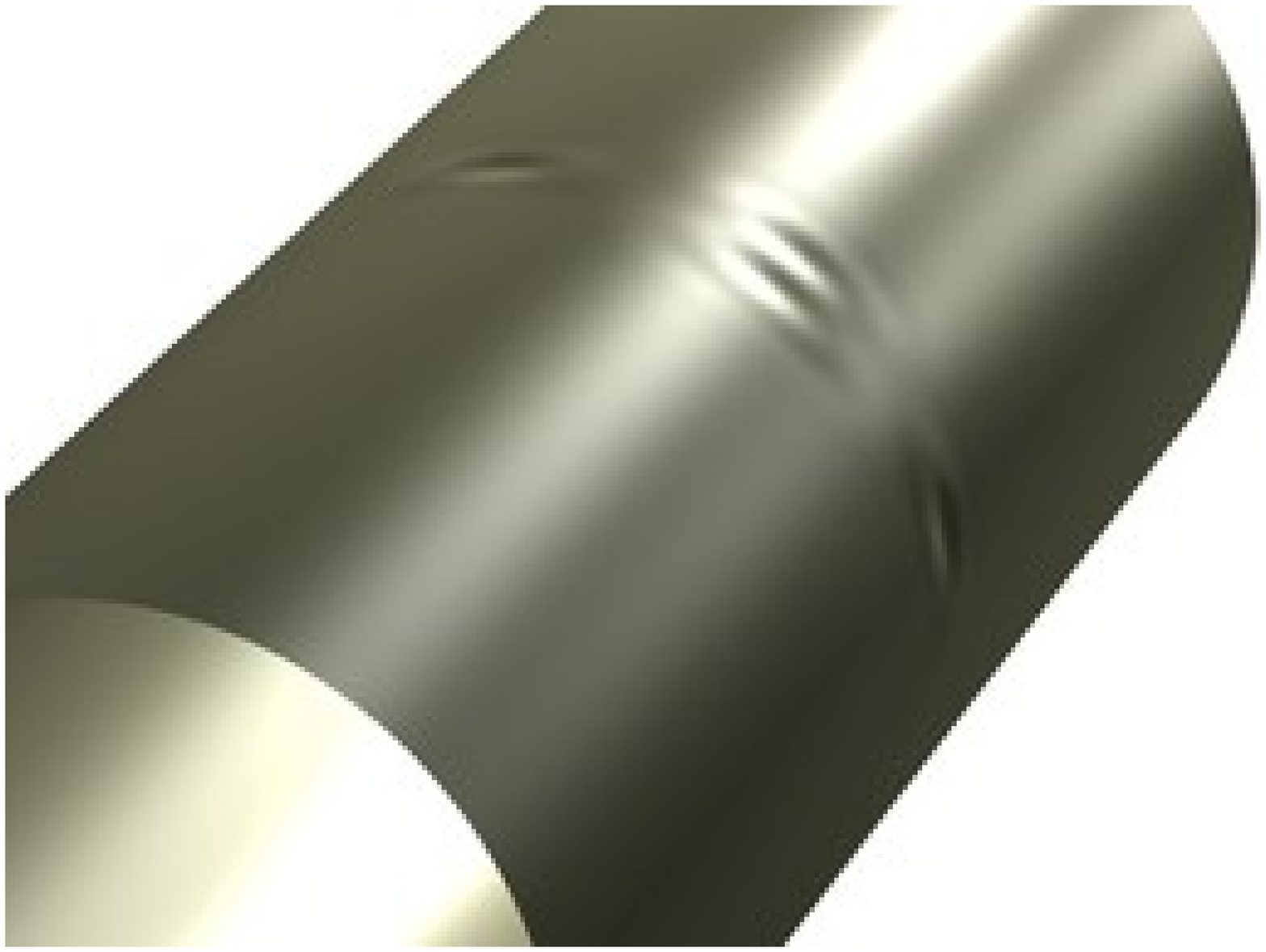}}}
\put(138,0){\frparbcenter{\includegraphics[width=4.233cm]{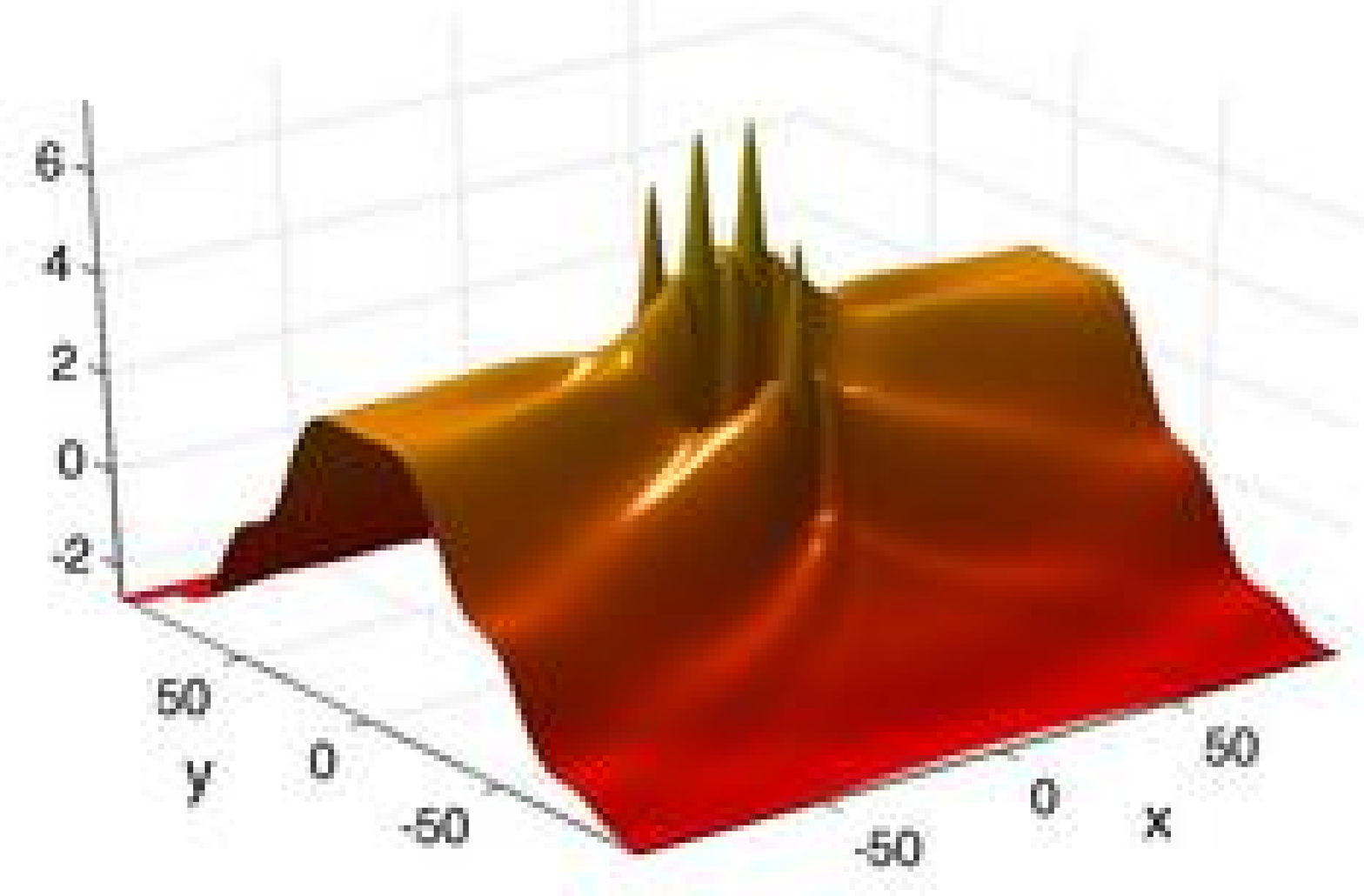}\\\includegraphics[width=4.064cm]{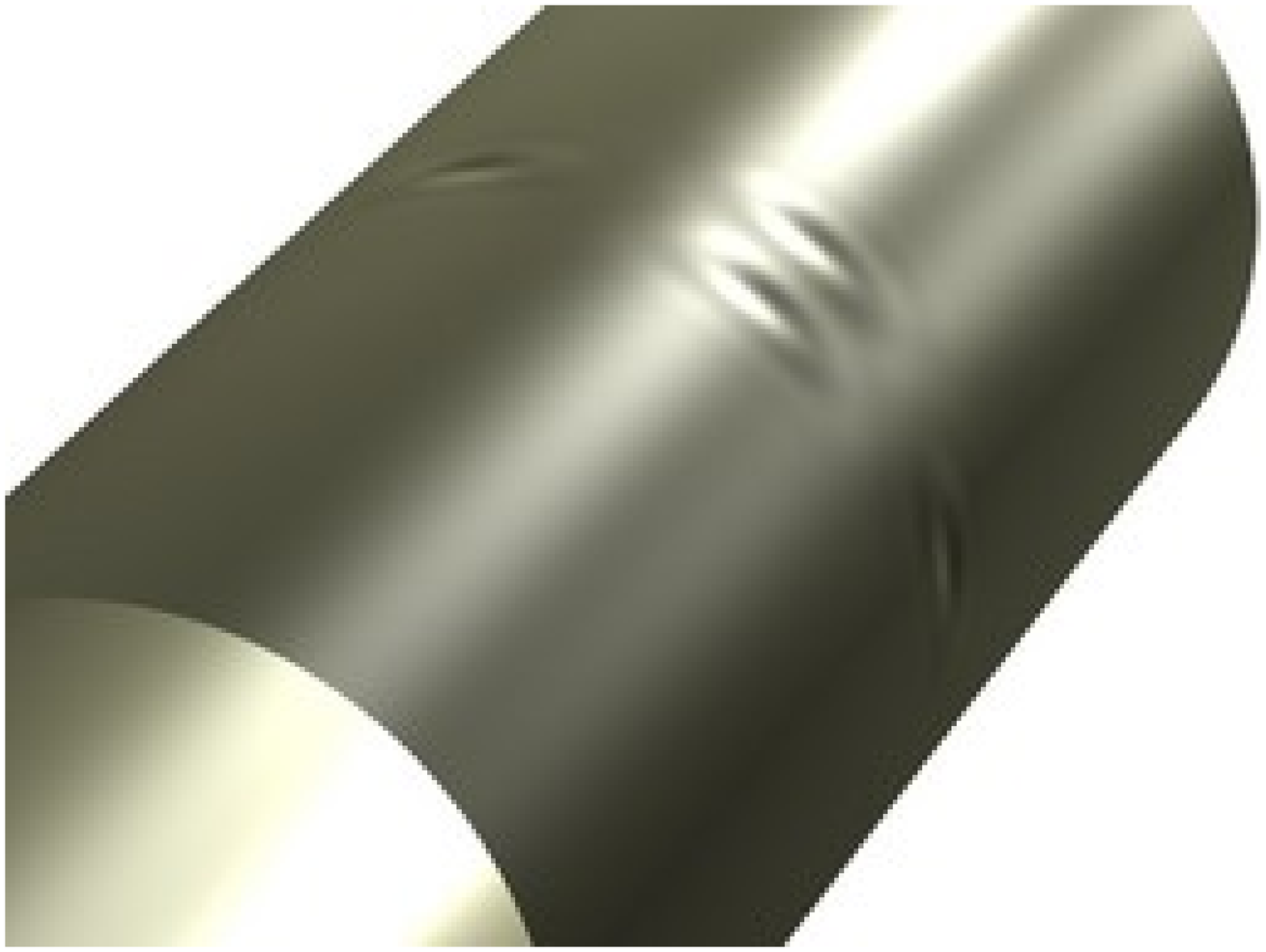}}}
\color{black}
\put(31,65){\makebox(11,5)[rb]{(3.11)}}
\put(77,65){\makebox(11,5)[rb]{(3.12)}}
\put(123,65){\makebox(11,5)[rb]{(3.13)}}
\put(169,65){\makebox(11,5)[rb]{(3.14)}}
\put(215,65){\makebox(11,5)[rb]{(3.15)}}
\put(31,1){\makebox(11,5)[rb]{(3.16)}}
\put(77,1){\makebox(11,5)[rb]{(3.17)}}
\put(123,1){\makebox(11,5)[rb]{(3.18)}}
\put(169,1){\makebox(11,5)[rb]{(3.19)}}
\end{picture}
\end{center}
\caption{Numerical solutions found using the CMPA/Newton with $S=40$. More details are in Table~\ref{tab:sol_cmpa}.}
\label{fig:sol_cmpa2}
\end{figure}
\end{landscape}

We then fixed $S=40$  and applied the CMPA to obtain constrained local
mountain passes of $E$ described in~Table~\ref{tab:sol_cmpa}. They are
again ordered according to the increasing value of stored energy $E$.
Their graphs and renderings on a cylinder are shown
in~Figs.~\ref{fig:sol_cmpa1} and~\ref{fig:sol_cmpa2}. As endpoints
$w_1$, $w_2$ of the path in the CMPA we used the constrained local
minimizers 1.1--1.7.

There are 21 possible pairs $(w_1,w_2)$ to be used but only 19
solutions in Table~\ref{tab:sol_cmpa}. The algorithm did not converge
for the following three pairs: (1.1, 1.3), (1.5, 1.6), and (1.4, 1.6),
most likely due to the complicated nature of the energy landscape
between these endpoints. On the other hand, two choices of pairs
denoted by $\ast$ and $\dagger$ in the Table yielded two solutions
each. When the path is deformed it sometimes comes close to another
critical point of $E$ which is not a constrained mountain pass. In
that case the algorithm slows down and one can apply Newton's method
to such a point. It is a matter of luck whether Newton's method
converges. The CMPA then runs further and might converge to another point,
this time a constrained mountain-pass point. And finally, two choices
of $(w_1,w_2)$ yielded the same solution 3.3.

\begin{table}[htbp]
  \centering
  \begin{tabular}{|c||c|c|c|c|c|}
    \hline
    {\bf CMPA} & $\l$ & $S$ & $E$ & $F_\l$ & end points \tabularnewline
    \hline
     3.1 & 1.310815 & 40 & 63.98996 & 11.55737 & 1.2, 1.4
     \tabularnewline 
    \hline
     3.2 & 1.332112 & 40 & 64.38609 & 11.10161 & 1.3, 1.6
     \tabularnewline 
    \hline
     3.3 & 1.447626 & 40 & 66.49032 & 8.585294 & 1.2, 1.3 or 1.3, 1.4
     \tabularnewline 
    \hline
     3.4 & 1.440841 & 40 & 66.72079 & 9.087129 & 1.1, 1.4
     \tabularnewline 
    \hline
     3.5 & 1.447594 & 40 & 66.97057 & 9.066810 & 1.1, 1.6
     \tabularnewline 
    \hline
     3.6 & 1.484790 & 40 & 68.20637 & 8.814758 & 1.3, 1.5
     \tabularnewline 
    \hline
     3.7 & 1.477769 & 40 & 68.23274 & 9.121955 & 1.1, 1.5$^\ast$
     \tabularnewline 
    \hline
     3.8 & 1.482261 & 40 & 68.41086 & 9.120428 & 1.1, 1.7
     \tabularnewline 
    \hline
     3.9 & 1.413917 & 40 & 68.56697 & 12.01028 & 1.1, 1.2
     \tabularnewline 
    \hline
     3.10 & 1.520975 & 40 & 68.83818 & 7.999162 & 1.2, 1.5
     \tabularnewline 
    \hline
     3.11 & 1.475705 & 40 & 69.00087 & 9.972652 & 1.1, 1.5$^\ast$
     \tabularnewline 
    \hline
     3.12 & 1.532000 & 40 & 69.27379 & 7.993781 & 1.4, 1.5$^\dagger$
     \tabularnewline 
    \hline
     3.13 & 1.527108 & 40 & 69.35834 & 8.274019 & 1.2, 1.6
     \tabularnewline 
    \hline
     3.14 & 1.551762 & 40 & 69.47838 & 7.407904 & 1.4, 1.5$^\dagger$
     \tabularnewline 
    \hline
     3.15 & 1.547955 & 40 & 69.68292 & 7.764712 & 1.6, 1.7
     \tabularnewline 
    \hline
     3.16 & 1.539785 & 40 & 69.78487 & 8.193480 & 1.5, 1.7
     \tabularnewline 
    \hline
     3.17 & 1.546480 & 40 & 69.85900 & 7.999795 & 1.3, 1.7
     \tabularnewline 
    \hline
     3.18 & 1.549780 & 40 & 70.16253 & 8.171339 & 1.2, 1.7
     \tabularnewline 
    \hline
     3.19 & 1.561117 & 40 & 70.74117 & 8.296474 & 1.4, 1.7
     \tabularnewline 
    \hline
  \end{tabular}
  \vspace{3mm}
  \caption{Numerical solutions obtained by the CMPA/Newton on $\Omq$ with $\partial_{xy}$ discretized using left-sided finite differences. Graphs are shown in Figs.~\ref{fig:sol_cmpa1},~\ref{fig:sol_cmpa2}.}
  \label{tab:sol_cmpa}
\end{table}

\section{Remarks on the numerics}
\subsection{Bias in the discretization of $\partial_{xy}$}
\label{sec:bias}
In this section we examine the influence of the
discretization of the mixed derivative $\partial_{xy}$ on the numerical
solution. We recall that the mixed derivative $\partial_{xy}$ can be
discretized using left/right-sided finite differences
\pref{eq:dx_RL}, \pref{eq:d_xy_RL} or using the fast
Fourier transform (Section~\ref{sec:d_xy}). For 
comparison we use the single-dimple solution on
$\Omega=(-100,100)^2$ at load $\l=1.4$ obtained by the MPA. 

Let $\Delta x=\Delta y=0.5$. Table~\ref{tab:mpa_xy} gives a list of
numerical experiments together with the values of shortening and
energy. Figure~\ref{fig:profiles_xy} shows a profile of the numerical
solutions in the circumferential direction at $x=0$.

\begin{table}[htbp]
  \centering\setlength{\unitlength}{1mm}
  \begin{tabular}{|c|c|c|c|c|c|c|}
    \hline
    domain & discretization of $\partial_{xy}$ & $\l$ & $S$ & $E$ &
    $F_\l$ & Figure~\ref{fig:profiles_xy} \tabularnewline
    \hline
    $\Omega$ or $\Omq$ & Fourier & 1.4 & 14.93529 & 24.71825 & 3.808850 &
    \begin{picture}(11.25,2)
      \linethickness{.5pt}
      \put(0,1){\line(1,0){11.25}}
    \end{picture} \tabularnewline
    \hline
    $\Omega$ & left/right-sided & 1.4 & 14.93617 & 24.70828 & 3.797636 &
    \begin{picture}(11.25,2)
      \linethickness{.1pt}
      \put(0,1){\line(1,0){11.25}}
    \end{picture} \tabularnewline
    \hline
    $\Omq$ & left-sided & 1.4 & 17.73822 & 29.42997 & 4.596460 &
    \begin{picture}(11.25,2)
      \linethickness{.15mm}
      \multiput(0,1)(1,0){12}{\line(1,0){.15}}
    \end{picture} \tabularnewline
    \hline
    $\Omq$ & right-sided & 1.4 & 12.81205 & 21.16342 & 3.226549 &
    \begin{picture}(11.25,2)
      \linethickness{.25pt}
      \multiput(0,1)(2.5,0){5}{\line(1,0){1.25}}
    \end{picture} \tabularnewline
    \hline
  \end{tabular}
  \vspace{3mm}
  \caption{Single-dimple numerical solution obtained by the MPA with and
    without the symmetry assumption \pref{eq:sym} and with various kinds
    of discretization of $\partial_{xy}$.}
  \label{tab:mpa_xy}
\end{table}

\begin{figure}[htbp]
  \centering\setlength{\unitlength}{1mm}
  \begin{picture}(120,57)
    \put(0,0){\includegraphics[width=120mm]{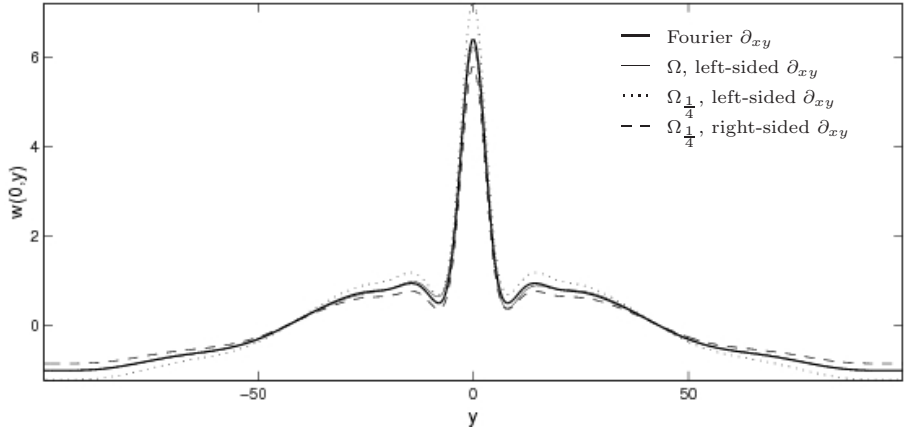}}
    \linethickness{.5pt}
    \put(82,52){\line(1,0){3.755}}
    \linethickness{.1pt}
    \put(82,48){\line(1,0){3.75}}
    \linethickness{.15mm}
    \multiput(82.3,44)(1,0){4}{\line(1,0){.15}}
    \linethickness{.25pt}
    \multiput(82,40)(2.5,0){2}{\line(1,0){1.25}}
    \put(88,51){\scriptsize Fourier $\partial_{xy}$}
    \put(88,47){\scriptsize $\Omega$, left-sided $\partial_{xy}$}
    \put(88,43){\scriptsize $\Omq$, left-sided $\partial_{xy}$}
    \put(88,39){\scriptsize $\Omq$, right-sided $\partial_{xy}$}
  \end{picture}
  \caption{Profile of the single-dimple numerical solution $w_\MP$ at
    $x=0$ obtained by the MPA with and without the symmetry assumption
    \pref{eq:sym} and with various kinds of discretization of
    $\partial_{xy}$.}
  \label{fig:profiles_xy}
\end{figure}

On the full domain $\Omega$ with no assumption on symmetry of
solutions the discretization of $\partial_{xy}$ using the
left/right-sided finite differences~\pref{eq:dx_RL} provides a
numerical solution that is slightly asymmetric
(Fig.~\ref{fig:profiles_xy}, thin solid line).  The Fourier transform
provides a symmetric solution (Fig.~\ref{fig:profiles_xy}, thick solid
line). The same numerical solution can be obtained on $\Omq$ under the
symmetry assumption \pref{eq:sym} with $\partial_{xy}$ discretized
using the fast cosine/sine transform.

On $\Omq$ the symmetry of numerical solutions is guaranteed by
assumption~\pref{eq:sym}. The use of left/right-sided discretization
of $\partial_{xy}$ does, however, have an influence on the shape of
the numerical solution, as Fig.~\ref{fig:profiles_xy} shows (the
dotted and the dashed line).

\subsection{Convergence}
\begin{figure}[htbp]
  \centering\setlength{\unitlength}{1mm}
  \begin{picture}(140,57)
    \put(0,0){\includegraphics[width=65mm]{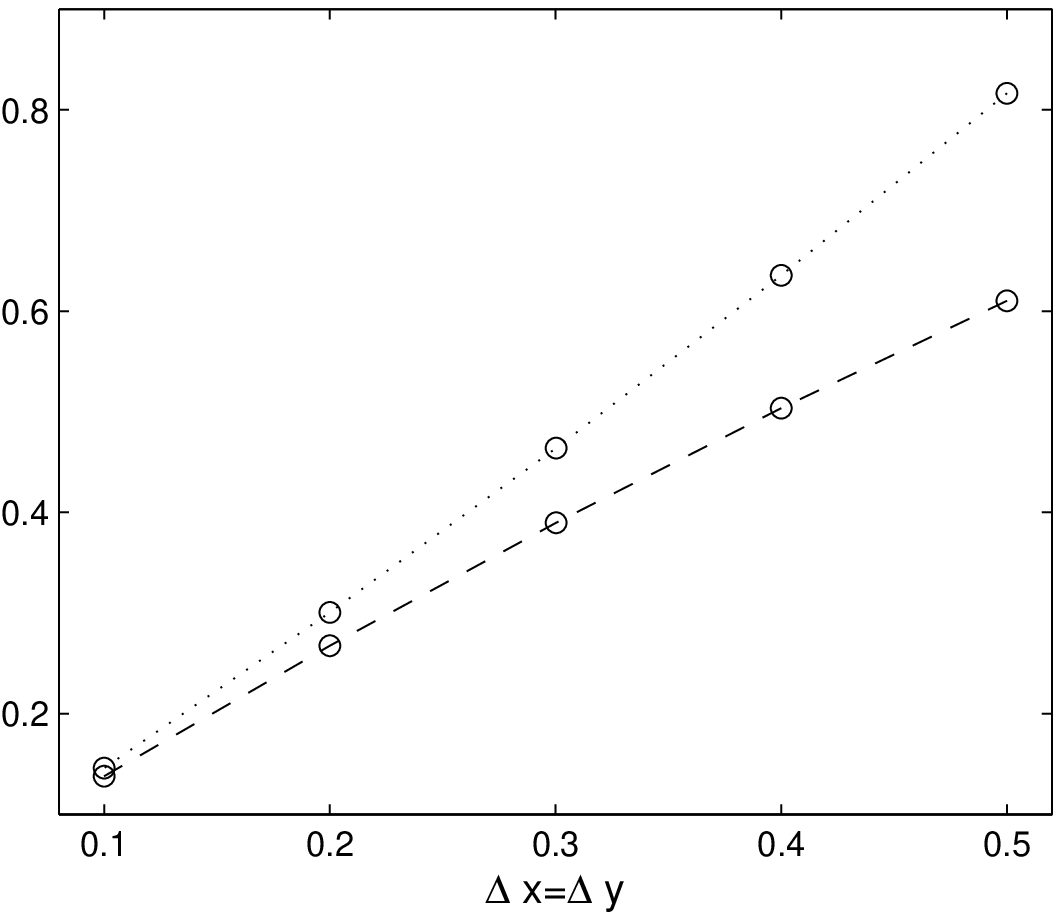}}
    \put(70,0){\includegraphics[width=65mm]{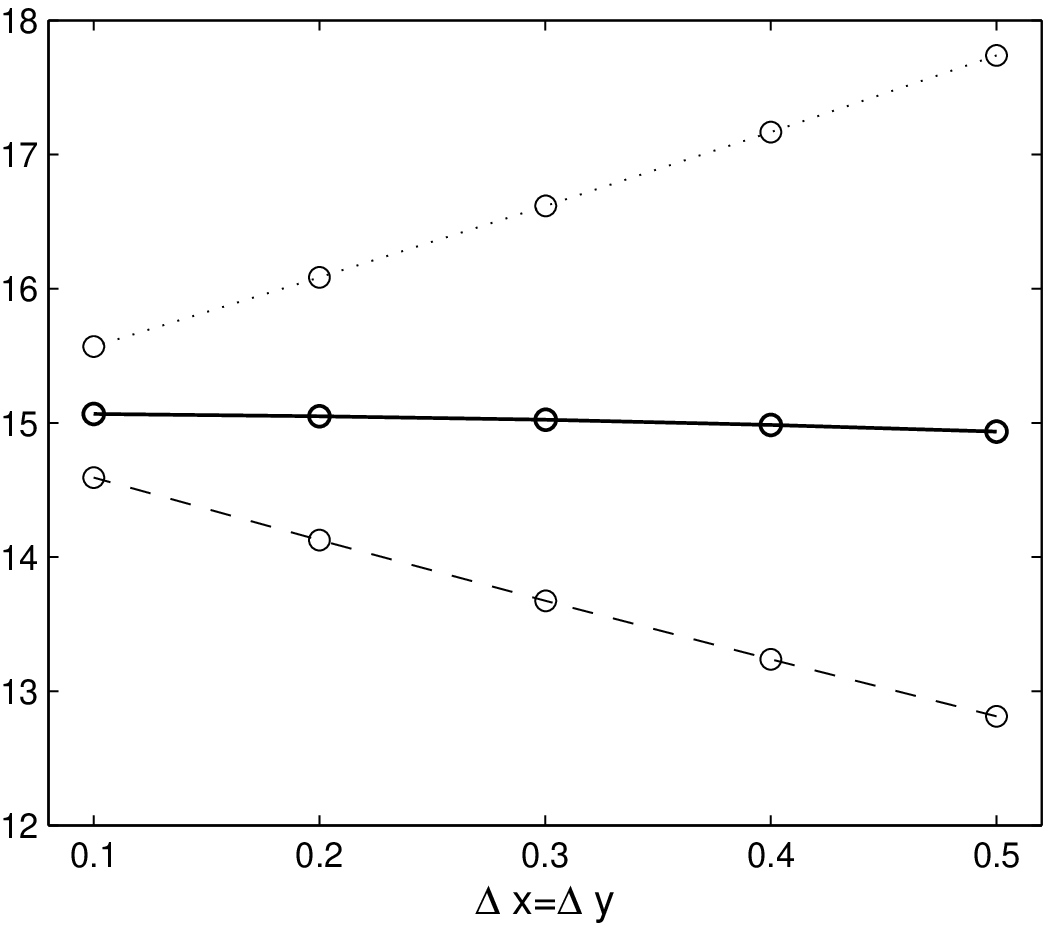}}
    \put(37,47){\scriptsize $\Mod{w^{}_L-w^{}_{CS}}_\infty$}
    \put(44,25){\scriptsize $\Mod{w^{}_R-w^{}_{CS}}_\infty$}
    \put(95,49){\scriptsize $S(w^{}_L)$}
    \put(95,15){\scriptsize $S(w^{}_R)$}
    \put(110,27){\scriptsize $S(w^{}_{CS})$}
  \end{picture}
  \caption{Influence of the size of the space step $\Delta x$, $\Delta y$
    on the numerical solution $w_\MP$ obtained by the MPA for three
    different kinds of discretization of $\partial_{xy}$. Let
    $w^{}_L$, $w^{}_R$ denote the numerical solutions obtained using
    the left and right-sided discretization of $\partial_{xy}$,
    respectively, $w^{}_{CS}$ using the fast cosine/sine transform.
    Left: comparison of the solutions in the maximum norm; right: the
    value of shortening $S$.}
  \label{fig:various_dxdy}  
\end{figure}

We now turn to the influence of the size of the space step $\Delta x$,
$\Delta y$ on the numerical solution. We run the MPA on $\Omq$ under the
symmetry assumption~\pref{eq:sym} with $\partial_{xy}$ discretized
using (a) the fast cosine/sine transform, (b) left-sided finite
differences, (c) right-sided finite differences. We consider $\Delta
x=\Delta y = 0.5, 0.4, 0.3, 0.2, 0.1$, i.e., we take 200, 250, 333,
500, and 1000 points in both axis directions, respectively.
Figure~\ref{fig:various_dxdy} illustrates convergence 
as $\Delta x, \Delta y\to
0$ of the numerical solutions obtained by various types of discretization
of $\partial_{xy}$.

\subsection{Dependence on the size of the domain}\label{sec:domain_size}
As observed in~\cite{HoLoPe1}, the localized nature of the solutions
suggests that they should be independent of domain size, in the sense
that for a sequence of domains of increasing size the solutions
converge (for instance uniformly on compact subsets). Such a
convergence would also imply convergence of the associated energy
levels. Similarly, we would expect that the aspect ratio of the domain
is of little importance in the limit of large domains.

We tested these hypotheses by computing the single-dimple solution on
domains of different sizes and aspect ratios. In all the computations
the space step $\Delta x = \Delta y = 0.5$ is fixed. In order 
to use the continuation method of Sec.~\ref{sec:continuation},
we discretized $\partial_{xy}$ using the left-sided finite differences.
We also assumed symmetry of solutions given by~\pref{eq:sym} and worked
on $\Omq$. We recall the notation of computational domains,
$\Omega=(-a,a) \times (-b,b)$, $\Omq=(-a,0)\times(-b,0)$.

\begin{figure}[htbp]
 \centering\setlength{\unitlength}{1mm}
 \begin{picture}(139,115)
   \put(0,80){\includegraphics[width=45mm]{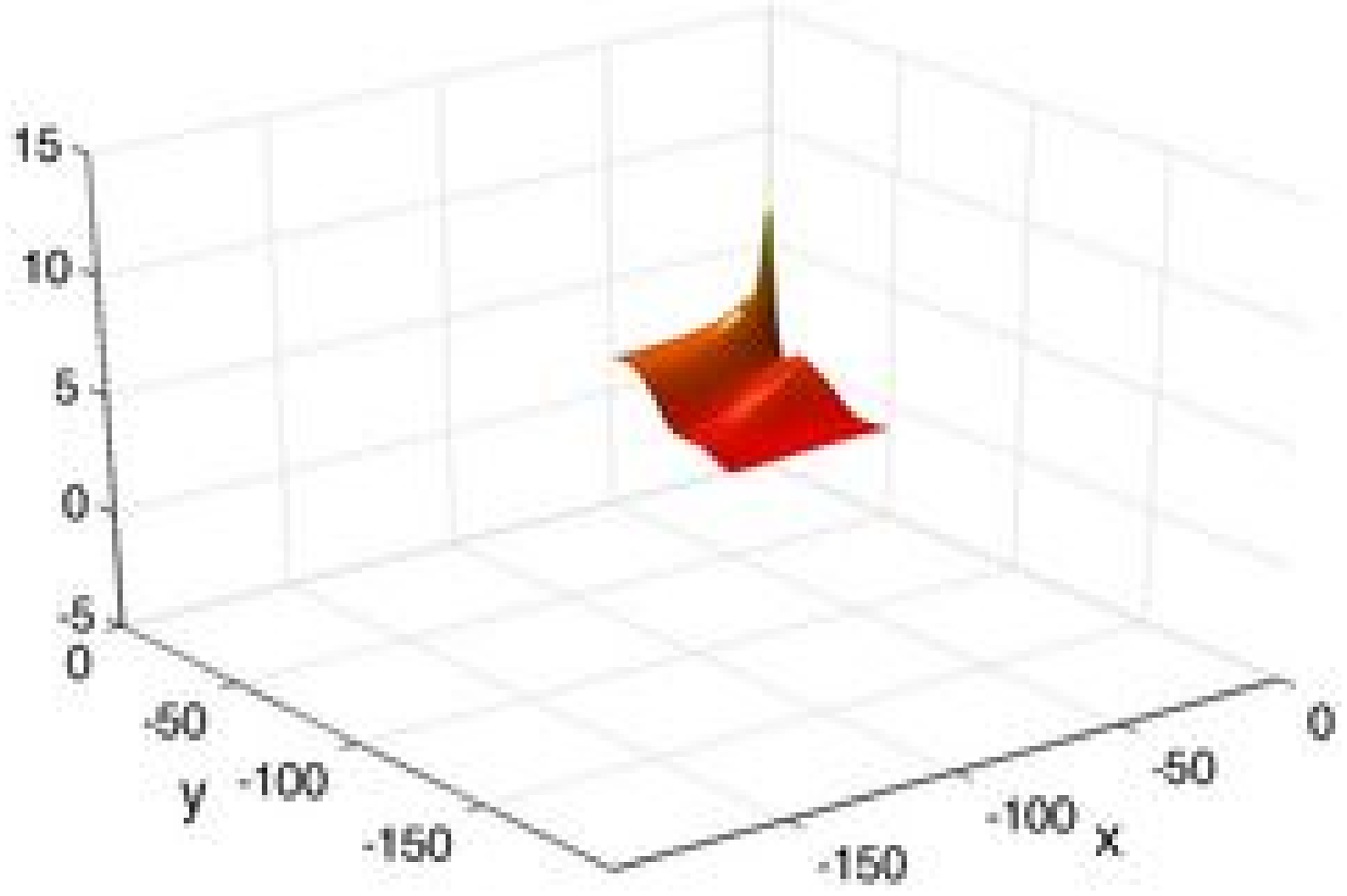}}
   \put(47,80){\includegraphics[width=45mm]{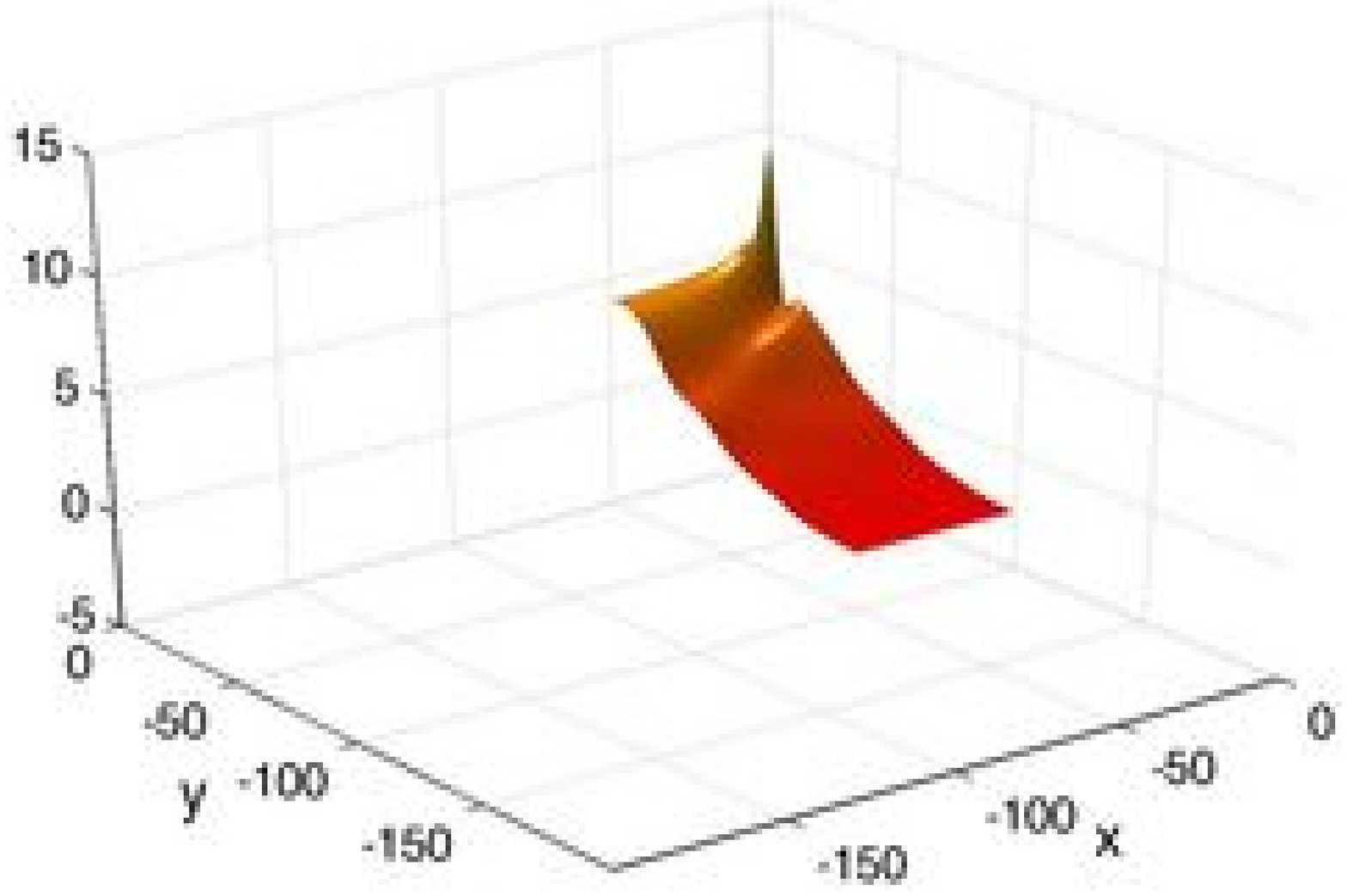}}
   \put(94,80){\includegraphics[width=45mm]{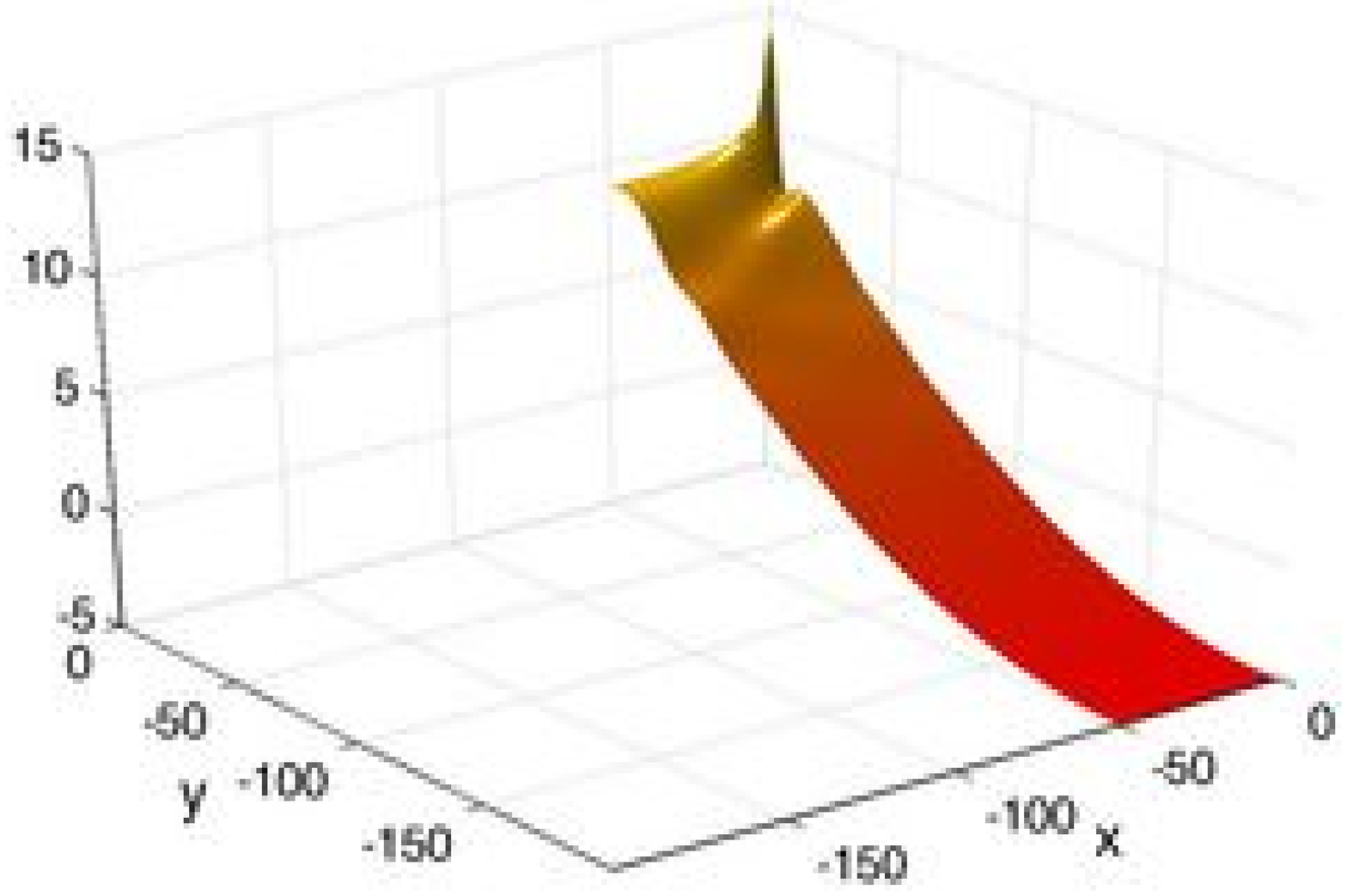}}
   \put(5,110){\scriptsize $(a,b)=(50,50)$}
   \put(52,110){\scriptsize $(a,b)=(50,100)$}
   \put(99,110){\scriptsize $(a,b)=(50,200)$}
   \put(0,40){\includegraphics[width=45mm]{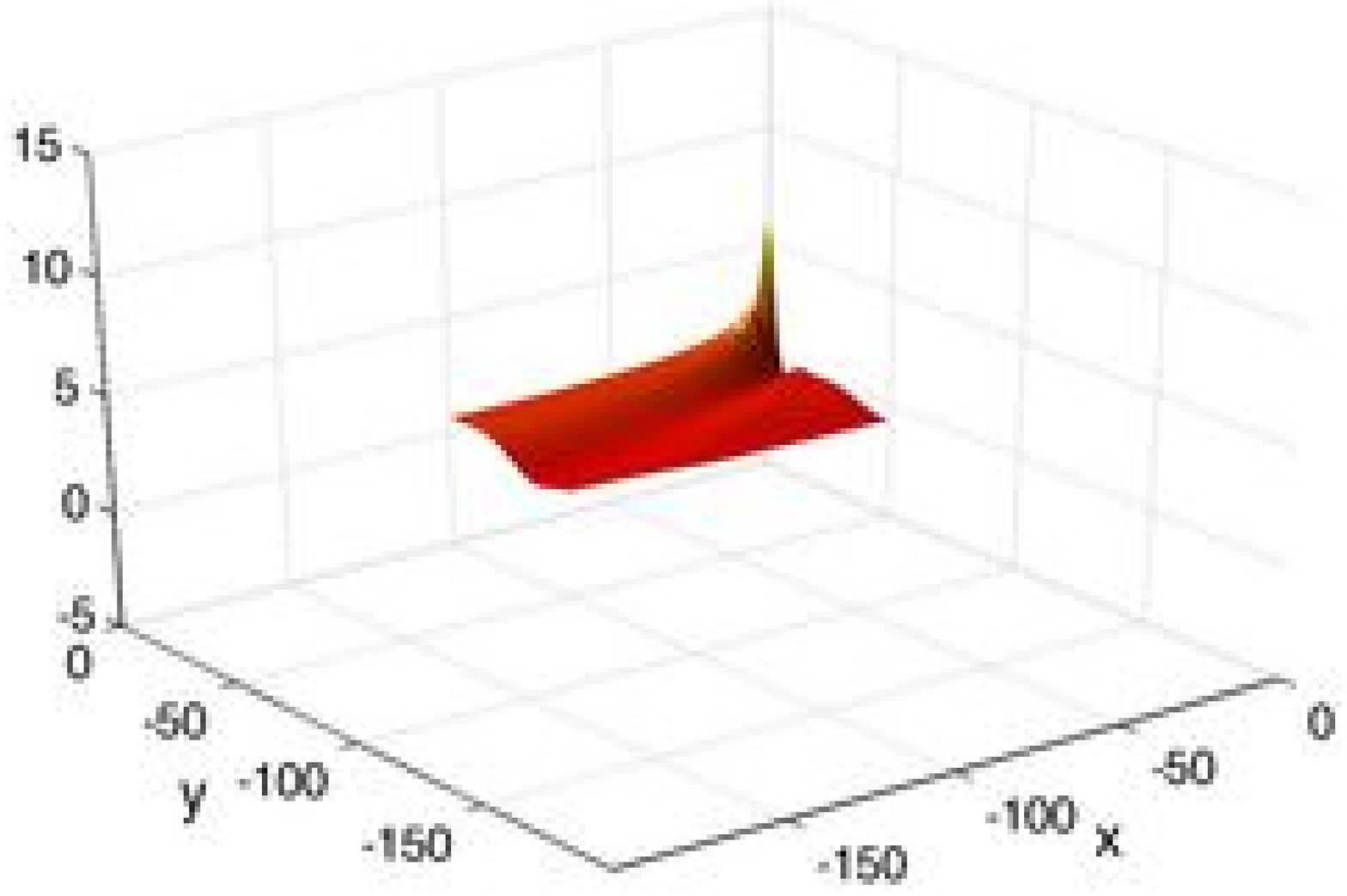}}
   \put(47,40){\includegraphics[width=45mm]{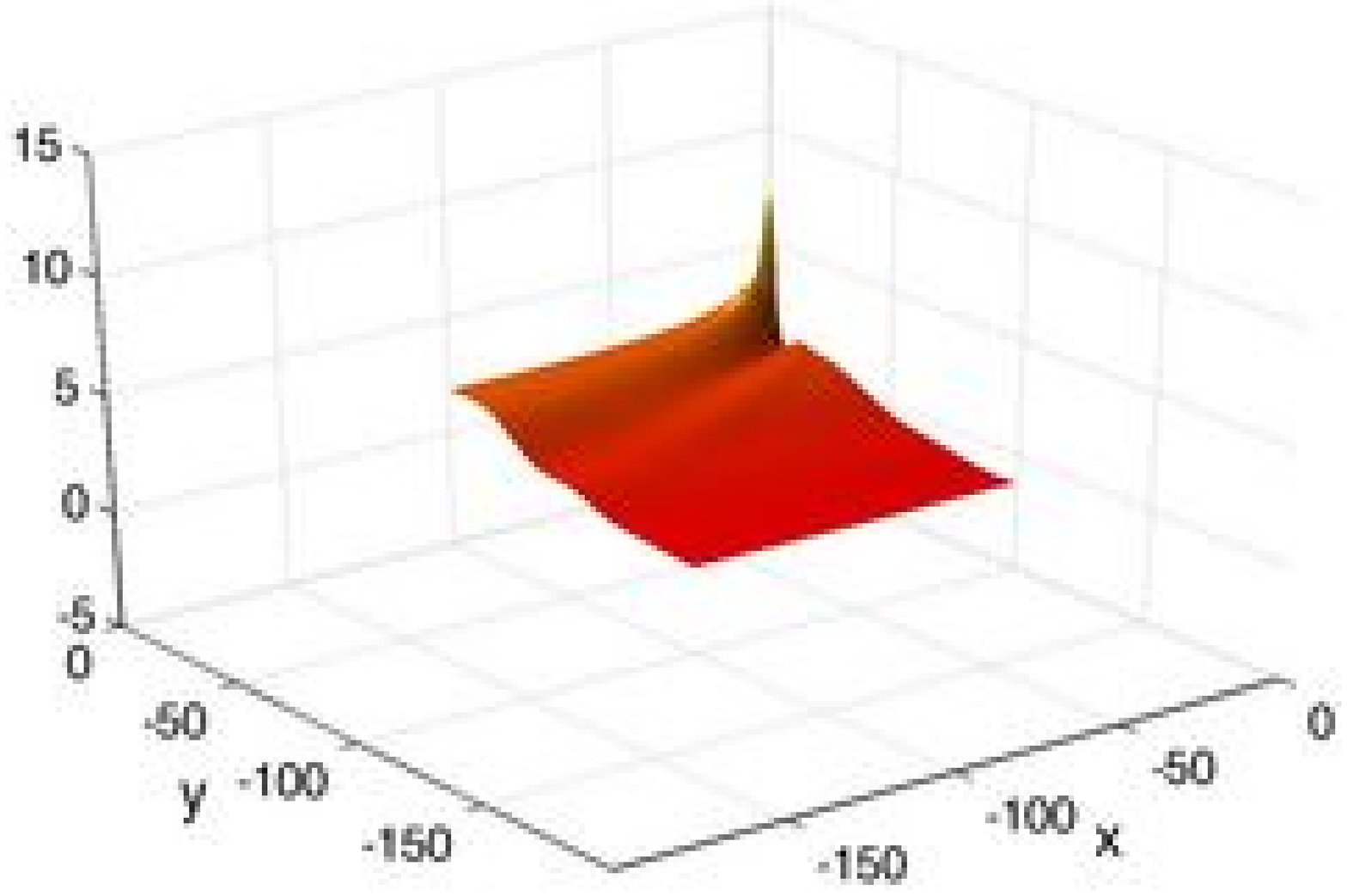}}
   \put(94,40){\includegraphics[width=45mm]{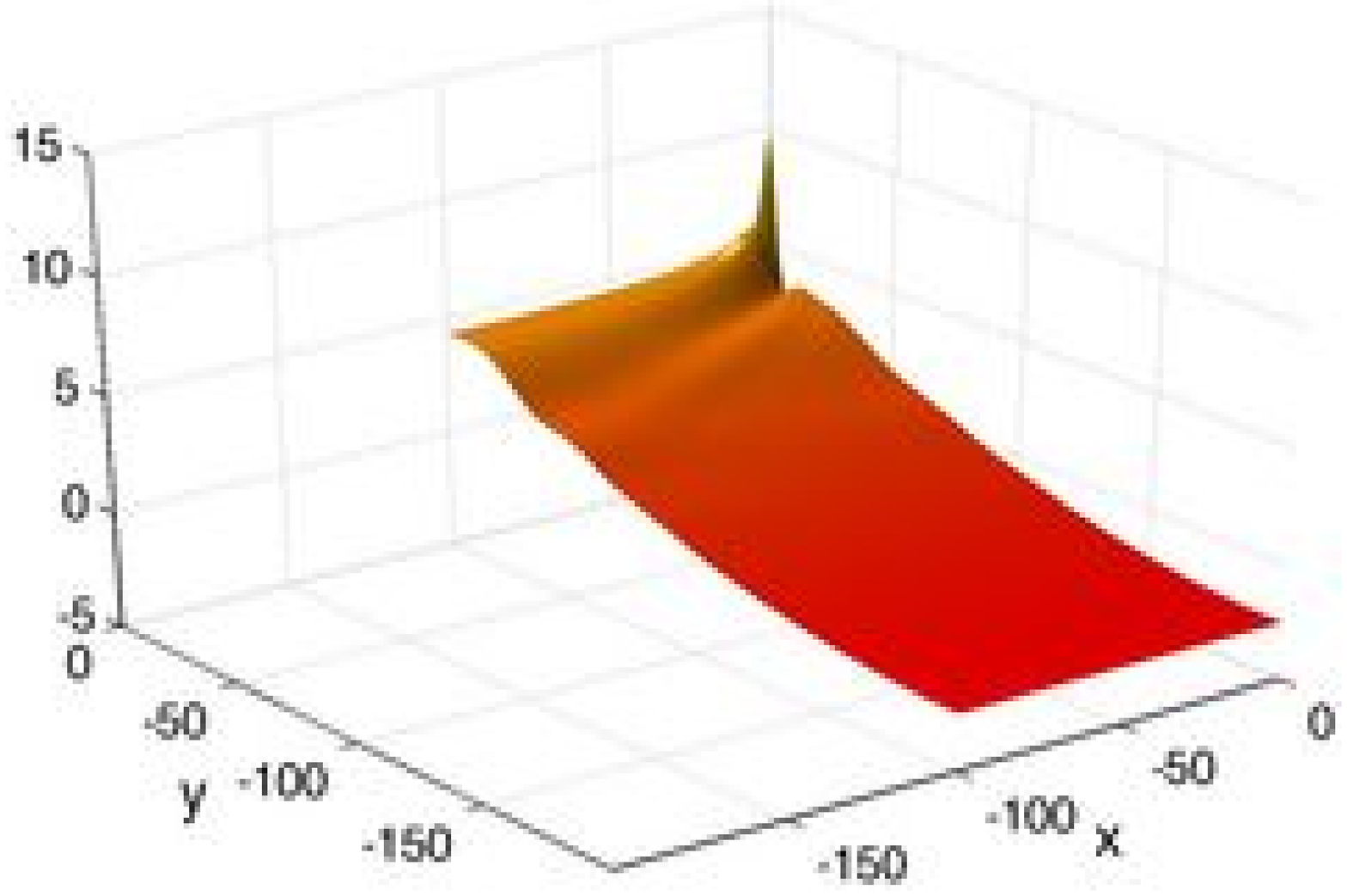}}
   \put(5,70){\scriptsize $(a,b)=(100,50)$}
   \put(52,70){\scriptsize $(a,b)=(100,100)$}
   \put(99,70){\scriptsize $(a,b)=(100,200)$}
   \put(0,0){\includegraphics[width=45mm]{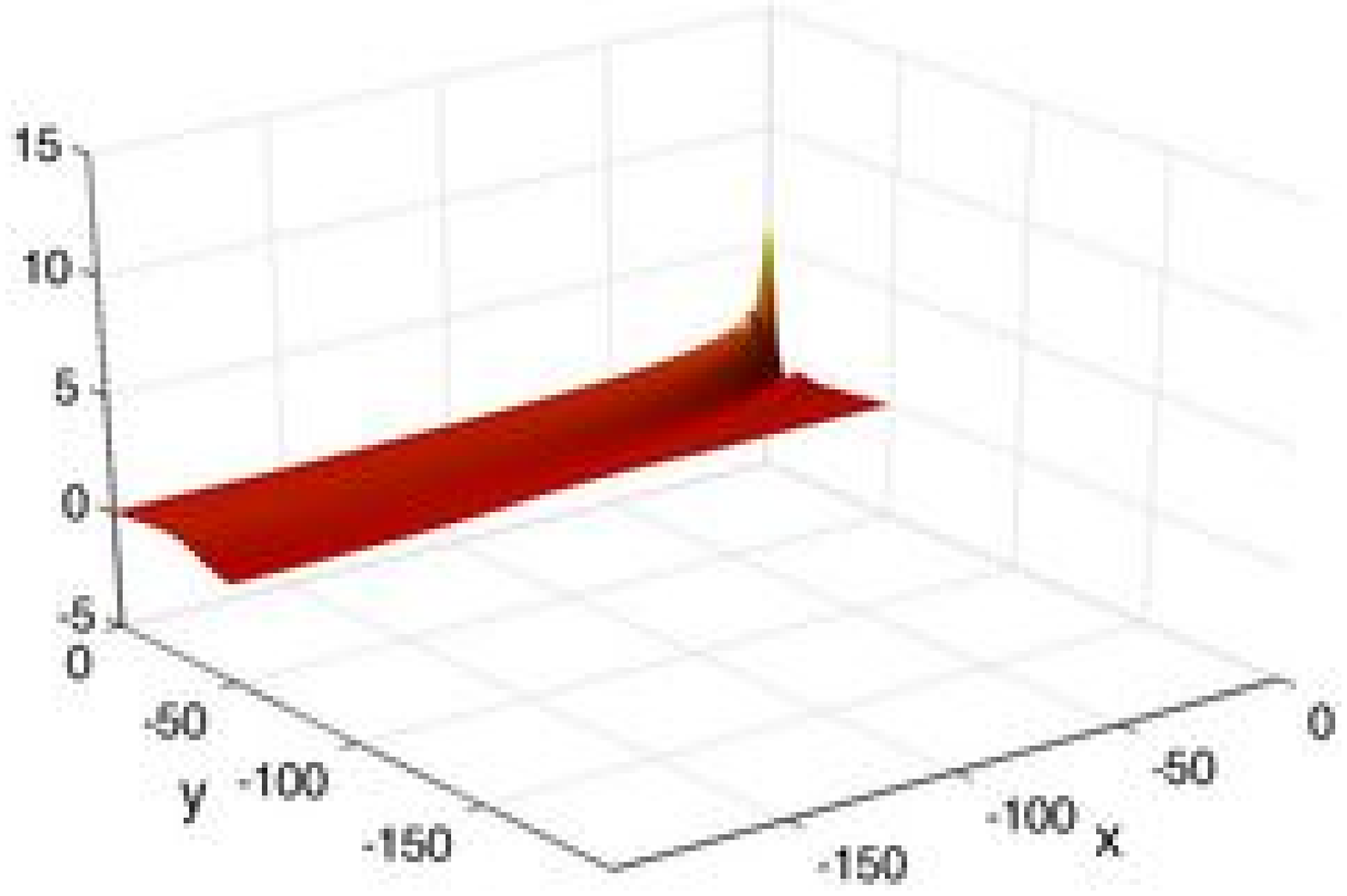}}
   \put(47,0){\includegraphics[width=45mm]{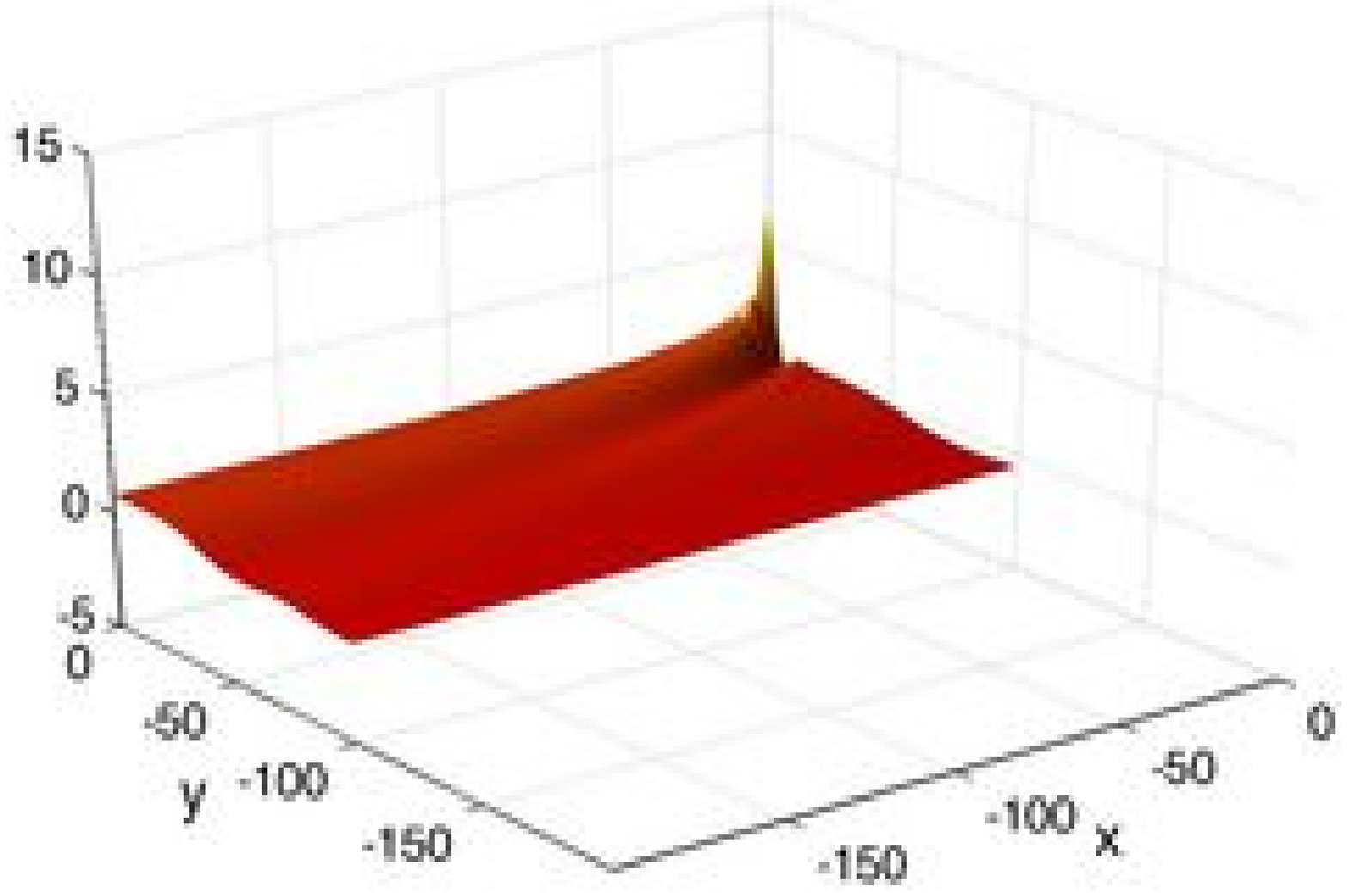}}
   \put(94,0){\includegraphics[width=45mm]{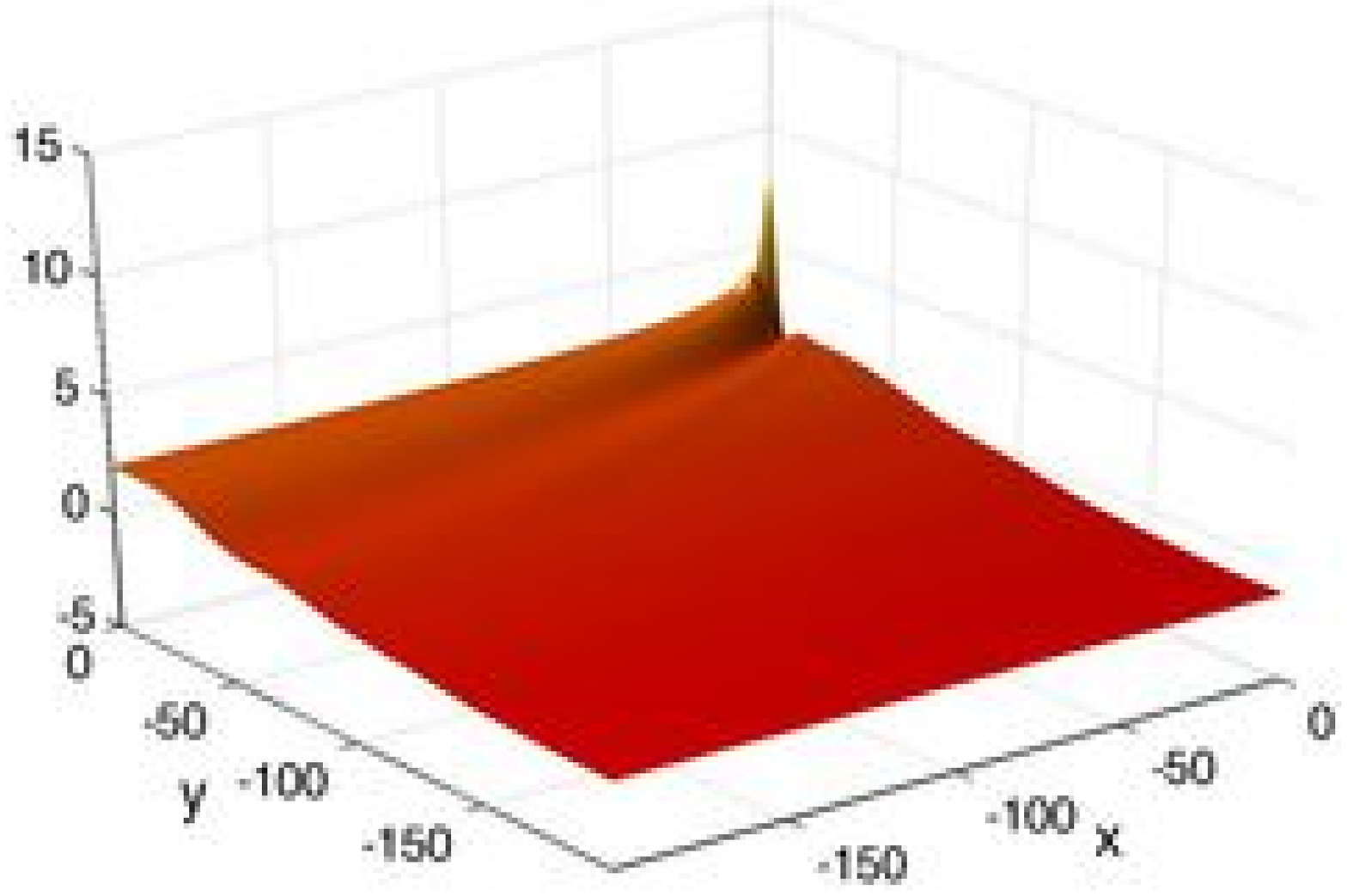}}
   \put(5,30){\scriptsize $(a,b)=(200,50)$}
   \put(52,30){\scriptsize $(a,b)=(200,100)$}
   \put(99,30){\scriptsize $(a,b)=(200,200)$}
 \end{picture}
 \caption{The single-dimple mountain-pass solution
   with $\l=1.4$ computed under the assumption of
   symmetry~\pref{eq:sym} with left-sided discretization of
   $\partial_{xy}$ for various domain sizes.
   }
 \label{fig:w_domains}
\end{figure}

Figure~\ref{fig:w_domains} shows the results for load $\l=1.4$. First
we notice that the central dimple has almost the same shape in all the
shown cases. But there seems to be a difference in the slope of the
``flat'' part leading to this dimple. On domains with small $a$ (short
cylinder) the derivative in the circumferential $y$-direction in this
part is larger than on domains with larger $a$ (longer cylinder). The
circumferential length $b$ seems to be less important for the shape of
the solution: for example, the cases (200,50) and (200,100) look like
restrictions of the case (200,200) to smaller domains.

We
take a closer look at domains of sizes $(a,b) = (100,100)$,
$(100,200)$, $(200,100)$, and $(200,200)$ and the corresponding
solutions $w_{100,100}$, $w_{100,200}$, $w_{200,100}$, and
$w_{200,200}$ shown in the figure. We compare the first three with the
last one, respectively. It does not make sense to compare the values
of $w$ itself since the energy functional $F_\l$ depends on
derivatives of $w$ only. We choose to compare $w_{xx}$ and
$w_{yy}$. Table~\ref{tab:w_100_200} gives the infinity norm of the
relative differences. Fig.~\ref{fig:w_100_200} shows graphs of the
difference $w_{100,100}-w_{200,200}$ and of the second derivatives
$(w_{200,200}-w_{100,100})_{xx}$, $(w_{200,200}-w_{100,100})_{yy}$ on
the subdomain $(-100,0)^2$.

We conclude that solutions on different domains compare well; the
maximal difference in the second derivatives of $w$ is three orders of
magnitude smaller than the supremum norm of the same derivative. We
also observe that varying the length parameter $a$ while keeping
the circumference parameter $b$ fixed causes larger changes in the numerical
solution than varying the cylinder circumference while keeping the length
fixed.

\begin{figure}[htbp]
  \centering\setlength{\unitlength}{1mm}
  \begin{picture}(139,35)
    \put(0,0){\includegraphics[width=45mm]{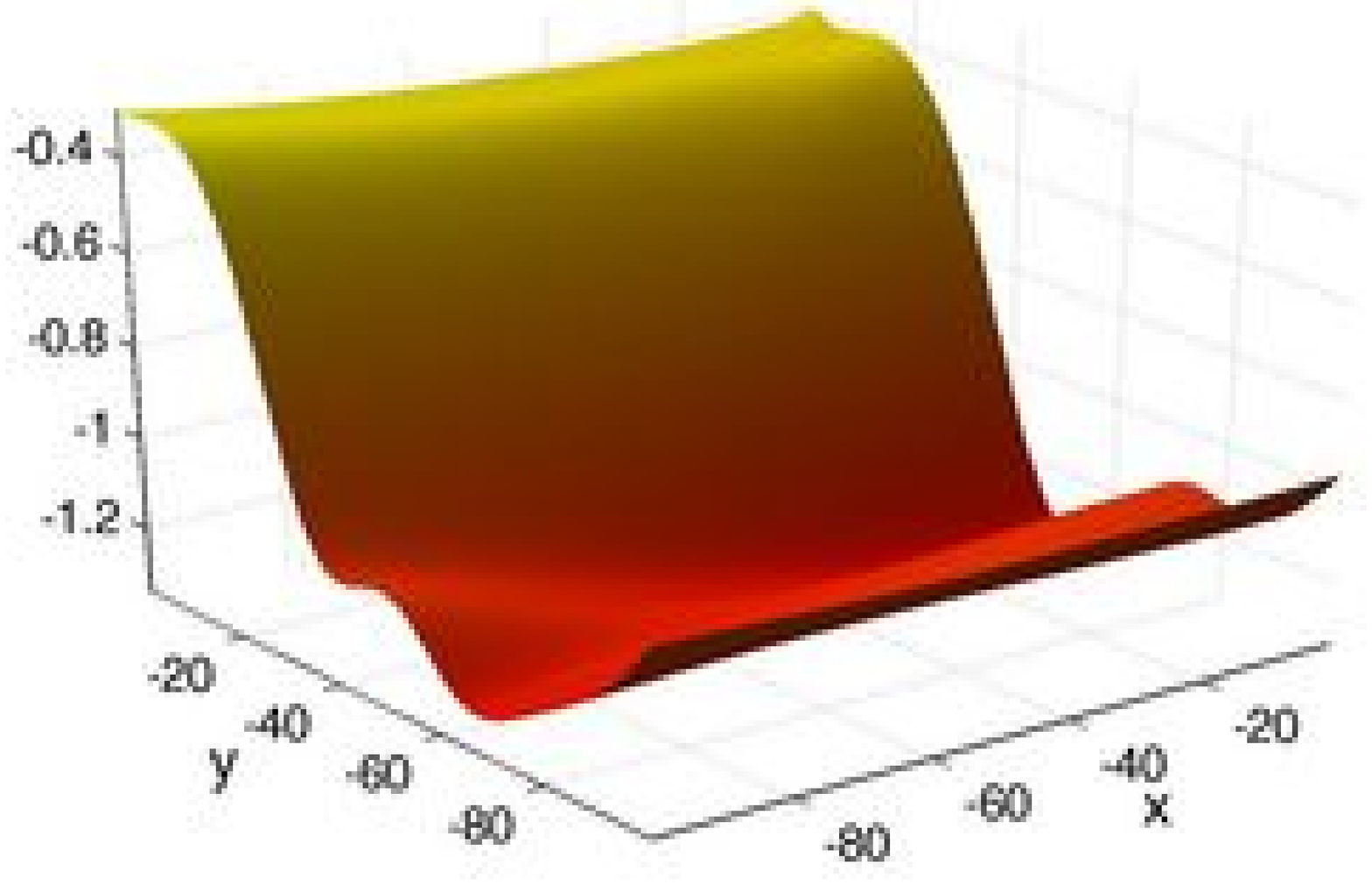}}
    \put(47,0){\includegraphics[width=45mm]{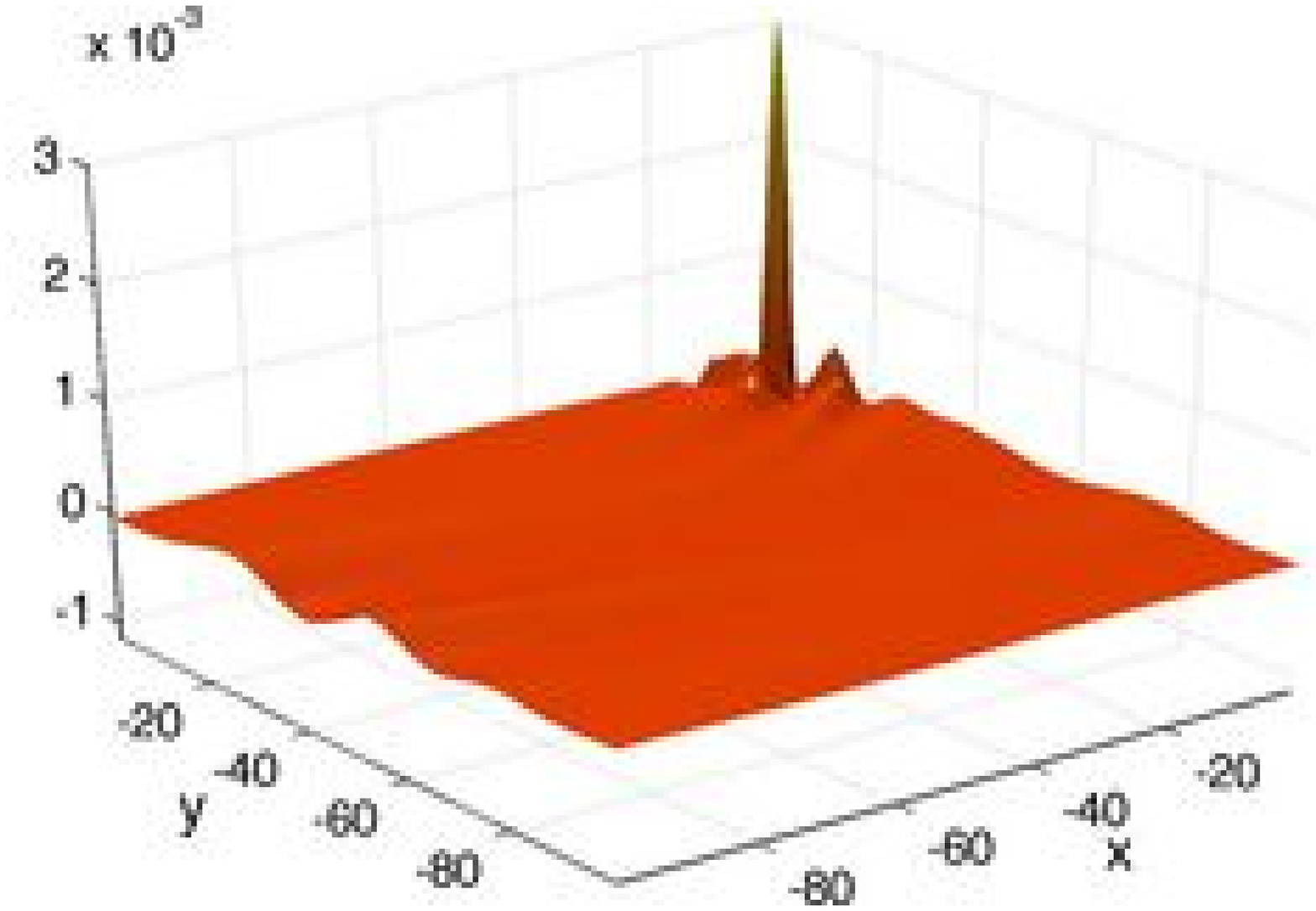}}
    \put(94,0){\includegraphics[width=45mm]{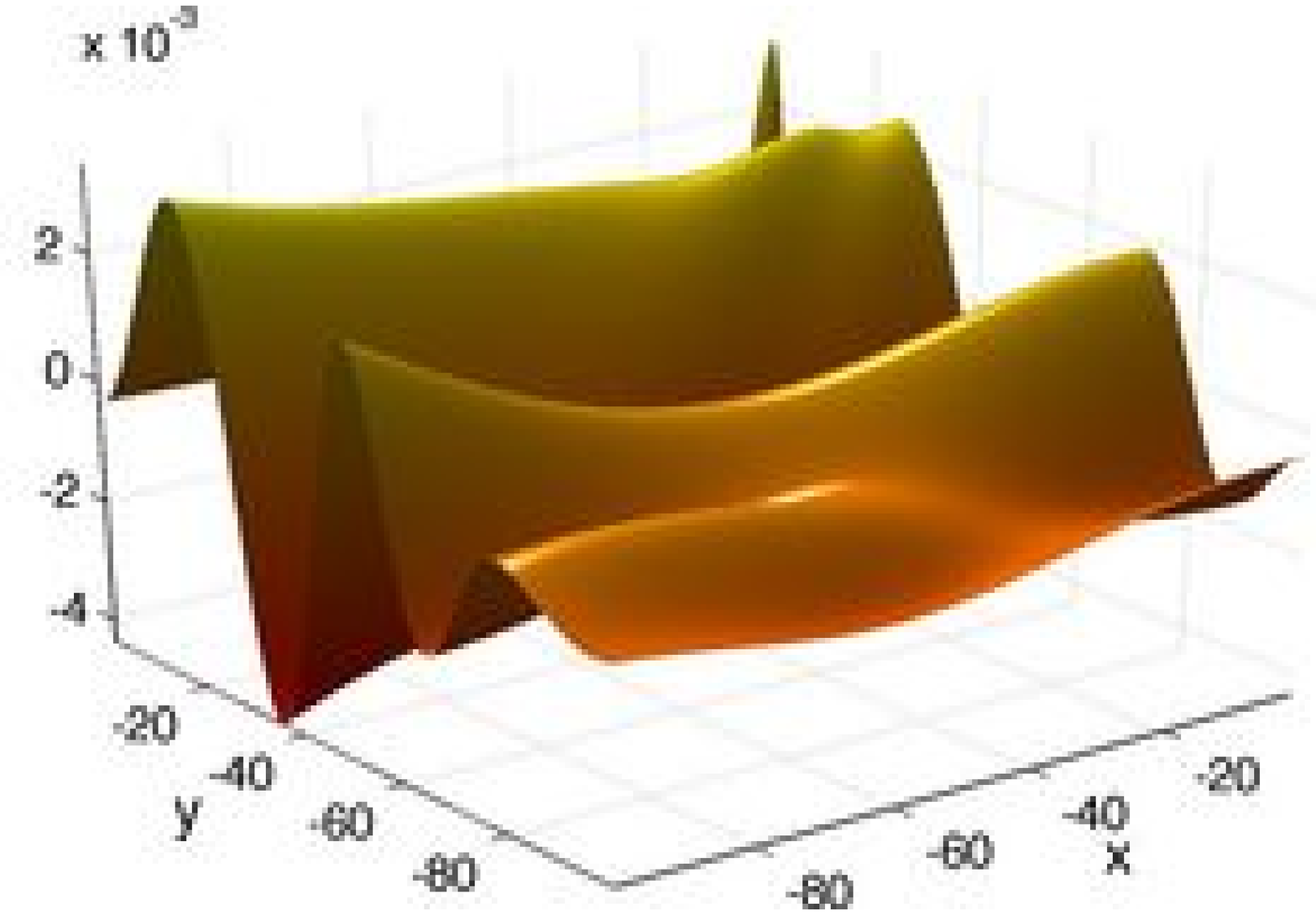}}
    \put(5,32){\scriptsize $w_{100,100}-w_{200,200}$}
    \put(52,32){\scriptsize $(w_{200,200}-w_{100,100})_{xx}$}
    \put(99,32){\scriptsize $(w_{200,200}-w_{100,100})_{yy}$}
  \end{picture}
  \caption{Comparison of solutions $w_{100,100}$, $w_{200,200}$
    from Fig.~\ref{fig:w_domains} obtained on
    square domains with $a=b=100$ and $a=b=200$, respectively, and their
    second derivatives. For a reference, we note that
    $\Mod{(w_{200,200})_{xx}}_\infty = 1.064522$ and
    $\Mod{(w_{200,200})_{yy}}_\infty = 0.8242491$.}
  \label{fig:w_100_200}
\end{figure}

\begin{table}[htbp]
  \centering
  \renewcommand{\arraystretch}{1.2}
  \begin{tabular}{c|c|c|}
    \cline{2-3}
    & $\frac{\Mod{(w-w_{200,200})_{xx}}_\infty}{\Mod{(w_{200,200})_{xx}}_\infty}$
    & $\frac{\Mod{(w-w_{200,200})_{yy}}_\infty}{\Mod{(w_{200,200})_{yy}}_\infty}$
    \tabularnewline
    \hline
    \multicolumn{1}{|c|}{$w=w_{100,100}$} & $2.835\cdot 10^{-3}$ &
    $5.313\cdot 10^{-3}$ \tabularnewline
    \hline
    \multicolumn{1}{|c|}{$w=w_{100,200}$} & $1.943\cdot 10^{-3}$ &
    $4.917\cdot 10^{-3}$ \tabularnewline
    \hline
    \multicolumn{1}{|c|}{$w=w_{200,100}$} & $1.827\cdot 10^{-4}$ &
    $9.638\cdot 10^{-4}$ \tabularnewline
    \hline
  \end{tabular}
  \renewcommand{\arraystretch}{1}
  \vspace{3mm}
  \caption{Comparison of the second derivatives of solutions
    from Fig.~\ref{fig:w_domains}
    computed on domains with $(a,b)=(100,100)$,
    $(100,200)$, $(200,100)$, and $(200,200)$.}
  \label{tab:w_100_200}
\end{table}

Another way of studying the influence of the domain size on the
numerical solution is comparing solution branches obtained by
continuation as described in
Sec.~\ref{sec:continuation}. We start with
the mountain-pass solution for $\l=1.4$ shown in
Fig.~\ref{fig:w_domains} and continue it for both $\l>1.4$ and
$\l<1.4$. The results are presented in Fig.~\ref{fig:cont}. We observe
that the branches corresponding to the considered domains do not
differ much for the range of $\l$ between approximately $0.71$ and
$2$. Below $\l\approx 0.71$ the size of the domain, particularly the
length of the cylinder described by $a$, has a strong influence. The
graph on the right shows that the larger (longer) the domain $\Omega$
the smaller the value of $\l$ at which the norm $\Mod{w}_X$ starts to
rapidly increase for decreasing $\l$. The graph on the left shows the
energy $F_\l(w)$ along a solution branch. The data shown here
correspond to the ones in the graph on the right marked by a solid
line. The dashed line in the right graph shows also some data after
the first limit point is passed.

\begin{figure}[htbp]
  \centering\setlength{\unitlength}{1mm}
  \begin{picture}(140,50)(0,-1)
    \put(0,0){\includegraphics[width=67mm]{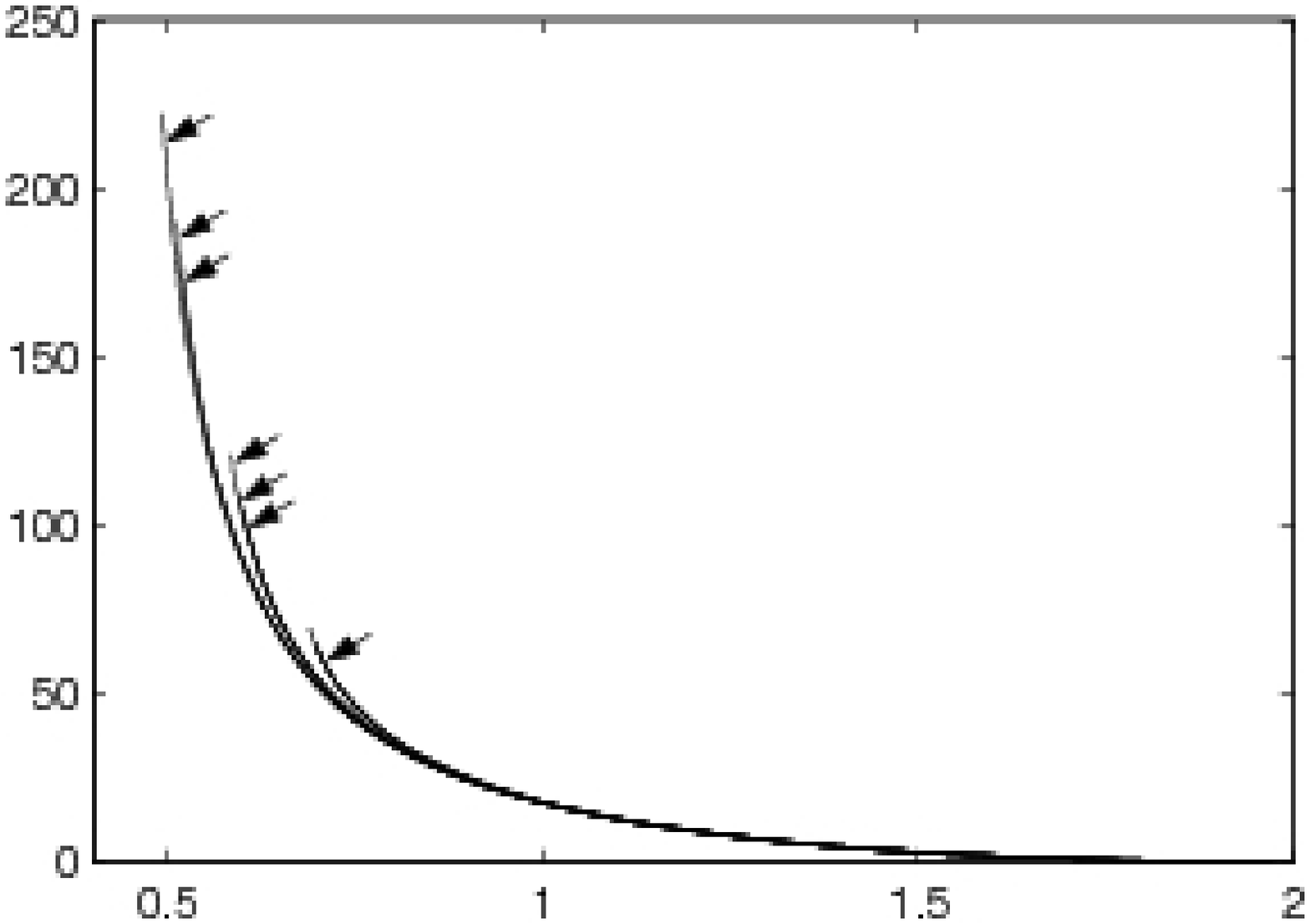}}
    \put(75,0){\includegraphics[width=65mm]{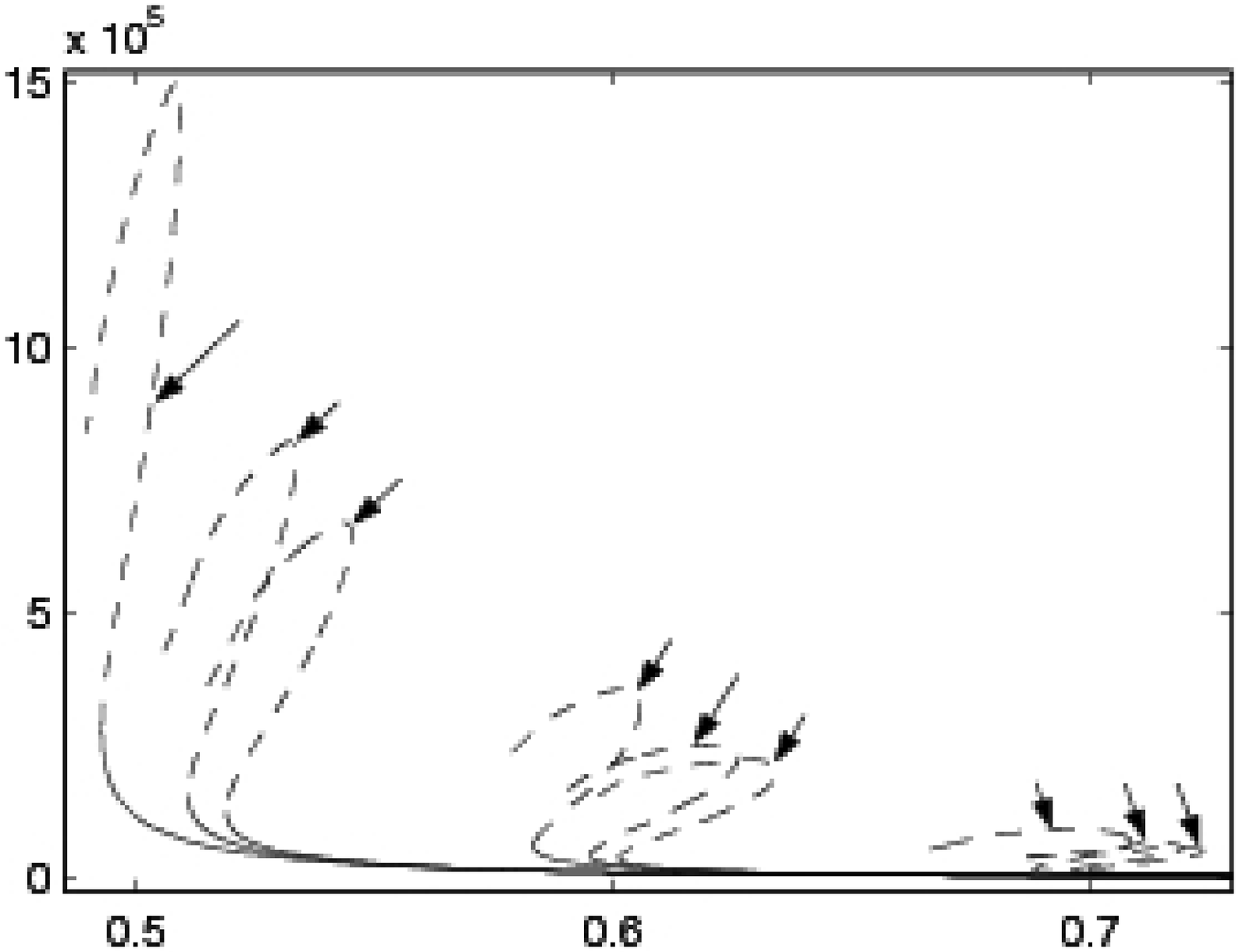}}

    \put(73,20.5){\rotatebox{90}{\scriptsize $\Mod{w}_X^2$}}
    \put(135,-1){\scriptsize $\l$}

    \put(93.5,38){\scriptsize $a=200$}
    \put(88,37){\line(1,0){15}}
    \put(88,34){\scriptsize $b=50$}
    \put(91,30){\scriptsize $b=100$}
    \put(94,26){\scriptsize $b=200$}

    \put(116,21){\scriptsize $a=100$}
    \put(110.5,20.5){\line(1,0){15}}
    \put(110.5,18){\scriptsize $b=50$}
    \put(113,16){\scriptsize $b=100$}
    \put(116.5,14){\scriptsize $b=200$}

    \put(132,12){\tiny $a=50$}
    \put(122.5,11.5){\line(1,0){16.5}}
    \put(122.5,9.5){\tiny $b=50,100,200$}

    \put(-2.5,20.5){\rotatebox{90}{\scriptsize $F_\l(w)$}}
    \put(61,-1){\scriptsize $\l$}

    \put(24.5,37){\scriptsize $a=200$}
    \put(23,33){\line(0,1){9.5}}
    \put(11.5,41){\scriptsize $b=50$}
    \put(12,36){\scriptsize $b=100$}
    \put(12.5,33){\scriptsize $b=200$}

    \put(28,23){\scriptsize $a=100$}
    \put(27,20.5){\line(0,1){6.5}}
    \put(15,25.5){\scriptsize $b=50$}
    \put(15.5,23){\scriptsize $b=100$}
    \put(16,20.5){\scriptsize $b=200$}

    \put(20,14.5){\scriptsize $b=50,100,200$, $a=50$}
  \end{picture}
  \caption{Continuation of the single-dimple solution found as a numerical
    mountain pass for $\l=1.4$ on domains of various sizes for a range
    of values $\l$. Left: $F_\l(w)$ as a function of $\l$, right:
    $\Mod{w}_{X}$ as a function of $\l$.}
  \label{fig:cont}
\end{figure}

\begin{figure}[htbp]
  \centering\setlength{\unitlength}{1mm}
  \begin{picture}(139,75)
    \put(0,40){\includegraphics[width=45mm]{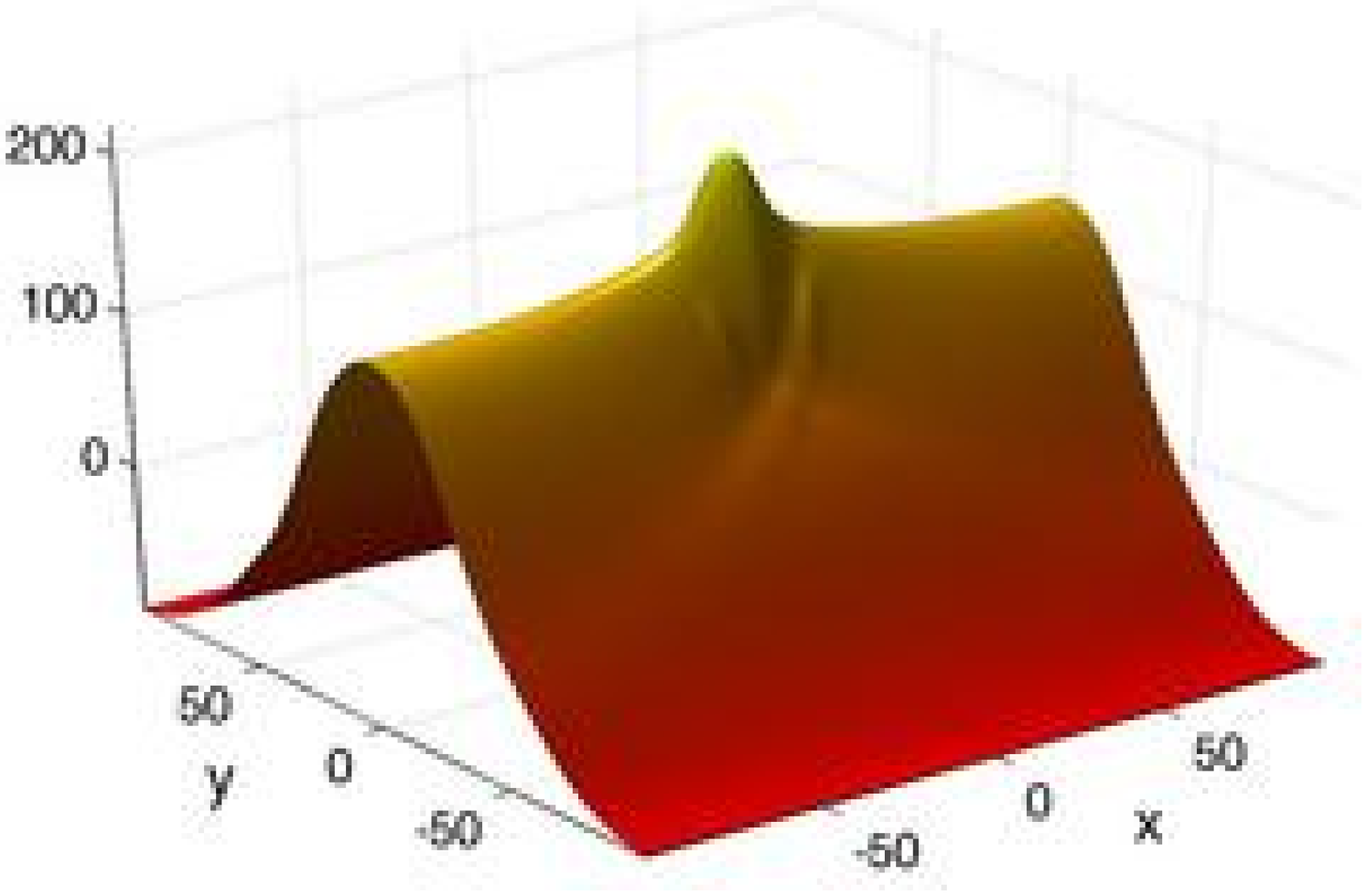}}
    \put(47,40){\includegraphics[width=45mm]{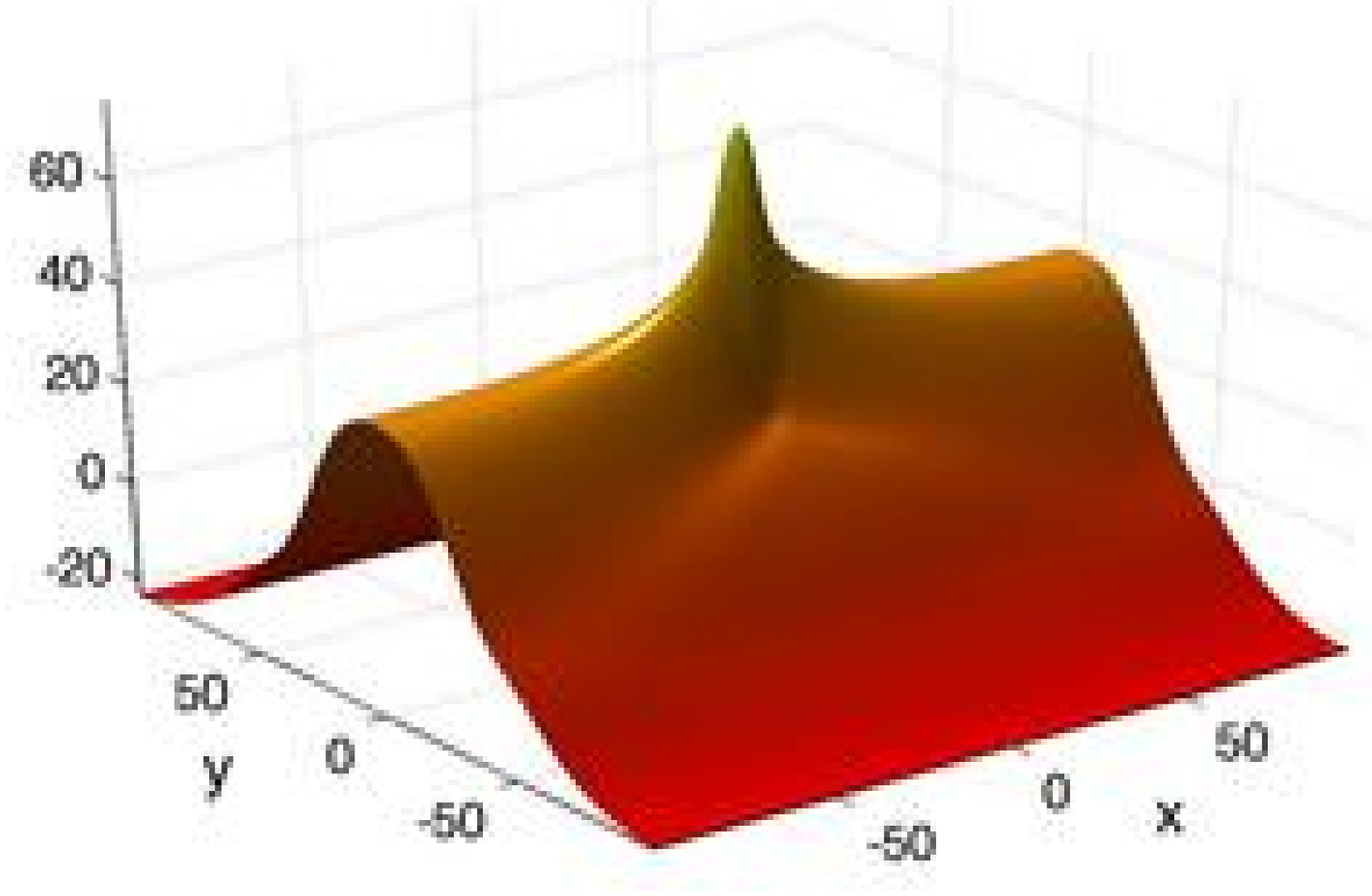}}
    \put(94,40){\includegraphics[width=45mm]{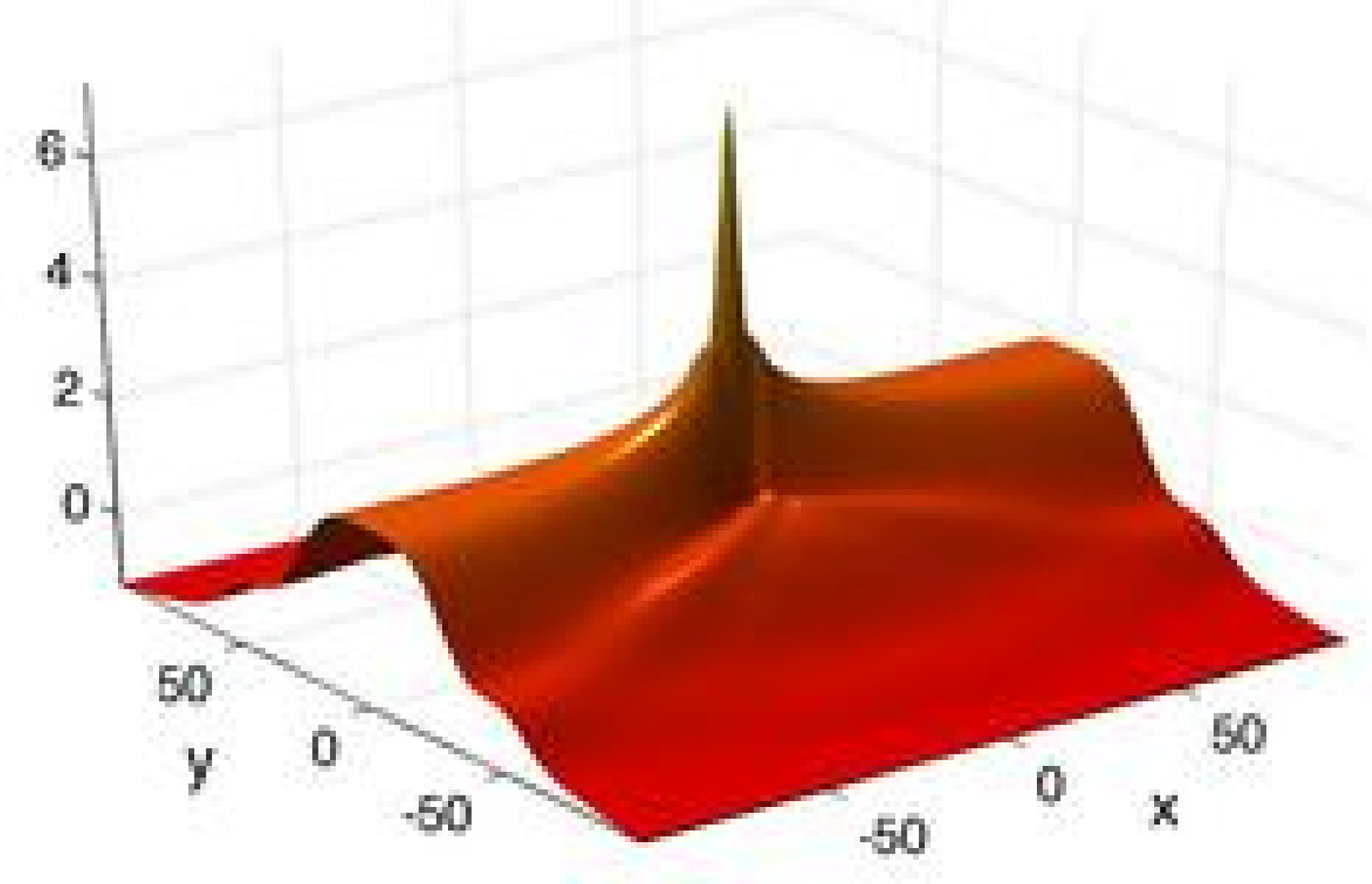}}
    \put(5,70){\scriptsize (a) $\l\approx0.593$}
    \put(52,70){\scriptsize (b) $\l=0.61$ (mountain pass)}
    \put(99,70){\scriptsize (c) $\l=1.4$ (mountain pass)}
    \put(94,0){\includegraphics[width=45mm]{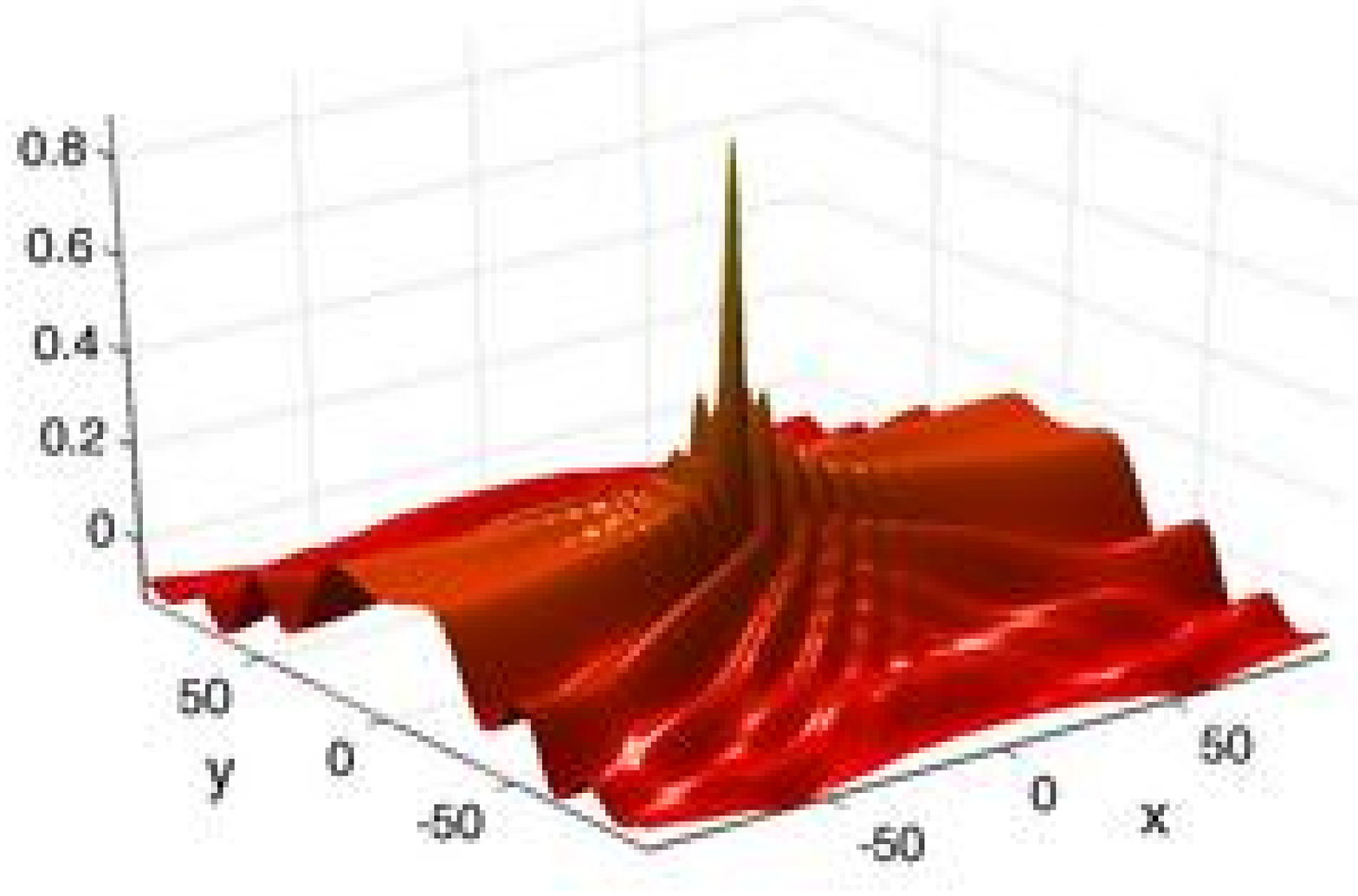}}
    \put(9.5,0){\includegraphics[width=70mm]{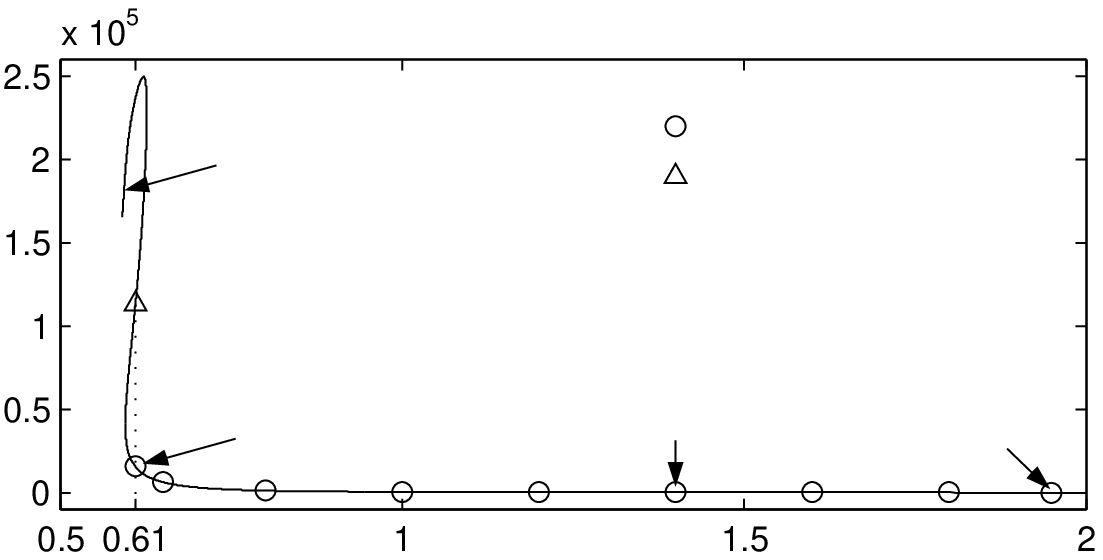}}
    \put(99,30){\scriptsize (d) $\l=1.95$ (mountain pass)}
    \put(73,-1){\scriptsize $\l$}
    \put(6,22){\rotatebox{90}{\scriptsize $\Mod{w}_X^2$}}
    \put(40,35){\scriptsize $a=b=100$}
    \put(54.5,27){\scriptsize mountain pass}
    \put(54.5,24){\scriptsize local minimizer}
    \put(54.5,21){\scriptsize (of $F_\l$)}
    \put(24,25){\scriptsize (a)}
    \put(25,8){\scriptsize (b)}
    \put(50.5,9){\scriptsize (c)}
    \put(69.2,8.3){\scriptsize (d)}
  \end{picture}
  \caption{A detailed look at the continuation of the single-dimple solution on the domain with $a=b=100$.}
  \label{fig:cont_100}
\end{figure}

Figure~\ref{fig:cont_100} shows how the graph of $w(x,y)$ changes as a
solution branch is followed. Here we chose a square domain with
$a=b=100$ and plotted the solution for four values of $\lambda$ (note
that Figs.~\ref{fig:cont_100}(c), \ref{fig:w_domains}(100,100), and
\ref{fig:profiles_xy}(dotted line) show the same numerical solution).
We observe that with decreasing $\l$ the height of the central dimple
increases, the dimple becomes wider, and the ripples (present at $\l$
close to 2) disappear. In Fig.~\ref{fig:cont_100}(a) we observe that
new dimples are being formed next to the central dimple.

It should be remarked that although we started the continuation at
$\l=1.4$ at a mountain-pass point, not all the points along a
continuation branch are mountain passes. Since it is not feasible to
use the MPA to verify this for each point, we chose just a few. Still
on the example of the domain with $a=b=100$ in
Fig.~\ref{fig:cont_100}, the circles on the continuation branch mark
those points that have also been found by the MPA (for $\l=0.61, 0.65,
0.8, 1.0, 1.2, 1.4, 1.6, 1.8, 1.95$). As described in
Sec.~\ref{sec:MPA}, in order to start the MP-algorithm a point $w_2$
is needed such that $F_\l(w_2)<0$. The analysis in~\cite{HoLoPe1}
shows that for a given $\l$ such a point exists provided the domain
$\Omega$ is large enough and in practice it is found by the SDM of
Sec.~\ref{sec:SDM}. This was, indeed, the case for all the chosen
values of $\l$ except for $\l=0.61$. In this case, starting from some
$w_0$ with a large norm, the SDM provides a trajectory $w(t)$ such
that $F_\l(w(t))>0$ for all $t>0$. In fact, the steepest descent
method converges to a local minimizer $w_\mathrm{M}^{}$ with
$F_\l(w_\mathrm{M}^{})\approx 76.1$. This is hence no mountain pass
but, nevertheless, lies on the same continuation branch and is marked
by a triangle in the figure. Despite $F_\l(w_\mathrm{M}^{})>0$ we can
still try to run the MPA with $w_2=w_\mathrm{M}^{}$. It converges and
yields $\wMP$ with $F_\l(\wMP)\approx94.8$ (marked by a circle at
$\l=0.61$ and shown in graph (b)).

The comparison of solutions computed on different domains and their
respective energies suggests that for each $\l$ we are indeed dealing
with a single, localized function defined on $\R^2$, of which our
computed solutions are finite-domain adaptations. Based on this
suggestion and the above discussion of the mountain-pass solutions we
could, for example, conclude that the mountain-pass energy
\[
V(\l,\Omega) := \inf_{w_2} \bigl\{\, F_\l\bigl(\wMP(\l,\Omega,w_2)\bigr):
F_\l(w_2)<0\, \bigr\}
\]
is a finite-domain approximation of a function $V(\l)$, whose graph
almost coincides with that of $V(\l,\Omega)$ for $\l$ not too small
(cf.~Fig.~\ref{fig:cont} left).

\section{Discussion}
\label{sec:conc}

\subsection{Variational numerical methods}
We have seen that given a complex energy surface many solutions may be 
found using
these variational techniques. For example, for a fixed end shortening 
of $S=40$,~Fig.~\ref{fig:sol_csdm}, Fig.~\ref{fig:sol_cmpa1}
and Fig.~\ref{fig:sol_cmpa2} are all solutions. Which of these solutions is of
greatest relevance depends on the question that is being asked.

In the context of the cylinder (and similar structural applications)
the mountain-pass solution from the unbuckled state ($w_1=0$) with 
minimal energy is of physical interest. Often the experimental buckling 
load may be at $20$--$30\%$ of the linear prediction from a bifucation 
analysis (in our scaling this corresponds to $\lambda=2$). 
This uncertainty in the buckling
load is a drawback for design and so ``knockdown'' factors have been
introduced, based on experimental data. It was argued in \cite{HoLoPe1}
that the energy of the mountain-pass solution $\wMP$ in fact provides a
lower bound on the energy required to buckle the cylinder and so these
solutions provide bounds on the (observed) buckling load of the
cylinder. 

This example illustrates an important aspect of the (constrained)
mountain-pass algorithm: its explicit non-locality. The algorithm
produces a saddle point which has an additional property: it is the
separating point (and level) between the basins of attraction of the
end points $w_1$ and $w_2$.

Another technique to investigate a complex energy surface is to perform 
a simulated annealing computation, essentially to solve the SDM (or the CSDM) 
problem with additive stochastic forcing. The aim in these techniques is 
often to find a global minimizer (if it exists) where there are a large 
number of local minimizers. Here by either the MPA or the CMPA we find the solution 
between two such minima and so get an estimate on the surplus energy needed to 
change between local minima.
  
\subsection{Numerical issues}

The numerical issues that we encountered are of two types. First there
are the requirements that are related to the specific problem of the
\vKD\ equations, such as the discretization of the mixed derivative
and the bracket, and the fact that the solutions are symmetric and
highly localized.

For other difficulties it is less clear. For smaller values of $\l$
each of the variational methods converged remarkably slowly. Newton's
method provides a way of improving the convergence, but the question
is relevant whether this slow convergence is typical for a whole class
of variational problems. It would be interesting to connect the rate
of convergence of, for instance, the SDM to certain easily measurable
features of the energy landscape.

\bibliography{refs2}

\end{document}